\documentclass{amsart}
\usepackage{amsmath,amssymb,amscd,latexsym}
\usepackage[all]{xy}

\newcommand{\C}{{\mathbb{C}}}

\newcommand{\N}{{\mathbb{N}}}

\newcommand{\R}{{\mathbb{R}}}

\newcommand{\Z}{{\mathbb{Z}}}

\newcommand{\Bh}{{\mathcal B}}
\newcommand{\Ch}{{\mathcal C}}

\newcommand{\Eh}{{\mathcal E}}
\newcommand{\Fh}{{\mathcal F}}
\newcommand{\Gh}{{\mathcal G}}

\newcommand{\Oh}{{\mathcal O}}

\newcommand{\Rh}{{\mathcal R}}

\newcommand{\Uh}{{\mathcal U}}

\newcommand{\Zh}{{\mathcal Z}}

\newcommand{\be}{\mathbf{1}}

\newcommand{\dr}{\mathrm{dr}\,}

\newcommand{\halb}{\frac{1}{2}}
\newcommand{\her}{\mathrm{her}}

\newcommand{\id}{\mathrm{id}}

\newcommand{\ord}{\mathrm{ord}\,}

\newcommand{\tr}{\mathrm{tr}}

\newcounter{number}[section]

\newenvironment{nummer}{\refstepcounter{number}{\noindent\arabic{section}.\arabic{number}}}{}

\newcommand{\bn}{\noindent \begin{nummer} \rm}
\newcommand{\en}{\end{nummer}}

\newenvironment{ntheorem}{\noindent {\sc Theorem:} \it}{}
\newenvironment{nlemma}{\noindent {\sc Lemma:} \it}{}
\newenvironment{nprop}{\noindent {\sc Proposition:} \it}{}
\newenvironment{ndefn}{\noindent {\sc Definition:} \it}{}
\newenvironment{ncor}{\noindent {\sc Corollary:} \it}{}
\newenvironment{nconj}{\noindent {\sc Conjecture:} \it}{}
\newenvironment{nremark}{\noindent {\sc Remark:}}{}

\newenvironment{nexamples}{\noindent {\sc Examples:} }{}

\newenvironment{nnotation}{\noindent {\sc Notation:} }{}

\newenvironment{nproof}{\noindent {\sc Proof:}}{\mbox{}\hfill 
\rule[-.2ex]{.25em}{1.8ex}}

\begin{document}

\title[Decomposition rank and $\mathcal{Z}$-stability]{{\sc Decomposition rank and $\mathcal{Z}$-stability}}

\author{Wilhelm Winter}
\address{School of Mathematical Sciences\\
University of Nottingham\\
Nottingham\\
United Kingdom}

\email{wilhelm.winter@nottingham.ac.uk}

\date{\today}
\thanks{Partially supported by EPSRC First Grant EP/G014019/1}

\setcounter{section}{-1}

\begin{abstract}
We show that separable, simple, nonelementary, unital $\mathrm{C}^{*}$-algebras with finite decomposition rank absorb the Jiang--Su algebra $\Zh$ tensorially. This has a number of consequences for Elliott's program to classify nuclear $\mathrm{C}^{*}$-algebras by their $\mathrm{K}$-theory data. In particular, it completes the classification of $\mathrm{C}^{*}$-algebras associated to uniquely ergodic, smooth, minimal dynamical systems  by their ordered $\mathrm{K}$-groups. 
\end{abstract}

\maketitle

\section{Introduction}

\noindent
The theory of $\mathrm{C}^{*}$-algebras is often thought of as noncommutative topology. This point of view was first suggested by the Gelfand--Naimark Theorem (which characterizes abelian $\mathrm{C}^{*}$-algebras as algebras of continuous functions on locally compact spaces) and has since then been a constant inspiration for many new developments, such as the Dixmier--Douady classification, Brown--Douglas--Fillmore theory and (bivariant) $\mathrm{K}$-theory, as well as for many important applications, like those related to dynamical systems. 

G.~Elliott was the first to suggest that large classes of $\mathrm{C}^{*}$-algebras might be completely classifiable in terms of $\mathrm{K}$-theoretic data. More precisely, he conjectured that separable nuclear $\mathrm{C}^{*}$-algebras are classified up to isomorphism by their so-called Elliott invariants; cf.\ \cite{Ell:classprob}. (A $\mathrm{C}^{*}$-algebra is nuclear if and only if it can be approximated by finite-dimensional ones in a suitable sense.) 

While we know today that Elliott's conjecture in its original form does not hold, there are a great number of partial verifications. Much more important, the conjecture has spurred many new developments in $\mathrm{C}^{*}$-algebra theory and its applications, revealing deep insights which go far beyond the realm of nuclear $\mathrm{C}^{*}$-algebras; cf.\ \cite{Ror:encyc} for an overview. 

Since the very beginnings of the classification program, various notions of noncommutative covering dimension have proven to be crucial for the theory; somewhat later the importance of strongly self-absorbing $\mathrm{C}^{*}$-algebras (or, more generally, of $\mathcal{D}$-stability for various strongly self-absorbing $\mathcal{D}$) was discovered. Kirchberg was the first to relate $\mathcal{D}$-stability to a regularity property of the Cuntz semigroup (and hence, in the broadest sense, to a $\mathrm{K}$-theoretic condition); cf.\ \cite{Kir:ICM} and \cite{Kir:CentralSequences}. Conditions of this type play a role in positive classification results as well as for the recently discovered counterexamples to Elliott's conjecture; cf.\ the articles \cite{KirPhi:classI}, \cite{KirPhi:classII}, \cite{BlanchardKirchberg:piHausdorff}, \cite{Toms:classproblem},  \cite{Win:Z-class}, \cite{Win:localizingEC}, \cite{HirRorWin:D-stable}, \cite{DadWin:trivial-fields} and  \cite{DadHirTomsWinter:example}, to mention but a few. 

As a first attempt to formalize the preceding observations, let us consider the following regularity conditions on a $\mathrm{C}^{*}$-algebra $A$, which, at first glance, do not seem to have much in common, but are satisfied for many of our stock-in-trade examples; even more surprisingly, there are large natural classes of $\mathrm{C}^{*}$-algebras for which these conditions are equivalent:

\begin{itemize}
\item[(A)] $A$ is topologically finite-dimensional.
\item[(B)] $A$ absorbs a suitable strongly self-absorbing $\mathrm{C}^{*}$-algebra tensorially.
\item[($\Gamma$)] $A$ allows comparison of its positive elements in the sense of Murray and von Neumann.
\item[($\Delta$)] The natural order structures on suitable homological invariants of $A$ are complete in the sense that they are sufficiently unperforated.
\end{itemize} 

Obviously, these conditions are of a somewhat philosophical nature, and require interpretation. All of them may be viewed as regularity properties, with (A) of a topological nature and (B) and ($\Gamma$) of a ($\mathrm{C}^{*}$-)algebraic type, thus approaching the homological condition ($\Delta$) from quite different directions.

There are a variety of noncommutative analogues of topological covering dimension, such as stable or real rank, dimension as an approximately (sub)homogeneous algebra, and decomposition rank (see Definition~\ref{d-dr} for the latter;  cf.\ also \cite{Rfl:sr}, \cite{BroPed:realrank}, \cite{KirWinter:dr} and \cite{Winter:subhomdr}). This results from various characterizations of dimension in the commutative case; these notions tend to agree for compact metrizable spaces, but exhibit rather different behaviour in the noncommutative setting. Often it is useful to combine various of these notions to  interpret (A), but in this paper we shall exclusively deal with the decomposition rank. This notion behaves more like a topological invariant than, say, the real or the stable rank, yet it is more flexible (being a \emph{local} concept in the sense that it does not depend on any particular increasing sequence of subalgebras) than dimension growth as an approximately (sub)homogeneous algebra.

The known strongly self-absorbing $\mathrm{C}^{*}$-algebras (cf.\ \cite{TomsWin:ssa}) form a certain hierarchy with the Jiang--Su algebra $\mathcal{Z}$ (cf.\ \ref{Z-intro}) at the bottom and the Cuntz algebra $\mathcal{O}_{2}$ at the top. The technical benefit of asking for $\mathcal{D}$-stability (with $\mathcal{D}$ one of these algebras) is to regularize $A$ in a manner that provides enough space for certain standard operations. Of course the bigger $\mathcal{D}$ is, the more space is available in $A \otimes \mathcal{D}$; on the other hand, $\mathcal{D}$-stability will then be a more restrictive condition. Accordingly, we are particularly interested in $\mathcal{Z}$-stability. 

Comparison theory for $\mathrm{C}^{*}$-algebras as introduced by Cuntz is largely modeled after that for von Neumann algebras, although one has to face subtle additional problems; cf.\ \cite{Cuntz:dimension}. Generally, positive elements will be compared up to Cuntz equivalence in terms of the values of dimension functions on their range projections; often it suffices to consider dimension functions induced by tracial states (see \cite{Ror:Z-absorbing} and \cite{Ror:UHFII}). 

Many invariants of $A$ which are homological in one sense or another, are built from equivalence relations on the cone $A_{+}$ of positive elements. Then, the natural order structure on $A_{+}$ often descends to an order structure on the invariant. In that case it turns out that comparison of positive elements corresponds to lack of perforation at the level of the (ordered) invariant; cf.\ \cite[Definition~3.3.2]{Ror:encyc} and \cite{Ror:Z-absorbing}.     

We should also remark that, whatever our concrete interpretation of (A)--($\Delta$) looks like, we will usually have to impose some additional requirements on $A$ to relate these conditions to each other. For example, our preferred notion of topological dimension might only make sense for nuclear $\mathrm{C}^{*}$-algebras,  and the relation between ($\Gamma$) and ($\Delta$) might rely on the subtle interplay between tracial states and dimension functions, which works particularly well for exact $\mathrm{C}^{*}$-algebras. 

Moreover, $\mathcal{D}$-stability  (at the least) implies nonexistence of nontrivial finite dimensional quotients. In fact, at its current stage the theory works best for $\mathrm{C}^{*}$-algebras that are sufficiently noncommutative, and attention will often be restricted to  simple, nonelementary $\mathrm{C}^{*}$-algebras (recall that a $\mathrm{C}^{*}$-algebra is elementary if it is isomorphic to the algebra of compact operators on some Hilbert space). In a sense,  we are working at the opposite end of the scale from Dixmier's and Douady's famous classification of continuous trace $\mathrm{C}^{*}$-algebras.

As a concrete interpretation of conditions (A)--($\Delta$), consider the following conjecture, variations and special cases (in many disguises) of which have been studied in abundance in the literature. In the form below, it was probably first  suggested by A.\ Toms; it was also the driving force behind earlier work of Toms and the author; cf.\ \cite[Remark~3.5]{TomsWinter:VI}. 

\bn
\label{finite-con}
\begin{nconj}
For a separable, finite, nonelementary, simple, unital and nuclear $\mathrm{C}^{*}$-algebra $A$, the following properties are equivalent:
\begin{enumerate}
\item $\dr A < \infty$
\item $A$ is $\mathcal{Z}$-stable
\item $A$ has strict comparison of positive elements.
\end{enumerate}
\end{nconj}
\en
These three properties  and their relationships were discussed extensively in \cite{EllToms:regularity}. Some of the implications are known, but none is trivial. The conjecture was verified in \cite{TomsWinter:VI} for the class of Villadsen algebras of the first type.  R{\o}rdam showed in \cite{Ror:UHFII} and \cite{Ror:Z-absorbing} that (iii) is equivalent to
\begin{enumerate}
\item[(iv)] \it $A$ has almost unperforated Cuntz semigroup. 
\end{enumerate} 
In \cite[Theorem~4.5]{Ror:Z-absorbing} he showed that (ii) implies (iv). Little is known about the implication (iii)$\Rightarrow$(ii); it was confirmed for $A$ strongly self-absorbing in \cite[Proposition~6.7]{RorWin:Z-revisited}. There are partial verifications of (ii)$\Rightarrow$(i) (cf.\ \cite[Corollary~8.6]{Win:localizingEC}), but only under additional hypotheses on $A$ (e.g., asking for locally finite decomposition rank and for many projections or few traces); these results factorize through classification theorems, and hence additionally require the Universal Coefficient Theorem (UCT) to hold for $A$. There are no direct proofs known. The situation for the implication (i)$\Rightarrow$(ii) up to now has been similar: it has been established  under additional structural conditions (e.g., for $A$ approximately homogeneous, or when $A$ has many projections and few traces, and only in the presence of the UCT). In the present paper, I shall give a direct proof of  (i)$\Rightarrow$(ii). 

Let me briefly comment on the proposed hypotheses of \ref{finite-con}. Concerning separability, most of the relevant results in the literature are stated in the separable case, and so this will be done here, mostly for convenience. However, the conjecture makes sense also for nonseparable $\mathrm{C}^{*}$-algebras,  at least on rephrasing $\mathcal{Z}$-stability in terms of the existence of a system of embeddings of $\mathcal{Z}$ into $A$ which are almost central with respect to finite subsets and arbitrarily small positive tolerances.  In fact, it can be shown that  the separable version of the conjecture implies the nonseparable one (we do not give any details here; an example of such a passage from the separable to the nonseparable situation can be found  in \cite[Proposition~2.6]{WinterZac:dimnuc}).

Unitality does not seem to be crucial for the equivalence of the properties \ref{finite-con}(i), (ii), and (iv), but (iii) would have to be reformulated in terms of (possibly) unbounded positive tracial functionals. To derive a nonunital version of (the validity of) the conjecture \ref{finite-con} directly from the unital case seems to be possible only to a limited extent (namely, when $A$ is Morita equivalent to a unital $\mathrm{C}^{*}$-algebra). 

Finiteness of $A$ is a necessary condition for $\dr A$ to be finite; in the infinite case, (ii) and (iii) (of \ref{finite-con}) are known to be equivalent by Kirchberg's characterization of purely infinite $\mathrm{C}^{*}$-algebras; cf.\ \cite[Theorem~5]{JiaSu:Z} and \cite[Theorems~4.1.10 and 7.2.6]{Ror:encyc}. In \cite{WinterZac:dimnuc}, J.\ Zacharias and the author will introduce the notion of nuclear dimension, which takes finite values also for infinite $\mathrm{C}^{*}$-algebras; this notion can be used to formulate a version of \ref{finite-con} for finite and infinite $\mathrm{C}^{*}$-algebras at the same time.

To formulate a version of \ref{finite-con} for non-simple $\mathrm{C}^{*}$-algebras seems possible, but subtle. A minimum requirement on $A$ would be that it has no nonelementary quotients or ideals (this is implied by $\mathcal{Z}$-stability). It is worthwhile mentioning that nuclear dimension, just like $\mathcal{Z}$-stability but unlike decomposition rank,  behaves well with respect to extensions, and so  might be more suitable for  formulating a non-simple version of \ref{finite-con}. 

Nuclearity, finally, is necessary for finite decomposition rank. However, the implications (ii)$\Rightarrow$(iii)$\Leftrightarrow$(iv) of \ref{finite-con} are known to hold for exact $\mathrm{C}^{*}$-algebras (a class strictly bigger than that of nuclear ones), and it is conceivable that there is a notion of topological dimension for other than nuclear $\mathrm{C}^{*}$-algebras which makes Conjecture~\ref{finite-con} valid. (The finiteness of known quantities such as stable or real rank is not enough to guarantee $\mathcal{Z}$-stability.)

It is interesting to compare our main result to that of \cite{EllGongLi:apprdiv}, which states that simple approximately homogeneous algebras with base spaces of bounded finite dimension are approximately divisible in the sense of \cite{BlaKumRor:apprdiv}. Clearly, this conclusion cannot quite hold for arbitrary unital simple $\mathrm{C}^{*}$-algebras with finite decomposition rank, since these do not need to have non-trivial projections (while approximately divisible unital algebras have many); however, on regarding $\mathcal{Z}$-stability as a generalization of approximate divisibility (cf.\ \cite{Jia:nonstable-K}), the present result is seen as the appropriate generalization of \cite{EllGongLi:apprdiv}. In fact, it seems clear that the present methods can be adapted to show that (in the simple and unital case) finite decomposition rank does imply approximate divisibility, provided that  this is compatible with the ordered $K_{0}$-group. Remarkably, our proof does not  in any way depend on  classification results, and does not rely on the UCT, whereas \cite{EllGongLi:apprdiv} uses the full force of the classification theorem of \cite{EllGongLi:simple_AH}.   Very recently, Dadarlat, Toms and Phillips have (independently and with very different methods) shown that simple unital AH algebras with bounded topological dimension are $\mathcal{Z}$-stable; their proof, just as ours, is direct and does not use  classification theorems. Whereas our result goes far beyond the class of AH algebras, it does not  subsume the method of \cite{DadPhiToms:AH-Z} since the latter can be generalized to the case of AH algebras with `exponentially slow dimension growth'.

The main result, Theorem~\ref{A}, has consequences in a number of directions; let us consider some of these briefly now.  

First, we see that with  \cite[Corollary~8.1]{Win:localizingEC} we may trade $\mathcal{Z}$-stability for finite decomposition rank; the result is  Corollary~\ref{few-traces-classification}, the classification of separable, simple, unital $\mathrm{C}^{*}$-algebras with finite decomposition rank which satisfy the UCT and for which projections separate tracial states. This generalizes results of \cite{Winter:fintopdim} and \cite{Win:Z-class}.

Moreover, we now see that $\mathcal{Z}$-stability is automatic in the situation of \cite[Corollary~8.4]{Win:localizingEC}, and so the present result completes the classification of $\mathrm{C}^{*}$-algebras associated to smooth, minimal, uniquely ergodic dynamical systems. This in particular covers the irrational rotation algebras as well as Connes's (projectionless) crossed products of odd spheres. 

Finally, we can improve the characterization \cite[Theorem~7.6]{RorWin:Z-revisited} of the Jiang--Su algebra, which can now be described as the uniquely determined separable, simple, nonelementary, unital, monotracial $\mathrm{C}^{*}$-algebra which has finite decomposition rank and is $\mathrm{KK}$-equivalent to the complex numbers. 

The paper is organized as follows. In Section~{\ref{preliminaries}},  the notions of decomposition rank and order zero maps are recalled, and some easy background  results which will be used frequently throughout the paper are provided.  Some useful notation is also introduced. Section~{\ref{Z-stability}} provides a criterion (in terms of generators and relations of dimension drop intervals) for when a $\mathrm{C}^{*}$-algebra is $\mathcal{Z}$-stable. 
The key result of Section~{\ref{tracial-matrix-cone-absorption}} is Lemma~\ref{F}, which says that finite decomposition rank ensures the existence of large matrix cones which are almost central and at the same time large in trace. Together with a result from \cite{TomsWinter:VI}, this  essentially shows that finite decomposition rank implies strict comparison of positive elements. These results are expanded in Section~{\ref{almost-central-dimension-drop-embeddings}} to yield almost central elements which almost satisfy the relations of dimension drop intervals. This result and the result from Section~{\ref{Z-stability}} are  combined in Section~{\ref{main-result}} to yield the main result, Theorem~\ref{A}. Also in Section~{\ref{main-result}},  some corollaries are given and a number of applications  described.  

We are indebted to  D.\ Archey, B.\ Jacelon, A.\ Toms and S.\ White for carefully reading an earlier version of the manuscript and for pointing out some typos and small mistakes. We  would also like to thank the referee for suggesting a number of mathematical and stylistic improvements.


\section{Preliminaries}
\label{preliminaries}

\noindent
Below we recall the concepts of decomposition rank and order zero maps, introduce some new notation in this context and derive some technical results using functional calculus for order zero maps.  We also recall some background information on approximate multiplicativity of completely positive contractive (c.p.c.) approximations, and on tracial states and dimension functions.

\bn
First, recall the following definition from \cite{KirWinter:dr}:

\label{d-dr} 
\begin{ndefn} (cf.\ \cite{KirWinter:dr}, Definitions 2.2 and 3.1) Let $A$ be a separable $\mathrm{C}^{*}$-algebra.
\begin{itemize}
\item[(i)] A completely positive map $\varphi :  F \to A$ has   order zero, $\ord \varphi = 0$, if it preserves orthogonality, i.e., $\varphi(e) \varphi(f) = \varphi(f) \varphi(e) = 0$ for all $e,f \in F$ with $ef = fe = 0$.
\item[(ii)] A completely positive map $\varphi : F \to A$ ($F$ a finite-dimensional $\mathrm{C}^{*}$-algebra) is $n$-decomposable, if there is a decomposition $F=F^{(0)} \oplus \ldots \oplus F^{(n)}$ such that the restriction of $\varphi$ to $F^{(i)}$ has   order zero for each $i \in \{0, \ldots, n\}$; we say $\varphi$ is $n$-decomposable with respect to $F=F^{(0)} \oplus \ldots \oplus F^{(n)}$.
\item[(iii)] $A$ has decomposition rank $n$, $\dr A = n$, if $n$ is the least integer such that the following holds: For any finite subset $\Gh  \subset A$ and $\varepsilon > 0$, there is a completely positive contractive (c.p.c.) approximation $(F, \psi, \varphi)$ for $\Gh$ to within $\varepsilon$ (i.e., $\psi:A \to F$ and $\varphi:F \to A$ are completely positive contractive and $\|\varphi \psi (b) - b\| < \varepsilon \; \forall \, b \in \Gh$) such that $\varphi$ is $n$-decomposable. If no such $n$ exists, we write $\dr A = \infty$.  
\end{itemize}
\end{ndefn}
\en

\bn
We collect below some well known facts about  
order zero maps (see \cite[Proposition~3.2(a)]{Win:cpr}  and 
\cite[1.2]{Winter:fintopdim} for Proposition~\ref{order-zero-facts}, and  
\cite[1.2.3]{Win:cpr2} for Proposition~\ref{order-zero-facts-2}). We let
$e_{ij}$ denote the canonical $(i,j)$-th 
matrix unit in $M_p$.

\begin{nprop}
\label{order-zero-facts}
Let $A$ be a $\mathrm{C}^{*}$-algebra, let $p \in \mathbb{N}$, and let
$\varphi \colon M_{p} \to A$ be a c.p.c.\ order zero map.
\begin{enumerate}
\item There is a unique $*$-homomorphism $\tilde{\varphi} \colon 
\mathcal{C}_{0}((0,1]) \otimes M_{p} \to A$ such that 
$\varphi(x) = \tilde{\varphi}( \iota\otimes x)$
for all $x \in M_{p}$, where $\iota(t) = t$.
\item There is a unique $*$-homomorphism $\pi \colon
M_{p} \to A''$ 
given by sending the matrix unit $e_{ij}$ in $M_p$ to the partial
isometry in $A''$ in the polar decomposition of $\varphi(e_{ij})$. We have 
\begin{equation}
\label{order-zero-facts-w1}
\varphi(x)= \pi(x) \varphi(1_{p}) = 
\varphi(1_{p}) \pi(x) 
\end{equation}
for all $x \in M_p$; and $\pi(1_{p})$ is the support
projection of $\varphi(1_{p})$. 
\item If, for some  $h  \in A''$ with $\|h\|\le 1$, the element
  $h^{*}h$ commutes  
with $\pi(M_{p})$ and satisfies 
$h^{*}h \pi(M_{p}) \subseteq A$, then the map
$\varphi_{h} \colon M_{p} \to A$
given by $ \varphi_{h}(x) = h \pi(x) h^{*}$, for $x \in M_p$,
is a well defined c.p.c.\ order zero map.
\end{enumerate}
\end{nprop}
\en

\bn
\label{order-zero-notation}
\begin{nnotation}
The map $\pi$ in \ref{order-zero-facts}(ii) above will be called the \emph{canonical supporting 
$*$-homomorphism} of $\varphi$. It is clear from \eqref{order-zero-facts-w1} that, whenever $\varphi:M_{p} \to A$ is a c.p.c.\ order zero map and $f \in \Ch_{0}((0,1])$ is a positive function of norm at most 1, then we may define a c.p.c.\ order zero map 
\[
f(\varphi): M_{p} \to A
\]
by setting
\begin{equation}
\label{oznw2}
f(\varphi)(x) := \pi(x) f(\varphi(\be_{M_{p}})).
\end{equation}
On approximating $f$ uniformly by polynomials, \eqref{oznw2} and \eqref{order-zero-facts-w1} yield 
\begin{equation}
\label{oznw1}
f(\varphi)(q) = f(\varphi(q))
\end{equation} 
whenever $q \in M_{p}$ is a projection.
\end{nnotation}
\en

\bn
\label{order-zero-traces}
\begin{nprop}
If $\varphi:F \to A$ is a c.p.c.\ order zero map, then, for  any $\tau \in T(A)$, the map $\tau \circ \varphi$ is a positive tracial functional  of norm at most 1 on $F$.
\end{nprop}

\begin{nproof}
This follows immediately from  \eqref{order-zero-facts-w1},  on observing that 
\begin{eqnarray*}
\tau\varphi(xy) & = & \tau(\pi(xy) \varphi(\be_{F})) \\
& = &  \tau(\pi(x) \pi(y) \varphi(\be_{F})) \\
& = &  \tau(\pi(x) \varphi(\be_{F})^{\halb }\pi(y) \varphi(\be_{F})^{\halb}) \\
& = &  \tau(\pi(y) \varphi(\be_{F})^{\halb }\pi(x) \varphi(\be_{F})^{\halb}) \\
& = & \tau(\varphi(yx)) 
\end{eqnarray*} 
for $x,y \in F$. Here, we have used that $\pi(x) \varphi(\be_{F})^{\halb }$ and $\pi(y) \varphi(\be_{F})^{\halb}$ are in $A$.
\end{nproof}
\en

\bn
\begin{nprop} \label{order-zero-facts-2}
Suppose $x_1, x_{2}, \ldots,x_{p} \in A$ satisfy the relations 
\begin{equation} \tag{$\mathcal{R}_{p}$}
\|x_{i}\| \le 1, \quad x_1 \ge 0, \quad x_ix_i^* = x_1^*x_1, \quad
x_j^*x_j \perp x_i^*x_i,
\end{equation}
for all $i,j=1, \dots, n$ with $i \ne j$. Then, the linear map $\psi
\colon M_p \to A$ given by $\psi(e_{ij}) = x_i^*x_j$ is a c.p.c.\
order zero map. 
\end{nprop}

Note that the original version of the above result was phrased in terms of elements 
of the form $e_{i1}$, $i=2,\ldots,p$. However, it is straightforward to check that the two versions 
are in fact equivalent. 
\en

\bn
\begin{nremark}
It is clear that \ref{order-zero-facts}, \ref{order-zero-notation} and \ref{order-zero-facts-2} can also be phrased for an arbitrary finite-dimensional $\mathrm{C}^{*}$-algebra in place of $M_{p}$, when replacing the matrix units $\{e_{ij}\}$ by matrix units of the form $\{e_{ij}^{(r)}\}$ and  changing $\mathcal{R}_{p}$ accordingly. We shall do so in the sequel without further comment.
\end{nremark}
\en

\bn
By \cite[Proposition~2.5]{KirWinter:dr}, the relations $\mathcal{R}_{p}$ are weakly stable. The following is a straightforward consequence of this fact.

\begin{nprop}
\label{almost-order-zero}
Let a finite-dimensional $\mathrm{C}^{*}$-algebra $F$ and $\gamma>0$ be given. Then, there is $\delta>0$ such that the following statement holds:

If $A$ is a $\mathrm{C}^{*}$-algebra, $h \in A$ with $\|h\|\le 1$ and 
\[
\widetilde{\Phi}:F \to A
\]
is a c.p.c.\ order zero map such that 
\[
\|[h,\widetilde{\Phi}(x)]\| \le \delta \|x\|, \, x \in F,
\]
then there is a c.p.c.\ order zero map
\[
\Phi:F \to \her(hh^{*}) \subset A
\]
such that
\[
\|\Phi(x) - h \widetilde{\Phi}(x) h^{*} \| \le \gamma \|x\|
\]
for all $x \in F$.
\end{nprop}
\en

\bn
\label{I}
\begin{nprop}
Let $\zeta>0$ and $f,g \in \Ch([0,1])$ be given. Then, there is $\delta>0$ such that the following statement holds:

If $d,a$ are positive elements of norm at most one in a unital $\mathrm{C}^{*}$-algebra, and if $\|[d,a\|]\le \delta$, then $\|[f(d),g(a)]\|\le \zeta$. 
\end{nprop}

\begin{nproof}
Immediate on using functional calculus and approximating $f$ uniformly by polynomials.
\end{nproof}
\en

\bn
\label{order-0-almost-commuting}
\begin{nprop}
Let $F_{1}, F_{2}$ be finite dimensional $\mathrm{C}^{*}$-algebras, let $\Gh \subset \Ch_{0}((0,1])$ be a finite subset of positive functions of norm at most 1 and let $r \in \N$ and $\gamma>0$ be given. 

Then, there is $\beta>0$ such that the following statement holds for any $\mathrm{C}^{*}$-algebra $A$:

If 
\[
\varphi_{i}:F_{i} \to A, \; i=1,2,
\]
are c.p.c.\ order zero maps satisfying
\[
\|[\varphi_{1}(x), \varphi_{2}(y)]\| \le \beta \|x\| \|y\|
\]
for all $x \in F_{1}$ and $y \in F_{2}$, then
\begin{enumerate}
\item $\|[ f_{1}(\varphi_{1})(x),f_{2}(\varphi_{2})(y)]\| \le \gamma \|x\| \|y\|$  for all $x \in F_{1}$, $y \in F_{2}$ and $f_{1},f_{2} \in \Gh$
\item $\| [ h(f(\varphi_{1})(x) g(\varphi_{2})(y) f(\varphi_{1})(x) ) , d(\varphi_{2})(y) ] \| \le \gamma$ for all positive normalized  $x \in F_{1}$, $y \in F_{2}$ and $d,f,g,h \in \Gh$
\item $| \tau(h(f(\varphi_{1})(p) g(\varphi_{2})(y) f(\varphi_{1})(p) )) - r\cdot \tau(h(f(\varphi_{1})(q) g(\varphi_{2})(y) f(\varphi_{1})(q) ))| < \gamma $ for all positive normalized $y \in F_{2}$, $f,g,h \in \Gh$, $\tau \in T(A)$ and projections $p,q \in F_{1}$ satisfying $\tr(p) = r \cdot \tr(q) $ for any trace on $F_{1}$
\item $\|d(\varphi_{2})(y) h(f(\varphi_{1})(x) g(\varphi_{2})(y) f(\varphi_{1})(x)) - h(f(\varphi_{1})(x) g(\varphi_{2})(y) f(\varphi_{1})(x))\| \le \gamma$ whenever $x \in F_{1}$, $y \in F_{2}$ are positive and normalized and $d,f,g,h \in \Gh$ satisfy $dg=g$.
\end{enumerate}
\end{nprop}

\begin{nproof}
It is clear that we may construct $\beta$ for (i), (ii), (iii) and (iv) separately.

(i) As 
\[
f_{i}(\varphi_{i})(x_{i}) \in \varphi_{i}(x_{i}) \cdot C^{*}(\varphi_{i}(\be_{F_{i}}))
\]
for $ x_{i} \in F_{i}$ and $i=1,2$, the proof becomes straightforward on approximating elements of $C^{*}(\varphi_{i}(\be_{F_{i}}))$ uniformly by polynomials in $\varphi_{i}(\be_{F_{i}})$; cf.\ Proposition~\ref{I}.

(ii) is a direct consequence of (i) on approximating each $h \in \Gh$ uniformly by polynomials. 

(iv) is proved in the same manner.

(iii) On approximating each $h \in \Gh$ uniformly by polynomials, we obtain $0<\delta$ such that, whenever $a,b \in A$ are positive of norm at most 1 with $\|a-b\|< \delta$, then $\|h(a) - h(b)\| < \frac{\gamma}{2r}$ for any $h \in \Gh$. 

From \cite[Proposition~2.5]{KirWinter:dr} we see that the relations defining order zero maps are weakly stable, so that  there is 
\[
0< \gamma'< \frac{\delta}{2}
\]
such that,  if 
\[
\|[ f(\varphi_{1})(x),g(\varphi_{2})(y)]\| \le \gamma' \|x\| \|y\|
\]
for all $x \in F_{1}$, $y \in F_{2}$ and $f,g \in \Gh$, then there are  c.p.c.\ order zero maps 
\[
\Phi_{f,g}: F_{1} \otimes F_{2} \to A
\]
satisfying 
\begin{equation}
\label{o0acw1}
\|\Phi_{f,g}(x \otimes y) - f^{2}(\varphi_{1})(x) g(\varphi_{2})(y) \| \le \frac{\delta}{2} \|x\| \|y\|
\end{equation}
for any $f,g \in \Gh$, $x \in F_{1}$ and $y \in F_{2}$. 

From (i) we obtain $0<\beta$ such that, if 
\[
\|[\varphi_{1}(x), \varphi_{2}(y)]\| \le \beta \|x\| \|y\|
\]
for all $x \in F_{1}$ and $y \in F_{2}$, then 
\begin{equation}
\label{o0acw2}
\|[ f(\varphi_{1})(x),g(\varphi_{2})(y)]\| \le \gamma' \|x\| \|y\|
\end{equation}
for all $x \in F_{1}$, $y \in F_{2}$ and $f,g \in \Gh$.

Now for $p,q,y,f,g,h$ and $\tau$ as in the assertion of (iii), we see that
\begin{eqnarray*}
\lefteqn{\|\Phi_{f,g}(p \otimes y) - f(\varphi_{1})(p) g(\varphi_{2})(y) f(\varphi_{1})(p) \| } \\
& \le & 
\|\Phi_{f,g}(p \otimes y) - f^{2}(\varphi_{1})(p) g(\varphi_{2})(y)\| \\
&& + \|f^{2}(\varphi_{1})(p) g(\varphi_{2})(y) - f(\varphi_{1})(p) g(\varphi_{2})(y) f(\varphi_{1})(p) \| \\
& \stackrel{\eqref{o0acw1},\eqref{o0acw2}}{\le} &  \frac{\delta}{2} + \gamma' \\
& < & \delta
\end{eqnarray*}
and, similarly,
\[
\|\Phi_{f,g}(q \otimes y) - f(\varphi_{1})(q) g(\varphi_{2})(y) f(\varphi_{1})(q) \| < \delta,
\]
whence
\begin{equation}
\label{o0acw3}
\|h(f(\varphi_{1})(p) g(\varphi_{2})(y) f(\varphi_{1})(p)) - h(\Phi_{f,g}(p \otimes y)) \| < \frac{\gamma}{2}
\end{equation}
and
\begin{equation}
\label{o0acw4}
\|r \cdot ( h(f(\varphi_{1})(q) g(\varphi_{2})(y) f(\varphi_{1})(q)) - h(\Phi_{f,g}(q \otimes y)) )\| < \frac{\gamma}{2}.
\end{equation}
Define a c.p.c.\ order zero map 
\[
\Phi_{f,g,y}(\, .\,) := \Phi_{f,g}(\, . \, \otimes y) : F_{1} \to A
\]
and observe that 
\begin{equation}
\label{o0acw6}
h(\Phi_{f,g,y})(x) \stackrel{\eqref{oznw1}}{=} h(\Phi_{f,g}(x \otimes y))
\end{equation}
whenever $x \in F_{1}$ is a projection. Since $\Phi_{f,g,y}$ has order zero, the map
\[
x \mapsto \tau(h(\Phi_{f,g,y})(x))
\]
defines a positive tracial functional of norm at most 1 on $F_{1}$ by Proposition~\ref{order-zero-traces}; we thus have 
\begin{equation}
\label{o0acw5}
\tau(h(\Phi_{f,g,y})(p)) = r \cdot \tau(h(\Phi_{f,g,y})(q))
\end{equation}
by our assumption on $p$ and $q$. Combining this with the preceding estimates, we obtain
\begin{eqnarray*}
\lefteqn{| \tau(h(f(\varphi_{1})(p) g(\varphi_{2})(y) f(\varphi_{1})(p) )) - r\cdot \tau(h(f(\varphi_{1})(q) g(\varphi_{2})(y) f(\varphi_{1})(q) ))| }\\
& \stackrel{\eqref{o0acw6}}{\le} & \|h(f(\varphi_{1})(p) g(\varphi_{2})(y) f(\varphi_{1})(p)) - h(\Phi_{f,g}(p \otimes y)) \| \\
& & + |\tau(h(\Phi_{f,g,y})(p)) - r \cdot \tau(h(\Phi_{f,g,y})(q))| \\
&& + \|r \cdot (  h(\Phi_{f,g}(q \otimes y)) - h(f(\varphi_{1})(q) g(\varphi_{2})(y) f(\varphi_{1})(q))  )\| \\
& \stackrel{\eqref{o0acw3},\eqref{o0acw4},\eqref{o0acw5}}{<} & \gamma .
\end{eqnarray*}
\end{nproof}
\en

\bn
\label{multiplicative-domain}
We shall have use for the following consequence of Stinespring's theorem, which  is a standard tool to analyze completely positive approximations of nuclear $\mathrm{C}^{*}$-algebras. See \cite{KirWinter:dr}, Lemma 3.5, for a proof.
 
\begin{nlemma}
Let $A$ and $F$ be $\mathrm{C}^{*}$-algebras, $b \in A$ a  positive element of norm at most 1 and $\eta>0$. If $A \stackrel{\psi}{\longrightarrow} F \stackrel{\varphi}{\longrightarrow} A$ are c.p.c.\  maps satisfying 
\[
\|\varphi \psi(b) - b\|, \, \|\varphi \psi(b^{2}) - b^{2}\| \le \eta \, ,
\]
then, for any $x \in F_{+}$, 
\[
\|\varphi(\psi(b)x)- \varphi \psi(b) \varphi(x)\|\le 2 \eta^{\halb} \|x\| \, .
\]
\end{nlemma}
\en

\bn
\label{trace-min}
\begin{nprop}
Let $A$ be a simple unital $\mathrm{C}^{*}$-algebra, and let $\Gh \subset A$ be a compact subset of nonzero positive elements. 

Then, the number 
\[
\min\{\tau(a) \mid \tau \in T(A),\, a \in \Gh\}
\]
exists and is strictly positive.
\end{nprop}

\begin{nproof}
We may regard any $a \in A_{+}$ as a positive continuous function on $T(A)$ via $(\tau \mapsto \tau(a))$. As $A$ is simple, each trace on $A$ is faithful, so 
\[
\Delta_{a}:=\{\tau(a) \mid \tau \in T(A)\}
\] 
is a set of strictly positive real numbers. Since $A$ is unital, $T(A)$ is compact, whence $T_{a}$ is also compact by continuity of $a$. Therefore, 
\[
\delta_{a}:= \min T_{a}
\]
exists and is nonzero for each $0 \neq a \in A_{+}$. One checks that 
\[
a \mapsto \delta_{a}
\]
defines a continuous function on $A_{+}$, the restriction of which to the compact subset $\Gh $ therefore attains a (strictly positive) minimum. 
\end{nproof}
\en

It will be convenient to introduce the following notation.

\bn
\begin{nnotation}
For positive numbers $0\le \eta < \varepsilon \le 1$ define continuous 
functions $$f_{\eta,\varepsilon}, g_{\eta,\varepsilon} \colon [0,1] \to \R^+$$ by
\[
g_{\eta,\varepsilon}(t)= \begin{cases}0, & t \le \eta,\\ 1, & \varepsilon \le t \le 1,
\\ \text{linear,} & \text{else,}\end{cases} 
\]
and
\[
f_{\eta,\varepsilon}(t)= \begin{cases}0, & t \le \eta,\\ t, & \varepsilon \le t \le 1,
\\ \text{linear,} & \text{else.}\end{cases} 
\]
\end{nnotation}
\en

\bn
\label{d-tau-notation}
\begin{nnotation}
If $A$ is a $\mathrm{C}^{*}$-algebra and $\tau$ is a positive tracial functional on $A$, we define the associated lowe semi-continuous dimension function $d_{\tau}$ on $A_{+}$ by 
\[
d_{\tau}(a):= \lim_{n\to \infty} \tau(a^{\frac{1}{n}});
\]
one checks 
\[
d_{\tau}(a) = \lim_{\varepsilon \to 0} \tau(g_{0,\varepsilon}(a))
\]
and
\begin{equation}
\label{dtnw1}
d_{\tau}((a-\beta)_{+}) \le \tau(g_{\beta/2,\beta}(a))
\end{equation}
whenever $\beta>0$. 

Any such dimension function induces a positive real-valued character (also denoted by $d_{\tau}$) on the Cuntz semigroup $W(A)$ with its natural order (cf.\ \cite{Cuntz:dimension}); if $A$ is unital and $\tau$ is a tracial state, then $d_{\tau}:W(A) \to \R_{+}$ is a state. 

We refer the reader to \cite{Ror:Z-absorbing} and \cite{Ror:UHFII} for
notation and background material on Cuntz comparison of positive elements and on the Cuntz semigroup. 
\end{nnotation}
\en

\section{$\mathcal{Z}$-stability}
\label{Z-stability}

\noindent
In this section we recall some facts about the Jiang--Su algebra and about $\mathcal{Z}$-stability; we also derive a new criterion for $\mathcal{Z}$-stability which will be useful for the proof of Theorem~\ref{A}.

\bn
\label{Z-intro}
Recall from \cite{JiaSu:Z} that the Jiang--Su algebra $\mathcal{Z}$ is the uniquely determined simple and  monotracial inductive limit of so-called prime dimension drop $\mathrm{C}^{*}$-algebras, 
\[
\mathcal{Z}= \lim_{\to} Z_{p_{\nu},q_{\nu}},
\]
where
\[
Z_{p_{\nu},q_{\nu}} = \{ f \in \Ch([0,1], M_{p_{\nu},q_{\nu}}) \mid f(0) \in M_{p_{\nu}} \otimes \be_{M_{q_{\nu}}} \mbox{ and }  f(1) \in \be_{M_{q_{\nu}}}  \otimes M_{q_{\nu}}\}, \, \nu \in \N,
\]
and $p_{\nu}$ and $q_{\nu}$ are prime. 

It was shown in \cite{JiaSu:Z} that $\mathcal{Z}$ is strongly self-absorbing in the sense of \cite{TomsWin:ssa}, and that it is $\mathrm{KK}$-equivalent to the complex numbers. 

In \cite{RorWin:Z-revisited}, several alternative characterizations of the Jiang--Su algebra were given; see also \cite{DadToms:Z}. 
In order to show that a unital $\mathrm{C}^{*}$-algebra $A$ is $\mathcal{Z}$-stable, it will suffice to construct approximately central unital $*$-homomorphisms from (arbitrarily large) prime dimension drop intervals to $A$; cf.\ \cite[Proposition~2.2]{TomsWin:ZASH}. Using \cite[Proposition~5.1]{RorWin:Z-revisited}, to this end it will be enough to realize the generators and relations of \cite[Proposition~5.1(iii)]{RorWin:Z-revisited} approximately -- and in an approximately central way (as formalized in \ref{R-relations} below), see Proposition~\ref{B}. In Proposition~\ref{C} we provide a slightly more flexible criterion for the relations  $\Rh(n, \Fh,\eta)$ to hold.
\en

\bn
\label{R-relations}
\begin{nnotation}
Let $A$ be a unital $\mathrm{C}^{*}$-algebra, $n \in \N$, $\eta>0$ and $\Fh \subset A$ a finite subset. If 
\[
\varphi:M_{n} \to A
\]
is a c.p.c.\ order zero map and $v \in A$ such that
\begin{enumerate}
\item $\| v^{*}v - (\be_{A}-\varphi(\be_{M_{n}}))\| < \eta $
\item $\| v v^{*} \varphi(e_{11}) - vv^{*}\|< \eta $
\item $\|[\varphi(x),a]\| < \eta $ for all $a \in \Fh$, and $ x \in M_{n}$ with $\|x\|=1$
\item $\|[v,a]\| < \eta$ for all $a \in \Fh$,
\end{enumerate}
then we say $\varphi$ and $v$ satisfy the relations $\Rh(n, \Fh,\eta)$.
\end{nnotation}
\en

\bn
\label{B}
\begin{nprop}
Let $A$ be a separable and unital $\mathrm{C}^{*}$-algebra. Suppose that, for any $n \in \N$, any finite subset $\Fh \subset A$ and any $0< \eta <1$, there are a c.p.c.\ order zero map 
\[
\varphi: M_{n} \to A
\]
and $v \in A$ satisfying the relations $\Rh(n, \Fh, \eta)$.

Then, $A$ is $\Zh$-stable. 
\end{nprop}

\begin{nproof}
Using separability of $A$, and increasing $\Fh$ and decreasing $\eta$, by $\Rh(n, \Fh,\eta)$ for each $n \in \N$ one obtains a c.p.c.\ order zero map
\[
\tilde{\varphi}: M_{n} \to A_{\infty} \cap A'
\]
and
\[
\tilde{v} \in A_{\infty} \cap A'
\]
such that 
\[
\tilde{v}^{*}\tilde{v} = \be_{A} - \tilde{\varphi}(\be_{M_{n}})
\]
and
\[
\tilde{v}\tilde{v}^{*} \tilde{\varphi}(e_{11}) = \tilde{v} \tilde{v}^{*}.
\]
Now by \cite[Proposition~5.1,(iii)$\Rightarrow$(iv)]{RorWin:Z-revisited}, $Z_{n,n+1}$ embeds unitally into $A_{\infty} \cap A'$. By \cite[Proposition~2.2]{TomsWin:ZASH}, this entails $\Zh$-stability of $A$.
\end{nproof}
\en

\bn
\label{C}
\begin{nprop}
Let $A$ be a  unital $\mathrm{C}^{*}$-algebra, $n \in \N$ and $\Fh \subset A$ a finite subset. Given $0<\eta$ there are $0<\delta<1$ and $0<\zeta<1$ such that the following holds:

If there are a a c.p.c.\ order zero map
\[
\varphi': M_{n} \to A
\]
and
\[
v' \in A
\]
satisfying
\begin{enumerate}
\item $\|(v')^{*} v' - (\be_{A} - \varphi'(\be_{M_{n}}))\| < \delta$,
\item $v' (v')^{*} \in \overline{(\varphi'(e_{11})- \zeta)_{+}A (\varphi'(e_{11})- \zeta)}_{+} $,
\item $\|[\varphi'(x),a]\| \le \delta \|x\|$ for all $ x \in M_{n}, \, a \in \Fh$,
\item $\|[v',a]\| < \delta$ for all $a \in \Fh$,
\end{enumerate}
then there are a c.p.c.\ order zero map 
\[
\varphi:M_{n} \to A
\]
and $v \in A$ of norm at most one satisfying the relations $\Rh(n,\Fh,\eta)$ of \ref{R-relations}. 
\end{nprop}

\begin{nproof}
We may clearly assume the elements of $\Fh$ to be positive and normalized.

Set 
\begin{equation}
\label{wwwC4}
\zeta:= \frac{\eta}{2}.
\end{equation}
Define $g_{\zeta/2,\zeta}^{-}\in \Ch((0,1])$ by
\[
g_{\zeta/2,\zeta}^{-}(t):= \left\{ 
\begin{array}{ll}
0 & t\le \zeta/2\\
\frac{2}{\zeta} - t^{-1} & \zeta/2 \le t \le \zeta \\
t^{-1} & t\ge \zeta.
\end{array}
\right.
\]
Choose $h^{-} \in \Ch_{0}([0,1))$ such that 
\[
h^{-}(t) \cdot (1-t) = 1- g_{\zeta/2,\zeta}(t) \mbox{ and } \|h^{-}\| \le 1+ \zeta.
\]
Obtain 
\begin{equation}
\label{wwwC5}0< \delta < \frac{\eta}{2}
\end{equation} 
such that the assertion of Proposition~\ref{I} holds for both $g_{\zeta/2,\zeta}^{-}$ and $\left(\frac{1}{1+\zeta}\right)^{\frac{1}{2}} (h^{-})^{\frac{1}{2}}$ in place of $f$.

Define $h_{0} \in C^{*}(\varphi'(\be_{M_{n}}))$ by
\[
h_{0}:= g_{\zeta/2,\zeta}^{-}(\varphi'(\be_{M_{n}}));
\]
we then have
\begin{equation}
\label{wwwC7}
h_{0} \varphi'(\be_{M_{n}}) = g_{\zeta/2,\zeta} (\varphi'(\be_{M_{n}}))  \mbox{ and } \|h_{0}\| \le \frac{1}{\zeta};
\end{equation}
Define $h_{1} \in C^{*}(\be_{A} - \varphi'(\be_{M_{n}}))$ by
\begin{equation}
\label{wwwC3}
h_{1}:= h^{-}(\be_{A} - \varphi'(\be_{M_{n}}));
\end{equation}
note that
\begin{equation}
\label{wwwC2}
h_{1} \cdot  (\be_{A} - \varphi'(\be_{M_{n}})) = \be_{A} - g_{\zeta/2,\zeta} (\varphi'(\be_{M_{n}}))  \mbox{ and } \|h_{1}\| \le 1 + \zeta.
\end{equation}
Set 
\begin{equation}
\label{wwwC1}
v:= \left(\frac{1}{1+\zeta}\right)^{\frac{1}{2}} \cdot v' h_{1}^{\frac{1}{2}}
\end{equation}
and define $\varphi:M_{n} \to A$ by
\begin{equation}
\label{wwwC6}
\varphi:= g_{\zeta/2,\zeta}(\varphi'),
\end{equation}
cf.\ \ref{order-zero-notation}. It is clear from \ref{order-zero-notation} that $\varphi$ is a c.p.c.\ order zero map. We continue to check \ref{R-relations}(i)--(iv). 

For \ref{R-relations}(i), note that 
\begin{eqnarray*}
\lefteqn{\|v^{*}v - (\be_{A} - \varphi(\be_{M_{n}}))\|}\\
 & \stackrel{\eqref{wwwC1}}{=} & \left\| \frac{1}{1+\zeta} \cdot h_{1}^{\frac{1}{2}} (v')^{*}v' h_{1}^{\frac{1}{2}} - (\be_{A} - \varphi(\be_{M_{n}})) \right\| \\
& \stackrel{\ref{C}\mathrm{(i)}}{\le} & \left\| \frac{1}{1+\zeta} \cdot h_{1}^{\frac{1}{2}} (\be_{A} - \varphi'(\be_{M_{n}})) h_{1}^{\frac{1}{2}} - (\be_{A} - \varphi(\be_{M_{n}})) \right\|\\
&& +  \delta \\
& \stackrel{\eqref{wwwC2},\eqref{wwwC3}}{\le} &  1- \frac{1}{1+\zeta} +  \delta \\
&  \stackrel{\eqref{wwwC4},\eqref{wwwC5}}{<} & \eta.
\end{eqnarray*}

For (ii), observe that 
\[
vv^{*} \varphi(e_{11}) - vv^{*} {=} vv^{*} g_{\zeta/2,\zeta}(\varphi'(e_{11})) - vv^{*} = 0,
\]
since 
\[
vv^{*} \stackrel{\eqref{wwwC1}}{=} v' h_{1}(v')^{*} \stackrel{\ref{C}\mathrm{(ii)}}{\in} \overline{(\varphi'(e_{11})- {\zeta})_{+}A (\varphi'(e_{11})-{\zeta})_{+} }
\]
and 
\[
\varphi(e_{11}) \stackrel{\eqref{wwwC6}}{=} g_{\zeta/2,\zeta}(\varphi'(e_{11})).
\]

\ref{R-relations}(iii) follows from our choice of $\delta$ and Proposition~\ref{I}, since
\begin{eqnarray*}
\|[\varphi(x),a]\| & \stackrel{\eqref{wwwC7}}{=} & \|[\varphi'(x)h_{0},a]\| \\
& \le & \|[\varphi'(x),a]\| + \|[h_{0},a]\| \|x\| \\
& \stackrel{\ref{C}\mathrm{(iii)},\ref{I}}{\le} & \delta \|x\| + \zeta \|x\| \\
& \stackrel{\eqref{wwwC4},\eqref{wwwC5}}{\le} & \eta \|x\|
\end{eqnarray*}
for any $ x \in M_{n}$ and $a \in \Fh$. 

\ref{R-relations}(iv) follows in essentially the same manner: 
\begin{eqnarray*}
\|[v,a]\| & \stackrel{\eqref{wwwC1}}{=} & \left\| \left[ \left(\frac{1}{1+\zeta}\right)^{\frac{1}{2}} v'h_{1}^{\frac{1}{2}},a\right] \right\| \\
& \le & \left\| \left[ v',a\right] \right\| + \left\| \left[ \left(\frac{1}{1+\zeta}\right)^{\frac{1}{2}} h_{1}^{\frac{1}{2}},a\right] \right\|  \\
& \stackrel{\ref{C}\mathrm{(iv)},\ref{I}}{<} & \delta  + \zeta  \\
& \stackrel{\eqref{wwwC4},\eqref{wwwC5}}{<} & \eta
\end{eqnarray*}
for any  $a \in \Fh$. 
\end{nproof}
\en


\section{Tracial matrix cone absorption}
\label{tracial-matrix-cone-absorption}

\noindent
The main result of this section is Lemma~\ref{F}, which will be a further technical key step in the proof of Theorem~\ref{A}. This will follow from Proposition~\ref{wA} by a geometric series argument resembling that of \cite{Win:Z-class} (cf.\ also \cite{Winter:fintopdim} and \cite{Winter:lfdr}).

Together with an earlier result of \cite{TomsWinter:VI} (Lemma~\ref{n-comparison} below), Lemma~\ref{F} will yield a version of comparison of positive elements `up to a factor' (Corollary~\ref{cor-dr-comparison}). In fact, with the results of this section one could directly prove that (for simple, unital $\mathrm{C}^{*}$-algebras) finite decomposition rank implies strict comparison of positive elements, i.e., implication (i) $\Longrightarrow$ (iii) of Conjecture~\ref{finite-con}.

\bn
\label{wD}
\begin{nprop}
Let $A$ be separable, simple, unital, nonelementary with $\dr A \le m < \infty$; let $f \in A$ a positive element of norm at most 1, $0<\tilde{\eta}<1$ and $k \in \N$ be given.

Then, there are an $m$-decomposable c.p.c.\ approximation (for $\her (f)$) $$(F,\psi,\varphi)$$ of $\{f,f^{2}\}$ to within $\tilde{\eta}$ and a $*$-homomorphism
\[
\widetilde{\Phi}: M_{k} \to F
\]
such that
\begin{equation}
\label{wD1}
\| [\widetilde{\Phi}(x),\psi(f)]\| \le \tilde{\eta} \|x\|
\end{equation}
and
\begin{equation}
\label{wD2}
\tau(\widetilde{\Phi}(\be_{M_{k}}) \psi(f)) \ge \frac{1}{2} \cdot \tau (\psi(f))
\end{equation}
for any  $x \in M_{k}$ and any positive tracial functional $\tau$ on $F$.
\end{nprop}

\begin{nproof}
We may clearly assume that $f$ is nonzero, so that 
\begin{equation}
\label{wD6}
\zeta:= \min \{\tau(f) \mid \tau \in T(A)\}
\end{equation}
exists and is strictly positive by Proposition~\ref{trace-min}. Choose
\[
0< \eta
\] 
such that
\begin{equation}
\label{wD8}
32(m+1) \eta^{\frac{1}{4}} + 4 \eta < \tilde{\eta} \mbox{ and } \eta < \frac{1}{32} \zeta, \, \halb \|f\|.
\end{equation}
We first prove the following

\noindent
{\sc Claim:} There is an $m$-decomposable c.p.c.\ approximation (for $A$), $(F',\psi',\varphi')$, of $\{f,f^{2},g_{0,\eta}(f)\}$ to within $\eta$ such that the following hold:
\begin{enumerate}
\item[1.] $\|\psi'(f^{2}) - \psi'(f)^{2}\| , \, \|\psi'(g_{0,\eta}(f)) - g_{0,\eta}(\psi'(f))\| < \eta$
\item[2.] $\tau(\psi'(f)) \ge \frac{\zeta}{2} \; \forall \, \tau \in T(F')$
\item[3.] if $t \in \sigma(f) \cap [2\eta,1]$, then the projections 
\[
\chi_{(t-\eta,t+\eta)}(\psi_{i}'(f)) \in M_{r_{i}}, \, i=1, \ldots,s,
\]
have rank at least $3k$, where $F'=M_{r_{1}} \oplus \ldots \oplus M_{r_{s}}$ and $\psi_{i}'$ denotes the $i$-th summand of $\psi'$.
\end{enumerate}

To prove this claim, let us employ \cite[Proposition~5.1 and Remark~5.2(ii)]{KirWinter:dr} to obtain a system $(F_{\lambda},\psi_{\lambda},\varphi_{\lambda})_{\lambda \in \N}$ of $m$-decomposable c.p.c.\ approximations (for $A$) such that
\begin{equation} 
\label{wD5}
\|\psi_{\lambda}(ab)-\psi_{\lambda}(a)\psi_{\lambda}(b)\| \stackrel{\lambda \to \infty}{\longrightarrow} 0
\end{equation} 
for all $a,b \in A$ and such that 
\[
\psi_{\lambda}(\be_{A}) = \be_{F_{\lambda}}
\] 
for all $\lambda$. 

Then, there is $\lambda_{0} \in \N$ such that the following hold for any $\lambda \ge \lambda_{0}$:
\begin{enumerate}
\item[(A)] for any $t \in \sigma(f) \cap [2 \eta,1]$, the projections
\[
\chi_{(t-\eta,t+\eta)}(\psi_{\lambda,i}(f))
\]
all have rank at least $3k$ (where $\psi_{\lambda,i}$ denotes the $i$-th component of $\psi_{\lambda}$)
\item[(B)] $\tau(\psi_{\lambda}(f)) \ge \frac{\zeta}{2}$ for all $\tau \in T(F_{\lambda})$.
\end{enumerate}
If no such $\lambda_{0}$ existed, there was a subsequence $(F_{\lambda_{\nu}},\psi_{\lambda_{\nu}},\varphi_{\lambda_{\nu}})_{\nu \in \N}$ of $(F_{\lambda},\psi_{\lambda},\varphi_{\lambda})_{\lambda \in \N}$ such that either (A) or (B) failed, i.e., there were $t_{\nu} \in \sigma(f) \cap [2 \eta,1]$ and indices $i_{\nu}$ such that the projections 
\[
p_{\nu} := \chi_{(t_{\nu}-\eta,t_{\nu}+ \eta)}(\psi_{i_{\nu},\lambda_{\nu}}(f))
\]
had rank less than $3k$ (in case (A) failed), or there were tracial states 
\[
\tau_{\nu} \in T(F_{\lambda_{\nu}})
\]
such that 
\begin{equation}
\label{wD7}
\tau_{\nu}(\psi_{\lambda_{\nu}}(f)) < \frac{\zeta}{2}
\end{equation} 
(in case (B) failed).

Now in the first case, the $t_{\nu}$ converge to some $\bar{t} \in \sigma(f) \cap [2\eta,1]$ by compactness. Let $\bar{p}$ denote the image of $(p_{\nu})_{\nu \in \N}$ in 
\[
C:= \prod F_{\lambda_{\nu}} / \bigoplus F_{\lambda_{\nu}},
\]
and let 
\[
\bar{\bar{\psi}}: A \to C
\]
be the c.p.c.\ map induced by the $\psi_{\lambda_{\nu}}$. By \eqref{wD5}, $\bar{\bar{\psi}}$ in fact is a $*$-homomorphism, and it is clear that
\[
d:= g_{\bar{t}- \eta/2,\bar{t}}(f) - g_{\bar{t},\bar{t}+\eta/2}(f)
\]
satisfies 
\[
\bar{\bar{\psi}}(d) \le \bar{p}.
\]
However, $d$ is nonzero as $\bar{t} \in \sigma(f)$ and $\bar{\bar{\psi}}$ is faithful ($A$ is simple and $\bar{\bar{\psi}}$ is unital, hence nonzero), so $\bar{\bar{\psi}}$ restricts to a faithful $*$-homomorphism 
\[
\tilde{\psi}:\her(d) \to \her(p) \subset C.
\]
On the other hand, $\her(p) = pCp$ is $(3k-1)$-subhomogeneous (this is well-known and follows from the fact that $\prod p_{\nu} F_{\lambda_{\nu}} p_{\nu}$ is $(3k-1)$-subhomogeneous), and $\her(d) \subset A$ is simple and nonelementary. But since $\her(d) \cong \tilde{\psi}(\her(d)) \subset \her(p)$, we see that $\her(d)$ is also $(3k-1)$-subhomogeneous, contradicting the fact that a nonzero hereditary subalgebra of a simple and nonelementary $\mathrm{C}^{*}$-algebra again is simple and nonelementary. It follows that case (A) does not fail.

Now suppose (B) fails. We may then choose a free ultrafilter $\omega \in \beta\N \setminus \N$ along $(\lambda_{\nu})_{\nu \in \N}$. As in \cite[2.4]{Voi:around-qd}, one checks that 
\[
\bar{\tau}(a):= \lim_{\nu \to \omega} \tau_{\nu} (\psi_{\lambda_{\nu}}(a)), \, a \in A,
\] 
defines a tracial state on $A$ for which 
\[
\bar{\tau}(f) \le \frac{\zeta}{2}
\]
by \eqref{wD7}, a contradiction to \eqref{wD6}. Therefore, case (B) can fail neither, and $\lambda_{0}$ as required exists. 

Using \eqref{wD5}, (A) and (B), we can choose $\lambda_{1} \ge \lambda_{0}$ large enough, so that 
\[
(F',\psi',\varphi'):= (F_{\lambda_{1}},\psi_{\lambda_{1}},\varphi_{\lambda_{1}})
\]
will satisfy 1., 2.\ and 3.\ above as desired. This verifies the claim.

Now that we have found $(F'= M_{r_{1}} \oplus \ldots \oplus M_{r_{s}},\psi',\varphi')$, for each $i=1,\ldots,s$ we inductively construct $L_{i} \in \N$, $t_{i,j} \in [2\eta,1]$ and projections $q_{i,j} \in M_{r_{i}}$ for $j=1, \ldots,L_{i}$, such that 
\begin{enumerate}
\item[a)] $q_{i,j}$ has rank at least $3k$ for each $j$
\item[b)] $q_{i,j} \le \chi_{(t_{i,j}-\eta,t_{i,j}+2 \eta]} (\psi_{i}'(f))$ for each $j$
\item[c)] $\sum_{j=1}^{L_{i}} q_{i,j} = \chi_{(\eta,1]}(\psi_{i}'(f))$.
\end{enumerate}
To this end, set 
\[
t_{i,1}:= \min((\sigma(f) \cap [2 \eta,1]) \cup \{1+\eta\});
\]
this number exists by \eqref{wD8} and compactness of the spectrum of $f$.  

If $t_{i,j}$ has been constructed for some $j$, set
\begin{equation}
\label{wD12}
t_{i,j+1}:= \min( (\sigma(f) \cap [t_{i,j}+2 \eta,1]) \cup \{1+ 2 \eta\})
\end{equation}
and
\begin{equation}
\label{wD10}
\bar{t}_{i,j}:= \min\{t_{i,j+1}-\eta,t_{i,j}+2\eta\};
\end{equation}
note that
\begin{equation}
\label{wD9}
t_{i,j} + \eta \le \bar{t}_{i,j} \le t_{i,j} + 2 \eta.
\end{equation}
Define
\begin{equation}
\label{wD11}
q_{i,j}:= \chi_{(t_{i,j}-\eta,\bar{t}_{i,j}]} (\psi_{i}'(f)).
\end{equation}
If $\sigma(f) \cap [t_{i,j}+2 \eta,1]=\emptyset$, set $L_{i}:= j$ and stop the induction process. Note that, by \eqref{wD12}, $t_{i,j+1} \ge t_{i,j}+2 \eta$, whence the induction indeed stops and $L_{i}$ exists. Furthermore, each $q_{i,j}$ has rank at least $3k$ by 3.\ above and \eqref{wD9}, so a) holds. From \eqref{wD10} and \eqref{wD11} we see that b) holds. Finally, using \eqref{wD12} and \eqref{wD10} it is now straightforward to check c).

Therefore, $L_{i}$, the $t_{i,j}$ and the $q_{i,j}$ are as desired.  

Since each $q_{i,j}$ has rank at least $3k$, for $i=1, \ldots,s$ and $j=1, \ldots,L_{i}$ there are $*$-homomorphisms
\begin{equation}
\label{wD13}
\Phi_{i,j}:M_{k} \to q_{i,j} M_{r_{i}} q_{i,j}
\end{equation}
such that
\begin{equation}
\label{wD3}
\tr(\Phi_{i,j}(\be_{M_{k}})) \ge \frac{3}{4} \cdot \tr(q_{i,j}),
\end{equation}
where we write $\tr$ for the normalized trace on any of the matrix blocks $M_{r_{i}}$. 

It is clear from b) that, for $x \in M_{k}$, $i \in \{1, \ldots,s\}$ and $j \in \{1, \ldots,L_{i}\}$, 
\begin{eqnarray*}
\|[\psi_{i}'(f),\Phi_{i,j}(x)]\| & \le & \| [t_{i,j} \cdot q_{i,j}, \Phi_{i,j}(x)]\| + 4 \eta \|x\| \\
& \stackrel{\eqref{wD13}}{=} & 4 \eta \|x\|.
\end{eqnarray*}
Set
\begin{equation}
\label{wD4}
p:= \sum_{i=1}^{s} \sum_{j=1}^{L_{i}} q_{i,j} \stackrel{\mathrm{c)}}{=} \chi_{(\eta,1]} ( \psi'(f))
\end{equation}
and
\begin{equation}
\label{wD14}
F:= pF'p, \, \bar{\psi}(\, . \,):= p \psi'(\, . \,)p \mbox{ and } \bar{\varphi}:= \varphi'|_{F}.
\end{equation}
Then, $\bar{\psi}:A \to F$ and $\bar{\varphi}:F \to A$   are obviously c.p.c.\ maps; $\bar{\varphi}$ is $m$-decomposable, say w.r.t.\ 
\[
\bar{\varphi}= \bar{\varphi}^{(0)} \oplus \ldots \oplus \bar{\varphi}^{(m)} \mbox{ and } F= F^{(0)} \oplus \ldots \oplus F^{(m)}.
\] 
We have
\begin{eqnarray}
\|\bar{\varphi} \bar{\psi}(f) - f\| & \le & \| \bar{\varphi} \bar{\psi}(f) - \varphi' \psi'(f)\| + \| \varphi' \psi'(f) - f\| \nonumber \\
& \stackrel{\eqref{wD14}}{\le} & \|\bar{\psi}(f) - \psi'(f) \| +  \| \varphi' \psi'(f) - f\| \nonumber \\
& < & 2 \eta; \label{wD15}
\end{eqnarray}
moreover,
\begin{eqnarray}
\|\bar{\varphi} \bar{\psi}(f^{2}) - f^{2}\| & \stackrel{\eqref{wD14}}{=} & \|\varphi'(p\psi'(f^{2})p) - f^{2} \| \nonumber  \\
& \stackrel{\mathrm{1.}}{<} &  \|\varphi'(p\psi'(f)^{2}p) - f^{2} \| + \eta \nonumber \\
& \stackrel{\eqref{wD4}}{\le} &  \|\varphi'(\psi'(f)^{2}) - f^{2} \| + 2\eta \nonumber \\
& \stackrel{\mathrm{1.}}{<} &  \|\varphi'\psi'(f^{2}) - f^{2} \| + 3\eta \nonumber  \\
& < & 4 \eta, \label{wD16}
\end{eqnarray}
where for the last estimates of \eqref{wD15} and \eqref{wD16} we have used \eqref{wD4} and the fact that $(F',\psi',\varphi')$ approximates $\{f,f^{2}\}$  to within $\eta$; cf.\ the claim above.

Before constructing $\widetilde{\Phi}$, we need to slightly modify the c.p.c.\ approximation $(F,\bar{\psi},\bar{\varphi})$ to actually obtain a c.p.c.\ approximation for $\her(f)$;  this will be accomplished by restricting $\bar{\psi}$ to $\her(f)$ and slightly perturbing $\bar{\varphi}$ to a map with range in $\her(f)$.  

So, choose 
\[
0< \gamma
\]
such that
\begin{equation}
\label{wD18}
h:= g_{0,\gamma}(f)
\end{equation}
satisfies 
\begin{equation}
\label{wD17}
\|(\be_{A}-h) g_{0,\eta}(f)\| < \eta.
\end{equation} 
We have, for $l=0, \ldots,m$, 
\begin{eqnarray}
\lefteqn{\|\bar{\varphi}^{(l)}(\be_{F^{(l)}}) - \bar{\varphi}^{(l)}(\be_{F^{(l)}})h\|} \nonumber \\
& \le & \|(\be_{A}-h) \bar{\varphi}^{(l)}(\be_{F^{(l)}}) (\be_{A}-h) \|^{\halb} \nonumber \\
& \stackrel{\eqref{wD4},\eqref{wD14}}{\le} & \|(\be_{A}-h) \varphi'(g_{0,\eta}(\psi'(f))) (\be_{A}-h) \|^{\halb} \nonumber \\
& \stackrel{\mathrm{1.}}{\le} & (\|(\be_{A}-h) \varphi'\psi' (g_{0,\eta}(f)) (\be_{A}-h) \| + \eta)^{\halb} \nonumber \\
& \le & (\|(\be_{A}-h) g_{0,\eta}(f) (\be_{A}-h) \| + 2\eta)^{\halb} \nonumber \\
& \stackrel{\eqref{wD17}}{\le} & (3 \eta)^{\halb} \nonumber \\
& < & 2 \eta^{\halb}, \label{wD17a}
\end{eqnarray}
where for the third last inequality we have used that $(F',\psi',\varphi')$ also approximates $g_{0,\eta}(f)$ to within $\eta$. 

Now by \cite[Lemma~3.6]{KirWinter:dr}, we see from \eqref{wD17a} that there are c.p.c.\ order zero maps
\[
\hat{\varphi}^{(l)}: F^{(l)} \to \her(h) (\stackrel{\eqref{wD18}}{=}\her(f))
\]
such that 
\begin{equation}
\label{wD19}
\|\hat{\varphi}^{(l)}(x) - \bar{\varphi}^{(l)}(x) \| \le 16 \eta^{\frac{1}{4}} \|x\| \; \forall \, x \in F^{(l)}, \, l=0, \ldots,m.
\end{equation}
The $\hat{\varphi}^{(l)}$ sum up to an $m$-decomposable c.p.\ map
\[
\hat{\varphi}: F \to \her(h)
\]
satisfying
\begin{equation}
\label{wD20}
\| \hat{\varphi}(x) - \bar{\varphi}(x)\| \stackrel{\eqref{wD19}}{\le} (m+1) 16 \eta^{\frac{1}{4}} \|x\|, \, x \in F.
\end{equation}
The map $\hat{\varphi}$ is not necessarily a contraction, but it has norm at most $1+(m+1) 16 \eta^{\frac{1}{4}}$ (by \eqref{wD20}, and since $\bar{\varphi}$ is contractive), so
\[
\varphi := \frac{1}{1+(m+1)16 \eta^{\frac{1}{4}}} \cdot \hat{\varphi}: F \to \her(f)
\]
will be an $m$-decomposable c.p.c.\ map; one checks that
\begin{equation}
\label{wD21}
\|\varphi(x) - \hat{\varphi}(x)\| \le (m+1) 16 \eta^{\frac{1}{4}} \|x\|, \, x \in F.
\end{equation}
Letting $\psi$ denote the restriction of $\bar{\psi}$ to $\her(f)$, 
\begin{equation}
\label{wD23}
\psi:= \bar{\psi}|_{\her(f)},
\end{equation}
we obtain an $m$-decomposable c.p.c.\ approximation $(F, \psi,\varphi)$ for $\her(f)$, which satisfies
\begin{eqnarray*}
\|\varphi \psi(f) - f\| & \le & \|\varphi \psi(f) - \hat{\varphi}\psi(f)\| + \|\hat{\varphi}\psi(f) - \bar{\varphi}\psi(f)\| \\
&& + \| \bar{\varphi}\bar{\psi}(f) -   f\| \\
& \stackrel{\eqref{wD21},\eqref{wD20},\eqref{wD14},\eqref{wD4},\eqref{wD15}}{<} & 32(m+1) \eta^{\frac{1}{4}} + 2 \eta \\
& \stackrel{\eqref{wD8}}{<} & \tilde{\eta} 
\end{eqnarray*}
and, similarly (but using \eqref{wD16} in place of \eqref{wD15}),
\[
\|\varphi \psi(f^{2}) - f^{2} \| < 32(m+1) \eta^{\frac{1}{4}} + 4 \eta  \stackrel{\eqref{wD8}}{<}  \tilde{\eta}. 
\]
Therefore, $(F,\psi,\varphi)$ approximates $\{f,f^{2}\}$ to within $\eta$.

Finally, we may define a $*$-homomorphism 
\[
\widetilde{\Phi}: M_{k} \to F
\]
by setting
\begin{equation}
\label{wD22}
\widetilde{\Phi}(x):= \bigoplus_{i=1}^{s} \bigoplus_{j=1}^{L_{i}} \Phi_{i,j}(x), \, x \in M_{k}.
\end{equation}
From \eqref{wD3} and \eqref{wD4}, we see that 
\begin{eqnarray*}
\tau(\widetilde{\Phi}(\be_{M_{k}}) \psi(f)) & \stackrel{\eqref{wD4}}{\ge} & \tau \left(\widetilde{\Phi}(\be_{M_{k}}) \sum_{i,j} t_{i,j} \cdot q_{i,j} \right) - \eta \\
& \stackrel{\eqref{wD22},\eqref{wD3}}{\ge} & \frac{3}{4} \cdot \tau \left(\sum_{i,j} t_{i,j} \cdot q_{i,j} \right) - \eta \\
& \stackrel{\mathrm{b),c)}}{\ge} & \frac{3}{4} \cdot \tau(\psi'(f)) - 3 \eta \\
& \stackrel{\eqref{wD8},\eqref{wD6}}{\ge} & \halb \cdot \tau(\psi'(f)) + \eta \\
& \stackrel{\eqref{wD23},\eqref{wD14},\eqref{wD4}}{\ge} & \halb \cdot \tau(\psi(f))
\end{eqnarray*}
for any tracial state $\tau \in T(F')$. The estimate clearly remains true for any positive tracial functional on $F'$. Since any positive tracial functional on $F$ extends to one on $F'$, we can thus conclude that
\[
\tau(\widetilde{\Phi}(\be_{M_{k}}) \psi(f)) \ge \halb \cdot \tau(\psi(f))
\]
for any positive tracial functional on $F$.
\end{nproof}
\en

\bn
\label{wF}
\begin{nprop}
Let $A$ be separable, simple, unital, nonelementary with $\dr A \le m < \infty$; let $F$ be a finite-dimensional $\mathrm{C}^{*}$-algebra, $\varphi:F \to A$ a c.p.c.\ order zero map, $0<\bar{\eta}<1$ and $k \in \N$ be given.

Then, there are an $m$-decomposable c.p.c.\ approximation (for $\her (\varphi(\be_{F}))$) $$(\bar{F} = \bar{F}^{(0)}\oplus \ldots \oplus \bar{F}^{(m)},\bar{\psi},\bar{\varphi})$$ of $\varphi(\Bh_{1}(F)) \cup \varphi(\Bh_{1}(F))^{2}$ to within $\bar{\eta}$ and  a $*$-homomorphism
\[
\bar{\Phi}: M_{k} \to \bar{F}
\]
such that
\begin{equation}
\label{wwF7}
\| [\bar{\Phi}(x),\bar{\psi}\varphi(b)]\| \le \bar{\eta} \|x\| \|b\|
\end{equation}
and
\begin{equation}
\label{wwF8}
\tau\left(\sum_{j=0}^{m}\bar{\varphi}^{(j)}(\bar{\Phi}^{(j)}(\be_{M_{k}}) \bar{\psi}^{(i)} \varphi(c)) \right) \ge \frac{1}{4} \cdot \tau \left(\sum_{j=0}^{m} \bar{\varphi}^{(j)} \bar{\psi}^{(j)} \varphi(c) \right)
\end{equation}
for any $x \in M_{k}$, $b \in F$, $c \in F_{+}$ and $\tau \in T(A)$.
\end{nprop}

\begin{nproof}
A moment's thought shows that, since $\varphi$ has order zero, it is enough to prove the assertion when $F$ is of the form $M_{n}$ for some $n \in \N$, and that we may assume $\varphi$ to be nonzero. Let 
\[
\{e_{rs}\}
\]
be a set of matrix units for $M_{n}$, and let 
\begin{equation}
\label{wwF9}
\pi: M_{n} \to A''
\end{equation}
be the canonical supporting $*$-homomorphism for $\varphi$. Set
\[
f:= \varphi(e_{11})
\]
and choose
\begin{equation}
\label{wwF6}
0< \tilde{\eta} < \frac{\bar{\eta}}{n^{2}}, \, \frac{\zeta}{16n^{2}},
\end{equation}
where
\begin{equation}
\label{wwF15}
\zeta:= \min\{\tau(\varphi(e_{11})) \mid \tau \in T(A) \}
\end{equation}
exists and is nonzero by Proposition~\ref{trace-min}.
Choose
\begin{equation}
\label{wwF16}
0 < \gamma \le \frac{\tilde{\eta}}{4}.
\end{equation}
By \ref{order-zero-facts}, we have 
\begin{equation}
\label{wwF1}
C^{*}(\varphi(M_{n})) \cong C^{*}(f) \otimes M_{n}.
\end{equation} 
Since the tensor product of two approximately multiplicative maps is approximately multiplicative, it is not hard to show  that there is 
\[
0< \hat{\eta} < \tilde{\eta}
\]
such that the following holds: If, for some $\mathrm{C}^{*}$-algebra $F'$,  
\[
\psi': \her(f) \to F'
\]
is c.p.c.\ such that 
\begin{equation}
\label{wwF10}
\|\psi'(f^{2}) - \psi'(f)^{2}\| < \hat{\eta},
\end{equation}
then 
\[
\bar{\psi}:= \psi' \otimes \id_{M_{n}}: \her(\varphi(\be_{M_{n}})) \stackrel{\eqref{wwF1}}{\cong} \her(f) \otimes M_{n} \to F' \otimes M_{n}
\]
satisfies 
\begin{equation}
\label{wwF11}
\| \bar{\psi}(g_{0,\gamma}(\varphi)(u) \varphi(v)) - \bar{\psi}(g_{0,\gamma}(\varphi)(u)) \bar{\psi}(\varphi(v))\| < \gamma
\end{equation}
for all $u,v \in M_{n}$ with $\|u\|,\|v\| \le 1$. (In other words, if $\psi'$ is sufficiently multiplicative on $C^{*}(f)$, then $\bar{\psi}$ is almost multiplicative on $C^{*}(\varphi)$).

Now, apply Proposition~\ref{wD} to obtain an $m$-decomposable c.p.c.\ approximation (for $\her(f)$)
\begin{equation}
\label{wwF2}
\left(F'=F^{(0)}\oplus \ldots \oplus F^{(m)}, \psi'= \bigoplus_{i=0}^{m} \psi^{(i)}, \varphi'= \sum_{i=0}^{m} \varphi^{(i)}\right)
\end{equation}
of $\{f,f^{2}\}$ to within $\tilde{\eta}$ and a $*$-homomorphism
\[
\widetilde{\Phi}: M_{k} \to F
\]
satisfying \eqref{wD1} and \eqref{wD2}. Define
\[
\bar{F}:=  F \otimes M_{n}, \, \bar{F}^{(i)}:= F^{(i)} \otimes M_{n},\,  \bar{\varphi}^{(i)}:= \varphi^{(i)} \otimes \pi : \bar{F}^{(i)} \to \her(\varphi(\be_{M_{n}})) \subset A
\]
and
\[
\bar{\psi}^{(i)}: \her(\varphi(\be_{M_{n}})) \stackrel{\eqref{wwF1}}{\cong} \her(f) \otimes M_{n} \to \bar{F}^{(i)}
\]
by
\begin{equation}
\label{wwF3}
\bar{\psi}^{(i)}:= \psi^{(i)} \otimes \id_{M_{n}}.
\end{equation}
Note that by \ref{order-zero-facts}, each $\bar{\varphi}^{(i)}$ indeed is a c.p.c.\ order zero map from $\bar{F}^{(i)}$ to $A$; note also that 
\[
\bar{\varphi}:= \sum_{i=0}^{m} \bar{\varphi}^{(i)} \mbox{ and } \bar{\psi}:= \sum_{i=0}^{m} \bar{\psi}^{(i)}
\]
are c.p.c.\ maps, and that 
\[
\label{wwF14}
\| \bar{\varphi} \bar{\psi} \varphi(b) - \varphi(b) \| \le n^{2} \tilde{\eta} \|b\| \mbox{ and } \| \bar{\varphi} \bar{\psi} (\varphi(b) \varphi(b')) - \varphi(b) \varphi(b')\| \le n^{2} \tilde{\eta} \|b\| \|b'\|
\]
for $b,b' \in F$ by \eqref{wwF2}. Therefore, $(\bar{F},\bar{\psi},\bar{\varphi})$ indeed is a c.p.c.\ approximation with the desired properties.

Next, define a $*$-homomorphism 
\[
\bar{\Phi}:M_{k} \to \bar{F}
\]
by
\begin{equation}
\label{wwF4}
\bar{\Phi}(x) :=  \widetilde{\Phi}(x) \otimes \be_{M_{n}};
\end{equation}
let
\[
\bar{\Phi}^{(i)}:M_{k} \to \bar{F}^{(i)}
\]
denote the $i$th component. For $r,s \in \{1, \ldots,n\}$ and $x \in M_{k}$ we have 
\begin{eqnarray}
\|[\bar{\Phi}(x),\bar{\psi}\varphi(e_{r,s})]\| &\le &  \max_{i \in \{0, \ldots,m\}} \|[\bar{\Phi}^{(i)}(x), \bar{\psi}^{(i)}\varphi(e_{r,s})]\| \nonumber \\
& \stackrel{\eqref{wwF3},\eqref{wwF4},\eqref{wwF1}}{=} & \max_{i \in \{0, \ldots,m\}} \|[\widetilde{\Phi}^{(i)}(x) \otimes \be_{M_{n}}, \psi^{(i)}(f) \otimes e_{rs})]\| \nonumber \\
& \stackrel{\eqref{wD1}}{\le} & \tilde{\eta} \|x\|, \label{wwF5}
\end{eqnarray}
whence
\begin{eqnarray}
\|[\bar{\Phi}(x),\bar{\psi}\varphi(b)]\| & \stackrel{\eqref{wwF5}}{\le} & n^{2} \tilde{\eta} \|x\| \|b\| \label{wwF12} \\
& \stackrel{\eqref{wwF6}}{\le} & \bar{\eta} \|x\| \|b\| \nonumber 
\end{eqnarray}
for $x \in M_{k}$, $b \in M_{n}$. This confirms \eqref{wwF7}.

To verify \eqref{wwF8}, take $c \in (M_{n})_{+}$; we may assume $\|c\| = 1$. Let $u \in \Uh(M_{n})$ be a diagonalizing unitary for $c$, so that we have 
\begin{equation}
\label{wwF13}
u^{*} c u = \sum_{r=1}^{n} \lambda_{r} \cdot e_{r,r} \mbox{ with } 1 \ge \lambda_{1} \ge \ldots \ge \lambda_{n} \ge 0.
\end{equation}
For $\tau \in T(A)$ we then get
{
\allowdisplaybreaks
\begin{eqnarray*}
\lefteqn{\tau\left(\sum_{j=0}^{m} \bar{\varphi}^{(j)}(\bar{\Phi}^{(j)}(\be_{M_{k}}) \bar{\psi}^{(j)} \varphi(c))\right) } \\
& = & \tau\left(\sum_{j=0}^{m} \bar{\varphi}^{(j)}(\bar{\Phi}^{(j)}(\be_{M_{k}}) \bar{\psi}^{(j)} \varphi(c u u^{*}))\right) \\
& \stackrel{\eqref{wwF9}, \ref{order-zero-facts}}{=} &  \tau\left(\sum_{j=0}^{m} \bar{\varphi}^{(j)}(\bar{\Phi}^{(j)}(\be_{M_{k}}) \bar{\psi}^{(j)} \varphi(\be_{M_{n}}) \pi(c u u^{*}))\right) \\
& \ge & \left| \tau \left(\sum_{j=0}^{m} \bar{\varphi}^{(j)}(\bar{\Phi}^{(j)}(\be_{M_{k}}) \bar{\psi}^{(j)}    g_{0,\gamma}(\varphi(\be_{M_{n}}))\varphi(\be_{M_{n}}) \pi(c u u^{*}))\right)\right| - \gamma\\
& \stackrel{\ref{order-zero-facts}}{=} & \left|\tau \left(\sum_{j=0}^{m} \bar{\varphi}^{(j)}(\bar{\Phi}^{(j)}(\be_{M_{k}}) \bar{\psi}^{(j)} (g_{0,\gamma}(\varphi)(c u) \varphi(u^{*})))\right)\right| - \gamma \\
& \stackrel{\eqref{wwF2},\eqref{wwF10},\eqref{wwF11}}{\ge} & \left|\tau \left(\sum_{j=0}^{m} \bar{\varphi}^{(j)}(\bar{\Phi}^{(j)}(\be_{M_{k}}) \bar{\psi}^{(j)} (g_{0,\gamma}(\varphi)(c u) ) \bar{\psi}^{(j)}(\varphi(u^{*})))\right)\right| - 2\gamma \\
& \stackrel{\ref{order-zero-traces}}{=} & \left| \tau \left(\sum_{j=0}^{m} \bar{\varphi}^{(j)}(\bar{\psi}^{(j)} \varphi(u^{*})  \bar{\Phi}^{(j)}(\be_{M_{k}}) \bar{\psi}^{(j)} (g_{0,\gamma}(\varphi)(c u) ) )\right) \right| - 2\gamma \\
& \stackrel{\eqref{wwF12}}{\ge} &  \left| \tau \left(\sum_{j=0}^{m} \bar{\varphi}^{(j)}( \bar{\Phi}^{(j)}(\be_{M_{k}}) \bar{\psi}^{(j)} \varphi(u^{*}) \bar{\psi}^{(j)} (g_{0,\gamma}(\varphi)(c u) ) )\right) \right| - 2\gamma - n^{2} \tilde{\eta}\\
& \stackrel{\eqref{wwF2},\eqref{wwF10},\eqref{wwF11}}{\ge} &   \tau \left(\sum_{j=0}^{m} \bar{\varphi}^{(j)}( \bar{\Phi}^{(j)}(\be_{M_{k}}) \bar{\psi}^{(j)} (\varphi(u^{*}cu) ) ) \right) - 4\gamma - n^{2} \tilde{\eta}\\
& \stackrel{\eqref{wwF13},\eqref{wwF4}}{=} & \tau \left(\sum_{j=0}^{m} \sum_{r=1}^{n} \lambda_{r} \bar{\varphi}^{(j)}(\widetilde{\Phi}^{(j)}(\be_{M_{k}}) \bar{\psi}^{(j)}(\varphi(e_{11}))) \right) - 4\gamma - n^{2} \tilde{\eta}\\
& \stackrel{\eqref{wD2}}{\ge} & \frac{1}{2} \cdot \tau \left(\sum_{j=0}^{m} \sum_{r=1}^{n} \lambda_{r} \bar{\varphi}^{(j)} \bar{\psi}^{(j)}(\varphi(e_{11})) \right) - 4\gamma - n^{2} \tilde{\eta}\\
& \stackrel{\eqref{wwF14}}{\ge} & \frac{1}{2} \cdot \tau \left(\sum_{r=1}^{n} \lambda_{r} \varphi(e_{11}) \right) - 4\gamma - n^{2} \tilde{\eta}- n^{2} \tilde{\eta}\\
& = & \frac{1}{2} \cdot \tau \left(\varphi \left(\sum_{r=1}^{n} \lambda_{r} e_{r,r}\right) \right) - 4\gamma - n^{2} \tilde{\eta}- n^{2} \tilde{\eta}\\
& \stackrel{\eqref{wwF13}}{=} & \frac{1}{2} \cdot \tau(\varphi(c))  - 4\gamma - n^{2} \tilde{\eta}- n^{2} \tilde{\eta}\\
& \stackrel{\eqref{wwF6},\eqref{wwF15},\eqref{wwF16}}{\ge} & \frac{1}{4} \cdot \tau(\varphi(c)) + n^{2} \tilde{\eta} \\
& \stackrel{\eqref{wwF14}}{\ge} &  \frac{1}{4} \cdot  \tau(\bar{\varphi}\bar{\psi} \varphi(c)).
\end{eqnarray*}
}
\end{nproof}
\en

\bn
\begin{nnotation}
We denote by $T_{\infty}(A)$ those tracial states on $A_{\infty}$ which are of the form
\[
[(a_{l})_{l \in \N}] \mapsto \lim_{\omega} \tau_{l}(a_{l})
\]
for a free ultrafilter $\omega $ on $ \N$ and a sequence $(\tau_{l})_{l \in \N} \subset T(A)$. Note that any $\tau \in T_{\infty}(A)$ restricts to a tracial state on $A$.  
\end{nnotation}
\en

\bn
\label{wG}
\begin{nprop}
Let $k,m \in \N$, $\varepsilon >0$ and a separable, simple, unital and nonelementary $\mathrm{C}^{*}$-algebra $A$ with $\dr A \le m$ be given. 

Then, for $i,j \in \{0, \ldots, m\}$ there are elements $$\tilde{h}^{(i)} \in (A_{\infty})_{+}$$ and c.p.c.\ order zero maps $$\widetilde{\Phi}^{(i,j)}: M_{k} \to A_{\infty}$$ with the following properties:
\begin{enumerate}
\item $\sum_{i=0}^{m} \tilde{h}^{(i)} = \be_{A_{\infty}}$
\item $\tilde{h}^{(i)} \in A_{\infty} \cap A'$
\item $[\widetilde{\Phi}^{(i,j)}(M_{k}),\tilde{h}^{(i)}A] = 0$ (note that this implies $[\widetilde{\Phi}^{(i,j)}(M_{k}),\tilde{h}^{(i)}] = 0$ since $A$ is unital)
\item $\tau(\sum_{i,f=0}^{m} \widetilde{\Phi}^{(i,j)}(\be_{M_{k}}) \tilde{h}^{(i)}a) \ge \frac{1}{4} \cdot  \tau(a)$ for any $\tau \in T_{\infty}(A)$ and $a \in A_{+}$.
\end{enumerate}
\end{nprop}

\begin{nproof}
Let 
\begin{equation}
\label{wG2}
\be_{A} \in \Fh \subset A
\end{equation}  
be a compact subset of positive elements of norm at most 1, and let $0<\eta <1$. Let
\begin{equation}
\label{wG3}
(F=F^{(0)}\oplus \ldots \oplus F^{(m)},\psi,\varphi)
\end{equation}  
be an $m$-decomposable c.p.c.\ approximation of $\Fh \cup \Fh^{2}$ to within $\eta$. Set
\begin{equation}
\label{wG1}
h^{(i)}:= \varphi(\be_{F^{(i)}}), \, i=0, \ldots,m.
\end{equation}
For each $i$, apply Proposition~\ref{wF} to $F^{(i)}$ and $\varphi^{(i)}$ (in place of $F$ and $\varphi$, respectively) to obtain a c.p.c.\ approximation (for $\her(\varphi^{(i)}(\be_{F^{(i)}}))$)
\begin{equation}
\label{wG8}
(\bar{F}^{(i)}= \bar{F}^{(i,0)} \oplus \ldots \oplus \bar{F}^{(i,m)}, \bar{\psi}^{(i)},\bar{\varphi}^{(i)}) 
\end{equation}
of $\varphi^{(i)}(\Bh_{1}(F^{(i)})) \cup \varphi^{(i)}(\Bh_{1}(F^{(i)}))^{2}$ to within $\eta$ and a $*$-homomorphism 
\[
\bar{\Phi}^{(i)}: M_{k} \to \bar{F}^{(i)}
\]
such that 
\begin{equation}
\label{wG9}
\| [\bar{\Phi}^{(i)}(x), \bar{\psi}^{(i)} \varphi^{(i)}(b)]\| \le \eta \|x\| \|b\|
\end{equation}
for $x \in M_{k}$, $b \in F^{(i)}$ and
\begin{equation}
\label{wG10}
\tau \left(\sum_{j=0}^{m} \bar{\varphi}^{(i,j)}(\bar{\Phi}^{(i,j)}(\be_{M_{k}}) \bar{\psi}^{(i,j)}\varphi^{(i)}(c)) \right) \ge \frac{1}{4} \cdot \tau \left( \sum_{j=0}^{m} \bar{\varphi}^{(i,j)} \bar{\psi}^{(i,j)} \varphi^{(i)}(c)\right)
\end{equation}
for $c \in (F^{(i)})_{+}$, $\tau \in T(A)$ and $i \in \{0, \ldots,m\}$. 

For $i,j \in \{0, \ldots,m\}$, we define a c.p.c.\ order zero map 
\[
\Phi^{(i,j)}:M_{k} \to A
\]
by
\begin{equation}
\label{wG7}
\Phi^{(i,j)}:= \bar{\varphi}^{(i,j)} \circ \bar{\Phi}^{(i,j)}.
\end{equation}
We have
\[
\be_{A} - \eta \stackrel{\eqref{wG2},\eqref{wG3}}{\le} \sum_{i=0}^{m} \varphi^{(i)}\psi^{(i)}(\be_{A}) \le \sum_{i=0}^{m} \varphi^{(i)}(\be_{F^{(i)}}) \stackrel{\eqref{wG1}}{=} \sum_{i=0}^{m} h^{(i)} \le \be_{A},
\]
whence
\begin{equation}
\label{wG11}
\left\| \be_{A} - \sum_{i=0}^{m} h^{(i)} \right\| \le \eta.
\end{equation}
Moreover, we have 
\begin{eqnarray}
\| \varphi^{(i)}(\be_{F^{(i)}}) a - \varphi^{(i)} \psi^{(i)}(a)\| & \stackrel{\eqref{wG3}}{\le} & \|\varphi^{(i)}(\be_{F^{(i)}}) \varphi \psi(a) - \varphi(\be_{F^{(i)}} \psi(a)) \| + \eta \nonumber \\
& \stackrel{\eqref{wG3}, \ref{multiplicative-domain}}{\le} & 2 \eta^{\halb} + \eta  \label{wG4}
\end{eqnarray}
for $i \in \{0, \ldots,m\}$ and $a \in \Fh$; by the same argument we have
\begin{equation}
\label{wG5}
\|a \varphi^{(i)}(\be_{F^{(i)}}) - \varphi^{(i)}\psi^{(i)}(a)\| \le 2 \eta^{\halb} + \eta.
\end{equation}
As a consequence,
\begin{equation}
\label{wG12}
\|[a,h^{(i)}]\| \stackrel{\eqref{wG1},\eqref{wG4},\eqref{wG5}}{<} 6 \eta^{\halb}
\end{equation}
for $i \in \{0, \ldots,m\}$ and $a \in \Fh$. 

Next, we check that
\begin{eqnarray}
\lefteqn{\|[\Phi^{(i,j)}(x), h^{(i)}a]\|}\\
& \stackrel{\eqref{wG7},\eqref{wG1},\eqref{wG3}}{\le} & \| \bar{\varphi}^{(i,j)} \bar{\Phi}^{(i,j)}(x) \varphi^{(i)}(\be_{F^{(i)}}) \varphi \psi(a) \nonumber \\
&& - \varphi^{(i)}(\be_{F^{(i)}}) \varphi \psi(a) \bar{\varphi}^{(i,j)} \bar{\Phi}^{(i,j)}(x) \| \nonumber \\
&& + 2 \eta \nonumber  \\
& \stackrel{\eqref{wG3},\ref{multiplicative-domain}}{\le} &  \| \bar{\varphi}^{(i,j)} \bar{\Phi}^{(i,j)}(x) \varphi^{(i)} \psi^{(i)}(a) - \varphi^{(i)}\psi^{(i)}(a) \bar{\varphi}^{(i,j)} \bar{\Phi}^{(i,j)}(x) \|   \nonumber \\
&& + 2 \eta  + 4 \eta^{\halb } \nonumber \\
& \stackrel{\eqref{wG8}}{\le} & \| \bar{\varphi}^{(i,j)} \bar{\Phi}^{(i,j)}(x) \bar{\varphi}^{(i)} \bar{\psi}^{(i)} (\varphi^{(i)} \psi^{(i)}(a)) \nonumber \\
&& - \bar{\varphi}^{(i)} \bar{\psi}^{(i)} (\varphi^{(i)}\psi^{(i)}(a)) \bar{\varphi}^{(i,j)} \bar{\Phi}^{(i,j)}(x) \| \nonumber \\
&& + 4 \eta  + 4 \eta^{\halb } \nonumber \\
& \stackrel{\eqref{wG8},\ref{multiplicative-domain}}{\le} &  \| \bar{\varphi}^{(i,j)} (\bar{\Phi}^{(i,j)}(x)  \bar{\psi}^{(i,j)} (\varphi^{(i)} \psi^{(i)}(a))) \nonumber \\
&& - \bar{\varphi}^{(i,j)} (\bar{\psi}^{(i,j)} (\varphi^{(i)}\psi^{(i)}(a))  \bar{\Phi}^{(i,j)}(x)) \| \nonumber \\
&& + 4 \eta  + 8 \eta^{\halb } \nonumber \\
& \stackrel{\eqref{wG9}}{\le} & 5\eta + 8 \eta^{\halb} \label{wG13}
\end{eqnarray}
for $x \in (M_{k})_{+}$ with $\|x\| \le 1$, $a \in \Fh$ and $i,j \in \{0, \ldots,m\}$.

Finally, we have
\begin{eqnarray}
\lefteqn{\tau\left(\sum_{i,j=0}^{m} \Phi^{(i,j)}(\be_{M_{k}}) h^{(i)}a\right)} \nonumber \\
& \stackrel{\eqref{wG3},\eqref{wG1},\ref{multiplicative-domain}}{\ge} & \tau\left(\sum_{i,j=0}^{m} \Phi^{(i,j)}(\be_{M_{k}}) \varphi^{(i)} \psi^{(i)}(a)\right) - (m+1)(2 \eta^{\halb} + \eta) \nonumber \\
& \stackrel{\eqref{wG8},\ref{multiplicative-domain},\eqref{wG7}}{\ge} & \tau\left(\sum_{i,j=0}^{m} \bar{\varphi}^{(i,j)}(\bar{\Phi}^{(i,j)}(\be_{M_{k}}) \bar{\psi}^{(i,j)} \varphi^{(i)} \psi^{(i)}(a))\right) \nonumber \\
&& - (m+1)(2 \eta^{\halb} + \eta) - (m+1)^{2}(2 \eta^{\halb} + \eta) \nonumber \\
& \stackrel{\eqref{wG10}}{\ge} & \frac{1}{4}  \cdot \tau\left(\sum_{i,j=0}^{m} \bar{\varphi}^{(i,j)} \bar{\psi}^{(i,j)} \varphi^{(i)} \psi^{(i)}(a)\right) \nonumber \\
&& - (m+1)(2 \eta^{\halb} + \eta) - (m+1)^{2}(2 \eta^{\halb} + \eta)  \nonumber \\
& \stackrel{\eqref{wG3},\eqref{wG8}}{\ge} & \frac{1}{4} \cdot \tau(a) - (m+1)(2 \eta^{\halb} + 2 \eta) - (m+1)^{2}(2 \eta^{\halb} + \eta) - \eta \nonumber \\
& \ge & \frac{1}{4} \cdot \tau(a) - 8(m+1)^{2} \eta^{\halb}  \label{wG14}
\end{eqnarray}
for all $a \in \Fh$ and $\tau \in T(A)$. 

Now if $(\Fh_{l})_{l \in \N}$ is an increasing sequence of compact sets exhausting $A_{+}$ and $(\eta_{l})_{l \in \N}$ is a decreasing null  sequence of strictly positive numbers, we can construct $h^{(i)}_{l}$ and $\Phi^{(i,j)}_{l}$ for $l \in \N$ following the above procedure. It is then clear from \eqref{wG11}, \eqref{wG12}, \eqref{wG13} and \eqref{wG14} that these give rise to elements $\tilde{h}^{(i)} \in (A_{\infty})_{+}$ and c.p.c.\ order zero maps $\widetilde{\Phi}^{(i,j)}:M_{k} \to A_{\infty}$ possessing properties (i) through (iv) of the proposition. 
\end{nproof}
\en

\bn
\label{wE}
\begin{nprop}
Given $m \in \N$, there are $K \in \N$, $0< \alpha$ and $0< \beta \le \frac{1}{2K}$ such that the following holds:

Let $A$ be separable, simple, unital and nonelementary with $\dr A \le m$; let $\Fh \subset A$ be a compact subset of nonzero positive elements of norm at most 1, and let $\eta, \varepsilon>0$ be given. 

Then, there are positive elements $d^{(1)}, \ldots, d^{(K)} \in A$ of norm at most 1 satisfying
\begin{enumerate}
\item[(a)] $\|[d^{(i)},a]\| < \eta$
\item[(b)] $\tau(\sum_{j=1}^{K} d^{(j)} a)  \ge \alpha \cdot \tau(a)   $
\item[(c)] $\tau((1-g_{0,\varepsilon})(d^{(i)})a) \ge (1-\beta) \tau(a)$
\end{enumerate}
for all $a \in \Fh$, $i \in \{1, \ldots, K\}$ and $\tau \in T(A)$.
\end{nprop}

\begin{nproof}
Define
\begin{equation}
\label{wE9}
K:=(m+1)^{2}, \, \beta:= \frac{1}{2K} \mbox{ and } \alpha:= \frac{\beta}{4}.
\end{equation}
We may apply Proposition~\ref{wG} with $K$ in place of $k$ to obtain elements 
\[
\tilde{h}^{(i)} \in (A_{\infty})_{+}
\]
and c.p.c.\ order zero maps
\[
\widetilde{\Phi}^{(i,j)}: M_{K} \to A_{\infty}
\]
for $i,j \in \{0,\ldots,m\}$ with the following properties:
\begin{enumerate}
\item $\sum_{i=0}^{m} \tilde{h}^{(i)} = \be_{A_{\infty}}$ (this implies that the $\tilde{h}^{(i)}$ have norm a t most 1)
\item $\tilde{h}^{(i)} \in A_{\infty} \cap A'$
\item $[\widetilde{\Phi}^{(i,j)}(M_{K}),\tilde{h}^{(i)}A] = 0$ (this implies $[\widetilde{\Phi}^{(i,j)}(M_{K}),\tilde{h}^{(i)}] = 0$ since $A$ is unital)
\item $\tau(\sum_{i,f=0}^{m} \widetilde{\Phi}^{(i,j)}(\be_{M_{K}}) \tilde{h}^{(i)}a) \ge \frac{1}{4} \tau(a)$ for any $\tau \in T_{\infty}(A)$ and $a \in A_{+}$.
\end{enumerate}
By (iii), as in \ref{order-zero-notation} we may define  c.p.c.\ order zero maps
\[
\hat{\Phi}^{(i,j)}: M_{K} \to A_{\infty}
\]
by setting 
\begin{equation}
\label{wE8}
\hat{\Phi}^{(i,j)}(y):= g_{0,\varepsilon}(\tilde{h}^{(i)}\widetilde{\Phi}^{(i,j)}(\be_{M_{K}})) \pi^{(i,j)}(y)
\end{equation}
for $y \in M_{K}$, where $\pi^{(i,j)}:M_{K} \to A''$ is the canonical supporting $*$-homomorphism for $\widetilde{\Phi}^{(i,j)}$. Note that
\begin{equation}
\label{wE6}
\hat{\Phi}^{(i,j)}(y)= g_{0,\varepsilon}((\tilde{h}^{i})^{\halb}\widetilde{\Phi}^{(i,j)}(y)(\tilde{h}^{i})^{\halb})
\end{equation}
if $y \in (M_{K})_{+}$.

It is clear from (ii), (iii), (iv) and \eqref{wE6}, that for $a \in A_{+}$ and $\tau \in T_{\infty}(A)$, 
\begin{equation}
\label{wE10}
\tilde{\tau}^{(i,j)}_{a}(\, . \,) := \tau(\widetilde{\Phi}^{(i,j)}(\, . \,) \tilde{h}^{(i)}a)
\end{equation}
and
\begin{equation}
\label{wE7}
\hat{\tau}^{(i,j)}_{a}(\, . \,) := \tau(\hat{\Phi}^{(i,j)}(\, . \,) a)
\end{equation}
define positive tracial functionals on $M_{K}$.

Set 
\begin{equation}
\label{wE5}
\tilde{d}^{(i,j)}:= \widetilde{\Phi}^{(i,j)}(q) \tilde{h}^{(i)},
\end{equation}
where $q \in M_{K}$ is a rank one projection; it is clear from (iii) that
\[
\tilde{d}^{(i,j)}= (\tilde{h}^{(i)})^{\halb} \widetilde{\Phi}^{(i,j)}(q) (\tilde{h}^{(i)})^{\halb}
\]
is positive of norm at most 1. We have
\begin{equation}
\label{wE3}
\tau(g_{0,\varepsilon}(\tilde{d}^{(i,j)})a) \stackrel{\eqref{wE6},\eqref{wE7},\eqref{wE8}}{=} \hat{\tau}^{(i,j)}_{a}(q) = \frac{1}{K}\cdot  \hat{\tau}^{(i,j)}_{a}(\be_{M_{K}}) \stackrel{\eqref{wE7}}{\le} \frac{1}{K} \cdot \tau(a) \stackrel{\eqref{wE9}}{=} \frac{\beta}{2} \cdot \tau(a)
\end{equation}
for any $\tau \in T_{\infty}(A)$ and $a \in A_{+}$. Moreover,
\begin{eqnarray}
\tau\left(\sum_{i,j=0}^{m} \tilde{d}^{(i,j)}a\right) & \stackrel{\eqref{wE5}, \eqref{wE10}}{=} & \sum_{i,j=0}^{m} \tilde{\tau}_{a}^{(i,j)}(q) \nonumber \\
& = & \frac{1}{K} \cdot \sum_{i,j=0}^{m} \tilde{\tau}^{(i,j)}_{a}(\be_{M_{K}}) \nonumber \\
& \stackrel{\eqref{wE10}}{=} & \frac{1}{K} \cdot \tau\left(\sum_{i,j=0}^{m} \widetilde{\Phi}^{(i,j)}(\be_{M_{K}})\tilde{h}^{(i)}a\right) \nonumber \\
& \stackrel{\mathrm{(iv)}}{\ge} & \frac{1}{4K} \cdot \tau(a) \nonumber \\
& \stackrel{\eqref{wE9}}{=} & \frac{\beta}{2} \cdot \tau(a) \label{wE4}
\end{eqnarray}
for $i,j \in \{0, \ldots,m\}$, $a \in A_{+}$ and $\tau \in T_{\infty}(A)$. 

Now lift each $\tilde{d}^{(i,j)}$ to a positive sequence 
\begin{equation}
\label{wE14}
(\tilde{d}^{(i,j)}_{l})_{l \in \N} \in \prod_{\N} A
\end{equation} 
of norm at most 1. 

We claim there is $l_{0} \in \N$ such that 
\begin{equation}
\label{wE1}
\tau(g_{0,\varepsilon}(\tilde{d}^{(i,j)}_{l})a) \le \beta \cdot \tau(a)
\end{equation}
and
\begin{equation}
\label{wE2}
\tau\left(\sum_{i,j=0}^{m} \tilde{d}^{(i,j)}_{l}a\right) \ge \alpha \cdot \tau(a)
\end{equation}
for all $\tau \in T(A)$, $a \in \Fh$ and $l \ge l_{0}$. 

To prove the claim, in order to obtain a contradiction let us assume that no such $l_{0}$ exists. Then either \eqref{wE1} or \eqref{wE2} fails, i.e., there are sequences 
\[
(a_{n})_{n \in \N} \subset \Fh , \,  (\tau_{n})_{n \in \N} \subset (T(A))\]
and a subsequence $(l_{n})_{n \in \N}$ of $(n)_{n \in \N}$ 
such that 
\[
\tau_{l_{n}}(g_{0,\varepsilon}(\tilde{d}_{l_{n}}^{(i,j)})a_{l_{n}}) > \beta \cdot \tau_{l_{n}}(a_{l_{n}})
\]
or
\[
\tau_{l_{n}}\left(\sum_{i,j=0}^{m} \tilde{d}^{(i,j)}_{l_{n}}a_{n}\right) < \alpha \cdot \tau_{l_{n}}(a_{l_{n}})
\]
for $n \in \N$. Now since $\Fh$ is compact, by passing to a subsequence we may even assume that there is $\bar{a} \in \Fh$ such that
\[
a_{l_{n}} \to \bar{a}.
\]
By Proposition~\ref{trace-min}, 
\begin{equation}
\label{wE15}
\zeta:= \min\{\tau(a) \mid \tau \in T(A), \, a \in \Fh\}
\end{equation}
exists and is strictly positive. Choose 
\begin{equation}
\label{wE16}
0< \bar{\varepsilon} < \frac{\alpha \zeta}{4}.
\end{equation}
But then, there is $n_{0} \in \N$ such that 
\[
\|a_{l_{n}} - \bar{a}\| < \bar{\varepsilon}
\]
for $n \ge n_{0}$; we may assume without loss of generality that $n_{0}=0$. We obtain
\begin{eqnarray}
\tau_{l_{n}}(g_{0,\varepsilon}(\tilde{d}^{(i,j)}_{l_{n}})\bar{a}) & \ge & \tau_{l_{n}}(g_{0,\varepsilon}(\tilde{d}^{(i,j)}_{l_{n}})a_{l_{n}}) - \bar{\varepsilon} \nonumber \\
& > & \beta \cdot \tau_{l_{n}}(a_{l_{n}}) - \bar{\varepsilon} \nonumber \\
& \ge &  \beta \cdot \tau_{l_{n}}(\bar{a}) - 2 \bar{\varepsilon} \label{wE11}
\end{eqnarray}
or 
\begin{eqnarray}
\tau_{l_{n}}\left(\sum_{i,j=0}^{m} \tilde{d}^{(i,j)}_{l_{n}} \bar{a}\right) & \le & \tau_{l_{n}}\left(\sum_{i,j=0}^{m} \tilde{d}^{(i,j)}_{l_{n}} a_{l_{n}}\right) + \bar{\varepsilon} \nonumber \\
& < & \alpha \cdot \tau_{l_{n}}(a_{l_{n}}) + \bar{\varepsilon} \nonumber \\
& \le & \alpha \cdot \tau_{l_{n}} (\bar{a}) + 2 \bar{\varepsilon} \label{wE12}
\end{eqnarray}
for $n \in \N$. Let $\omega \in \beta \N \setminus \N$ be a free ultrafilter along the subsequence $(l_{n})_{n \in \N}$ of $(n)_{n \in \N}$, and let $\bar{\tau} \in T_{\infty}(A)$ denote the trace given by
\begin{equation}
\label{wE13}
\bar{\tau}([(b_{n})_{n \in \N}]) := \lim_{n \to \omega} \tau_{n}(b_{n}).
\end{equation}
It is clear from \eqref{wE11} or \eqref{wE12}, respectively, that
\begin{eqnarray*}
\bar{\tau}(g_{0,\varepsilon}(\tilde{d}^{(i,j)})\bar{a}) & \stackrel{\eqref{wE11},\eqref{wE13},\eqref{wE14}}{\ge} & \beta \cdot \bar{\tau}(\bar{a}) - 2 \bar{\varepsilon} \\
& \stackrel{\eqref{wE15},\eqref{wE16},\eqref{wE9}}{\ge} & \frac{3}{4} \beta \cdot \bar{\tau}(\bar{a})
\end{eqnarray*}
or 
\begin{eqnarray*}
\bar{\tau}\left(\sum_{i,j=0}^{m} \tilde{d}^{(i,j)} \bar{a}\right) & \stackrel{\eqref{wE12},\eqref{wE13},\eqref{wE14}}{\le} & \alpha \cdot \bar{\tau}(\bar{a}) + 2 \bar{\varepsilon} \\
& \stackrel{\eqref{wE15},\eqref{wE16}}{\le} & \frac{3}{2} \alpha \cdot \bar{\tau}(\bar{a})\\
& \stackrel{\eqref{wE9}}{=}& \frac{3}{8} \beta \cdot \bar{\tau}(\bar{a}) .
\end{eqnarray*}
We thus obtain a contradiction to \eqref{wE3} or \eqref{wE4}, and  so $l_{0}$ as above indeed exists. 

By (ii), (iii) and \eqref{wE5}, there is $l_{1} \ge l_{0}$ such that 
\[
\|[\tilde{d}^{(i,j)}_{l_{1}},a]\| < \eta
\]
for $i,j \in \{0, \ldots, m\}$ and $a \in \Fh$. 

Let $d^{(1)}, \ldots, d^{(K)}$ be an enumeration of $\{\tilde{d}^{(i,j)}_{l_{1}} \mid i,j=0, \ldots,m\}$; then, by \eqref{wE1} and \eqref{wE2}, 
\[
\tau(g_{0,\varepsilon}(d^{(i)})a) \le \beta \cdot \tau(a) \mbox{ for } i=1, \ldots,K
\]
and
\[
\tau\left(\sum_{j=1}^{K} d^{(i)} a \right) \ge \alpha \cdot \tau(a)
\]
for all $\tau \in T(A)$ and $a \in \Fh$. So, the $d^{(i)}$ satisfy \ref{wE}(a), (b) and (c).
\end{nproof}
\en

\bn
\label{wB}
\begin{nprop}
Given $m,k \in \N$, there is $\omega_{k}>0$ such that the following holds:

Let $A$ be separable, simple, unital and nonelementary with $\dr A \le m$; let $\theta > 0$ and let $\Fh \subset A$ be a compact subset of nonzero positive elements of norm at most 1. 

Then, there are pairwise orthogonal positive elements $d_{1}, \ldots, d_{k} \in A$ of norm at most 1 such that
\[
\|[d_{i},f]\| < \theta
\]
and
\[
\tau(d_{i}f) \ge \omega_{k} \cdot \tau(f)
\]
for all $i \in \{1, \ldots,k\}$, $f \in \Fh$ and $\tau \in T(A)$.
\end{nprop}

\begin{nproof}
We shall use induction to prove the proposition when $k$ is of the form $2^{l}$, $l \in \N$; from here it will be easy to deduce the case of general $k$. So, let us first prove the assertion for $k=2$. Choose $\alpha$, $\beta$ and $K$ as in Proposition~\ref{wE}, and choose
\begin{equation}
\label{wB11}
0<\omega_{2} < \frac{\alpha}{4}, 1-2K\beta.
\end{equation}
Let $A$, $\Fh$ and $\theta$ as in the proposition be given. Set
\begin{equation}
\label{wB4}
\zeta:= \min\{\tau(a) \mid a\in \Fh, \, \tau\in T(A)\};
\end{equation}
$\zeta$ exists and is nonzero by Proposition~\ref{trace-min}. Choose 
\begin{equation}
\label{wB6}
0< \varepsilon< \frac{\zeta \alpha}{8K}.
\end{equation}

For $l\in \N$, apply Proposition~\ref{wE} with $\frac{1}{l}$ in place of $\eta$ to obtain positive elements $d^{(1)}_{l}, \ldots, d^{(K)}_{l} \in A$ of norm at most 1 satisfying
\begin{equation}
\label{wB3}
\|[d_{l}^{(i)},a]\| < \frac{1}{l}
\end{equation}
\[
\tau\left(\left(\sum_{j=1}^{K} d_{l}^{(j)}\right)a\right) \ge \alpha \cdot \tau(a)
\]
and
\[
\tau((1-g_{0,\varepsilon})(d^{(i)}_{l})a) \ge (1-\beta)\tau(a)
\]
for all $a \in \Fh$, $i \in \{1, \ldots,K\}$, $ \tau \in T(A)$ and $l \in \N$. 

We obtain
\[
\bar{d}^{(i)}:= [(d_{l}^{(i)})_{l \in \N}] \in A_{\infty} 
\]
satisfying 
\[
\bar{d}^{(i)} \in A_{\infty} \cap A',
\]
\begin{equation}
\label{wB5}
\tau_{\omega}\left(\left(\sum_{i=1}^{K} \bar{d}^{(i)}\right)a\right) \ge \alpha \cdot \tau(a)
\end{equation}
and
\begin{equation}
\label{wB7}
\tau_{\omega}((1-g_{0,\varepsilon})(\bar{d}^{(i)})a) \ge (1-\beta) \cdot \tau(a)
\end{equation}
for all $a \in \Fh$, $i \in \{1, \ldots,K\}$, $\tau \in T(A)$ and any free ultrafilter $\omega \in \beta \N \setminus \N$, where $\tau_{\omega} \in T_{\infty}(A)$ denotes the tracial state given by 
\[
\tau_{\omega}([(a_{l})_{l\in \N}]) := \lim_{l \to \omega} \tau(a_{l}).
\]
Set 
\[
C:= C^{*}(\bar{d}^{(1)}, \ldots, \bar{d}^{(K)}, \be_{A_{\infty}}) \subset A_{\infty} \cap A'.
\]
Note that, for $a \in A_{+}$, 
\[
\tau_{\omega,a}(c) := \tau_{\omega}(ca), \, c \in C,
\]
defines a positive tracial functional on $C$. We obtain
\begin{eqnarray}
\lefteqn{\frac{\alpha}{2} \cdot \tau(a)} \label{wB9} \\
& \le & \alpha \cdot \tau(a) - \frac{\alpha \zeta}{2} \nonumber \\
& \stackrel{\eqref{wB4}}{\le} & \tau_{\omega,a}\left(\sum_{j=1}^{K} \bar{d}^{(j)}\right) - \frac{\alpha \zeta}{2} \nonumber \\
& \stackrel{\eqref{wB5}}{\le} & \tau_{\omega,a}\left(\sum_{j=1}^{K} f_{\varepsilon,2\varepsilon}(\bar{d}^{(j)})\right) + K \varepsilon - \frac{\alpha \zeta}{2} \nonumber \\
& \le & \tau_{\omega,a}\left(\sum_{j=1}^{K} f_{\varepsilon,2\varepsilon}(\bar{d}^{(j)})\right) \label{wB8} \\
& \stackrel{\eqref{wB6}}{\le} & \tau_{\omega,a}\left(g_{0,\varepsilon}\left(\sum_{j=1}^{K} f_{\varepsilon,2\varepsilon}(\bar{d}^{(j)})\right)\right) \label{wB1} \\
& \le & d_{\tau_{\omega,a}}\left(\sum_{j=1}^{K} f_{\varepsilon,2\varepsilon}(\bar{d}^{(j)})\right) \nonumber \\
& \stackrel{}{\le} & \sum_{j=1}^{K} d_{\tau_{\omega ,a}}(f_{\varepsilon,2\varepsilon}(\bar{d}^{(j)})) \nonumber \\
& \le & \sum_{j=1}^{K} \tau_{\omega,a}(g_{0,\varepsilon}(\bar{d}^{(j)})) \nonumber \\
& \stackrel{\eqref{wB7}}{\le} & K \beta \cdot \tau(a) \label{wB2} 
\end{eqnarray}
for any $a \in \Fh$, $\tau \in T(A)$ and $\omega \in \beta\N \setminus \N$. Here, for the third last estimate we have used the fact that $d(a+b) \le d(a) + d(b)$ for any dimension function $d$ and positive elements $a,b$. 

Since these estimates hold for any free ultrafilter $\omega$ on $\N$, and since $\Fh$ is compact, there is  $L_{0} \in \N$ such that
\begin{equation}
\label{wB12}
\frac{\alpha}{2} \cdot \tau(a) - \frac{\alpha \zeta}{8} \stackrel{\eqref{wB9},\eqref{wB8}}{\le} \tau\left(\left(\sum_{j=1}^{K} f_{\varepsilon,2\varepsilon}(d_{l}^{(j)})\right)a\right)
\end{equation}
and
\begin{equation}
\label{wB10}
\tau\left(g_{0,\varepsilon}\left(\sum_{j=1}^{K} f_{\varepsilon,2\varepsilon}(d_{l}^{(j)})\right)a\right) \stackrel{\eqref{wB1},\eqref{wB2}}{\le} K \beta \cdot \tau(a) + K \beta \zeta \stackrel{\eqref{wB4}}{\le} 2 K \beta \cdot \tau(a)
\end{equation}
for any $a \in \Fh$ and $\tau \in T(A)$, if $l \ge L_{0}$. 

We then obtain
\begin{equation}
\label{wB14}
\tau\left(f_{\varepsilon,2\varepsilon}\left(\sum_{j=1}^{K} f_{\varepsilon,2\varepsilon}(d_{l}^{(j)})\right)a\right) \stackrel{\eqref{wB12}}{\ge} \frac{\alpha}{2} \cdot \tau(a) - \frac{2 \alpha \zeta}{8} \ge \frac{\alpha}{4} \cdot \tau(a) \ge \omega_{2} \cdot \tau(a)
\end{equation}
and
\begin{eqnarray}
\tau((1-g_{0,\varepsilon})\left(\sum_{j=1}^{K}f_{\varepsilon,2\varepsilon}(d_{l}^{(j)}))a\right) & \stackrel{\eqref{wB10},\eqref{wB6}}{\ge} & 1- 2K \beta \cdot \tau(a) \nonumber \\
& \ge & (1-2K\beta)\cdot \tau(a) \nonumber \\
& \stackrel{\eqref{wB11}}{\ge} & \omega_{2} \cdot \tau(a) \label{wB15}
\end{eqnarray}
if $a \in \Fh$ and $l \ge L_{0}$.

By \eqref{wB3}, and since $\Fh$ is compact, there is $L_{1} \in \N$ such that
\begin{equation}
\label{wB13}
\left\|\left[f_{\varepsilon,2\varepsilon}\left(\sum_{j=1}^{K}f_{\varepsilon,2\varepsilon}(d_{l}^{(j)})\right),a\right]\right\|, \, \left\|\left[g_{0,\varepsilon}\left(\sum_{j=1}^{K}f_{\varepsilon,2\varepsilon}(d_{l}^{(j)})\right),a\right]\right\| < \theta
\end{equation}
if $a \in \Fh$ and $l\ge L_{1}$. Set
\[
L:= \max\{L_{0},L_{1}\}
\]
and define
\[
d_{1}:= f_{\varepsilon,2\varepsilon}\left(\sum_{j=1}^{K} f_{\varepsilon,2\varepsilon}(d^{(j)}_{L})\right)
\]
and
\[
d_{2}:= (1-g_{0,\varepsilon})\left(\sum_{j=1}^{K} f_{\varepsilon,2\varepsilon}(d^{(j)}_{L})\right).
\]
Then, $d_{1}$ and $d_{2}$ are orthogonal, positive, of norm at most 1 and satisfy
\[
\|[d_{i},a]\| \stackrel{\eqref{wB13}}{<} \theta
\]
and
\[
\tau(d_{i}a) \stackrel{\eqref{wB14},\eqref{wB15}}{\ge} \omega_{2} \cdot \tau(a)
\]
for $i=1,2$, $a \in \Fh$ and $\tau \in T(A)$. We have thus verified the anchor step of our induction, i.e., the conclusion of the proposition in the case $k$$=$$2$.

For the induction step, suppose we have shown the proposition for $k=2^{l}$ for some $l \in \N$ and have some $0<\omega_{2^{l}}$. Set
\begin{equation}
\label{wB28}
\omega_{2^{l+1}}:= \frac{1}{2} \omega_{2} \omega_{2^{l}}.
\end{equation}
Again, let $A$, $\Fh$ and $\theta$ as in the proposition be given. Set 
\begin{equation}
\label{wB27}
\zeta:= \min\{\tau(a) \mid a \in \Fh,\, \tau \in T(A)\};
\end{equation}
as, before, $\zeta$ exists and is strictly positive  by Proposition~\ref{trace-min}. 

Apply the case $k=2^{l}$ of the proposition to obtain pairwise orthogonal positive elements $d_{1},\ldots,d_{2^{l}} \in A$ of norm at most 1 such that 
\begin{equation}
\label{wB20}
\|[d_{i},a]\| < \frac{\theta}{6}
\end{equation}
and 
\begin{equation}
\label{wB26}
\tau(d_{i}a) \ge \omega_{2^{l}} \cdot \tau(a)
\end{equation}
for all $i \in \{1, \ldots,2^{l}\}$, $a \in \Fh$ and $\tau \in T(A)$. Define 
\begin{equation}
\label{wB19}
\bar{\Fh}:= \Fh \cup \{d_{i}^{\halb} \mid i=1,\ldots,2^{l}\},
\end{equation}
then $\bar{\Fh}$ is compact and consists of nonzero positive elements of norm at most 1. Choose
\begin{equation}
\label{wB22}
0< \delta \le \omega_{2^{l+1}} \zeta , \, \frac{\theta}{6}.
\end{equation}
Choose 
\begin{equation}
\label{wB21}
0<\theta'< \frac{\theta}{6}
\end{equation} 
such that, whenever $b_{1},\ldots, b_{2^{l+1}} \in A$ are positive elements of norm at most 1 satisfying $\|b_{i}b_{i'}\| \le 2 \theta'$ if $i\neq i'$, then there are pairwise orthogonal $\tilde{b}_{i} \in A_{+}$ of norm at most 1 such that 
\begin{equation}
\label{wB17}
\|\tilde{b}_{i} - b_{i}\| \le \delta
\end{equation} 
for $i=1, \ldots,2^{l+1}$; this is possible by \cite[Proposition~2.5]{KirWinter:dr}.

Apply the proposition (for $k=2$) to obtain orthogonal positive elements $\hat{d}_{1},\hat{d}_{2} \in A$ of norm at most 1 such that
\begin{equation}
\label{wB16}
\|[\hat{d}_{j},a]\| < \theta'
\end{equation}
and
\begin{equation}
\label{wB25}
\tau(\hat{d}_{j}a) \ge \omega_{2} \cdot \tau(a)
\end{equation}
for $j=1,2$, $a \in \bar{\Fh}$ and $\tau \in T(A)$. 

Note that by \eqref{wB16} and by pairwise orthogonality of the $\hat{d}_{j}$ and the $d_{i}'$ we have 
\[
\|d_{i}^{\halb} \hat{d}_{j}d_{i}^{\halb} d_{i'}^{\halb} \hat{d}_{j'}d_{i'}^{\halb}\| \le  2 \theta' 
\]
if $i \neq i' $ or $j \neq j'$, so by our choice of $\theta'$ there are pairwise orthogonal positive elements $\tilde{d}_{i,j} \in A$ of norm at most 1 such that
\begin{equation}
\label{wB18}
\| \tilde{d}_{i,j} - d_{i}^{\halb} \hat{d}_{j} d_{i}^{\halb} \| \le \delta
\end{equation}
for $i \in \{1, \ldots, 2^{l}\}$, $j \in \{1,2\}$; cf.\ \eqref{wB17}. We now have
\begin{eqnarray}
\|[\tilde{d}_{i,j},a]\| & \stackrel{\eqref{wB18}}{\le} & \| [ d_{i}^{\halb} \hat{d}_{j} d_{i}^{\halb},a]\| + 2 \delta \nonumber \\
& \stackrel{\eqref{wB16},\eqref{wB19}}{\le} & \|[d_{i}\hat{d}_{j},a]\| + 2 \delta + 2 \theta' \nonumber \\
& \stackrel{\eqref{wB16},\eqref{wB20}}{<} & \theta' + \frac{\theta}{6} + 2 \delta + 2 \theta' \nonumber \\
& \stackrel{\eqref{wB22},\eqref{wB21}}{\le} & \theta \label{wB23}
\end{eqnarray}
for $a \in \Fh$, $i \in \{1,\ldots, 2^{l}\}$ and $j \in \{1,2\}$.

Moreover, 
\begin{eqnarray}
\tau(\tilde{d}_{i,j}a) & \stackrel{\eqref{wB18}}{\ge} & \tau(d_{i}^{\halb} \hat{d}_{j} d_{i}^{\halb} a) - \delta \nonumber \\
& = & \tau(d_{i}^{\halb} \hat{d}_{j} d_{i}^{\halb} a d_{i}^{\halb}) - \delta \nonumber \\
& \stackrel{\eqref{wB25},\eqref{wB19}}{\ge} & \omega_{2} \cdot \tau(d_{i}^{\halb} a d_{i}^{\halb}) - \delta \nonumber \\
& = & \omega_{2} \cdot \tau(d_{i}a) - \delta \nonumber \\
& \stackrel{\eqref{wB26}}{\ge} & \omega_{2} \omega_{2^{l}} \cdot \tau(a) - \delta \nonumber \\
& \stackrel{\eqref{wB22}}{\ge} & \omega_{2} \omega_{2^{l}} \cdot \tau(a) - \frac{\omega_{2} \omega_{2^{l}}\zeta}{2} \nonumber \\
& \stackrel{\eqref{wB27}}{\ge} & \frac{1}{2} \omega_{2} \omega_{2^{l}} \cdot \tau(a) \nonumber \\
& \stackrel{\eqref{wB28}}{=} & \omega_{2^{l+1}} \cdot \tau(a)  \label{wB24}
\end{eqnarray}
for $a \in \Fh$, $i \in \{1, \ldots,2^{l}\}$ and $j \in \{1,2\}$. Now \eqref{wB23} and \eqref{wB24} show that the elements $\tilde{d}_{ij}$ satisfy the assertion of the proposition for $k=2^{l+1}$.

We have thus inductively verified the proposition for $k=2^{l}$, $l \in \N$. For arbitrary $ k\in \N$, choose $l\in \N$ such that $2^{l} \ge k$ and set $\omega_{k}:= \omega_{2^{l}}$. Then the assertion for $2^{l}$ will obviously yield the assertion for $k$. 
\end{nproof}
\en

\bn
\label{wC}
\begin{nprop}
Given $m,k \in \N$, there is $\bar{\beta} > 0$ such that the following holds:

Let $A$ be separable, simple, unital and nonelementary with $\dr A \le m$; let $ f \in A$ be positive of norm at most 1 and $\bar{\varepsilon}>0$ be given. 

Then, there is a c.p.c.\ order zero map
\[
\bar{\Phi}: M_{k} \to \her(f) \subset A
\]
such that
\[
\|[\bar{\Phi}(x),f]\| \le \bar{\varepsilon} \cdot \|x\|
\]
for $x \in M_{k}$ and 
\[
\tau(\bar{\Phi}(\be_{M_{k}})f) \ge \bar{\beta} \cdot \tau(f).
\]
\end{nprop}

\begin{nproof}
Set $$\bar{\beta}:= \frac{1}{4} \omega_{m+1},$$ where $\omega_{m+1}$ comes from Proposition~\ref{wB}. 

Let $A$, $f$ and $\bar{\varepsilon}$ be as in the proposition. We may assume $f$ to be nonzero, for otherwise we can just take $\bar{\Phi}$ to be the zero map. Now by Proposition~\ref{trace-min}, 
\[
\zeta:= \min\{\tau(f) \mid \tau \in T(A) \}
\]
exists and is nonzero. Choose $$0<\eta <1$$ such that 
\begin{equation}
\label{wC17}
6 (m+1) \eta < \frac{1}{4} \omega_{m+1} \zeta , \, \bar{\varepsilon}.
\end{equation}
Choose 
\begin{equation}
\label{wC19}
0< \bar{\eta} <\eta 
\end{equation} such that 
\begin{equation}
\label{wC4}
7 \bar{\eta}^{\halb} <\eta
\end{equation}
and such that, if $\|[h,f]\| \le 7 \bar{\eta}^{\halb}$ for some positive $h$ of norm at most 1, then 
\begin{equation}
\label{wC6}
\|[h^{\halb},f]\|\le \eta; 
\end{equation}
this is possible by Proposition~\ref{I}.

By \cite[Proposition~3.8]{KirWinter:dr}, we have $\dr \her(f) = \dr A$. Let 
\[
(F=F^{(0)} \oplus \ldots \oplus F^{(m)}, \psi, \varphi)
\]
be an $m$-decomposable c.p.c.\ approximation as in Proposition~\ref{wD}, with $\bar{\eta}$ in place of $\tilde{\eta}$. Note that, by Lemma~\ref{multiplicative-domain}, 
\begin{equation}
\label{wC2}
\|\varphi(y\psi(f)) - \varphi(y) \varphi \psi(f) \| \le 2 \bar{\eta}^{\halb} \cdot \|y\|
\end{equation}
for any $y \in F$.

Set 
\[
f^{(i)}:= \psi^{(i)}(f) \in F^{(i)}.
\]
By Proposition~\ref{wD}, there are $*$-homomorphisms
\[
\tilde{\Phi}^{(i)}: M_{k} \to F^{(i)}
\]
such that
\begin{equation}
\label{wC5}
\|[\tilde{\Phi}^{(i)}(x),f^{(i)}]\| \le \bar{\eta} \|x\|
\end{equation}
for $x \in F^{(i)}$ and such that, for any positive tracial functional $\tau^{(i)}$ on $F^{(i)}$, we have
\begin{equation}
\label{wC1}
\tau^{(i)}(\tilde{\Phi}^{(i)}(\be_{M_{k}})f^{(i)}) \ge \frac{1}{2} \cdot \tau^{(i)}(f^{(i)}).
\end{equation}
Define c.p.c.\ order zero maps 
\begin{equation}
\label{wC7}
\hat{\Phi}^{(i)}: M_{k} \to \her(f)
\end{equation}
by
\begin{equation}
\label{wC3}
\hat{\Phi}^{(i)}(x):= \varphi^{(i)}(\tilde{\Phi}^{(i)}(x)).
\end{equation}
We have
\begin{eqnarray}
\lefteqn{\tau(\hat{\Phi}^{(i)}(\be_{M_{k}})f)} \nonumber \\
& \ge & \tau(\hat{\Phi}^{(i)}(\be_{M_{k}}) \varphi \psi(f)) - \bar{\eta} \nonumber \\
& \stackrel{\eqref{wC2},\eqref{wC3}}{\ge} & \tau(\varphi^{(i)}(\tilde{\Phi}^{(i)}(\be_{M_{k}})\psi^{(i)}(f))) - \bar{\eta} - 2 \bar{\eta}^{\halb} \nonumber \\
& \stackrel{\eqref{wC4}}{\ge} & \frac{1}{2} \cdot \tau(\varphi^{(i)}\psi^{(i)}(f)) - 3 \eta \label{wC16}
\end{eqnarray}
for any $\tau \in T(A)$ and $i \in \{ 0, \ldots,m\}$, where for the last estimate we have also used \eqref{wC1} and the fact that $\tau \circ \varphi^{(i)}$ is a tracial functional on $F^{(i)}$. 

Moreover, we have 
\begin{eqnarray}
\|[\hat{\Phi}^{(i)}(x),f]\| & \le & \|[\hat{\Phi}^{(i)}(x), \varphi \psi(f)]\| + 2 \bar{\eta} \|x\| \nonumber \\
& \stackrel{\eqref{wC2},\eqref{wC3}}{\le} & \|\varphi^{(i)}([\tilde{\Phi}^{(i)}(x), f^{(i)}])\| + (2 \bar{\eta} + 4 \bar{\eta}^{\halb}) \|x\| \nonumber \\
& \le & 7 \bar{\eta}^{\halb} \|x\|  \label{wC18}
\end{eqnarray}
for $x \in M_{k}$ and $ i \in \{0, \ldots, m\}$. Note that by our choice of $\bar{\eta}$, this implies 
\begin{equation}
\label{wC13}
\|[\hat{\Phi}^{(i)}(\be_{M_{k}})^{\halb},f]\| \le \eta.
\end{equation}
Next, choose 
\begin{equation}
\label{wC15} 
0< \delta < \eta
\end{equation} 
such that the following holds: If $\bar{d}^{(0)}, \ldots, \bar{d}^{(m)} \in \her(f)$ are positive elements of norm at most 1, satisfying
\[
\|\bar{d}^{(i)} \bar{d}^{(i')} \| \le \delta \mbox{ if } i \neq i' \in \{0, \ldots, m\}
\]
and
\[
\|[\bar{d}^{(i)}, \hat{\Phi}^{(i)}(x)]\| \le \delta  \|x\|, 
\]
then there are pairwise orthogonal c.p.c.\ order zero maps 
\[
\bar{\Phi}^{(i)}: M_{k} \to \her(f)
\]
such that 
\begin{equation}
\label{wC11}
\| \bar{\Phi}^{(i)}(x) - \bar{d}^{(i)} \hat{\Phi}^{(i)}(x) \bar{d}^{(i)} \| \le \eta \|x\|
\end{equation}
for $i = 0, \ldots, m$ and $x \in M_{k}$. Such a $\delta$ exists by \cite[Proposition~2.5]{KirWinter:dr}.

Choose 
\[
0 < \gamma < \eta
\]
such that 
\begin{equation}
\label{wC8}
\| g_{0,\gamma}(f) \hat{\Phi}^{(i)}(x) -  \hat{\Phi}^{(i)}(x) \| \le \frac{\delta}{36} \|x\|
\end{equation}
for $i=0, \ldots, m$ and $x \in M_{k}$; this is possible since the unit ball of $M_{k}$ is compact and the image of $\hat{\Phi}^{(i)}$ is in $\her(f)$; cf.\ \eqref{wC7}. 

Set
\[
\Fh:=  \{f, \, g_{0,\gamma}(f), \, \hat{\Phi}^{(i)}(x), \, \hat{\Phi}^{(i)}(\be_{M_{k}})^{\halb} f \hat{\Phi}^{(i)}(\be_{M_{k}})^{\halb} \mid 0 \le x \le \be_{M_{k}}, \, i \in \{0, \ldots, m\}\},
\]
then $\Fh \subset A$ is a compact subset of positive elements of norm at most 1. 

Choose $$0 < \theta < \delta $$ such that, if $\|[d,b]\| < \theta$, then 
\begin{equation}
\label{wC10}
\|[d^{\halb},b]\| < \frac{\delta}{36}
\end{equation} 
whenever $d,b \in A$ are positive elements of norm at most 1; cf.\ Proposition~\ref{I}. 

From Proposition~\ref{wB} we now obtain pairwise orthogonal positive elements $$d^{(0)}, \ldots, d^{(m)} \in A$$ such that
\begin{equation}
\label{wC9}
\|[d^{(i)},b]\| < \theta
\end{equation}
and
\begin{equation}
\label{wC14}
\tau(d^{(i)}b) \ge \omega \cdot \tau(b)
\end{equation}  
for $i \in \{0, \ldots, m\} $, $ b \in \Fh$ and $\tau \in T(A)$.

For each $i$, we then have
\begin{eqnarray*}
\lefteqn{\|[g_{0,\gamma}(f)(d^{(i)})^{\halb} g_{0,\gamma}(f), \hat{\Phi}^{(i)}(x)]\| } \\
 & \stackrel{\eqref{wC8}}{\le} &  \|  g_{0,\gamma}(f)(d^{(i)})^{\halb}  \hat{\Phi}^{(i)}(x) - \hat{\Phi}^{(i)}(x)   (d^{(i)})^{\halb}g_{0,\gamma}(f)\| + 2 \frac{\delta}{36} \|x\| \\
 & \stackrel{\eqref{wC9}, \eqref{wC10}, \eqref{wC8}}{\le} &  \|  (d^{(i)})^{\halb}  \hat{\Phi}^{(i)}(x) - \hat{\Phi}^{(i)}(x)   (d^{(i)})^{\halb}\| + 8 \frac{\delta}{36} \|x\| \\
 & \stackrel{\eqref{wC9}, \eqref{wC10}}{\le} &  9 \frac{\delta}{36} \|x\| 
\end{eqnarray*}
for $x \in (M_{k})_{+}$, from which we obtain 
\[
\|[g_{0,\gamma}(f)(d^{(i)})^{\halb} g_{0,\gamma}(f), \hat{\Phi}^{(i)}(x)]\| \le \delta \|x\|
\] 
for all $x \in M_{k}$. Moreover, if $i \neq i' \in \{0, \ldots,m\}$, 
\begin{eqnarray*}
\| g_{0, \gamma}(f)(d^{(i)})^{\halb} g_{0,\gamma}(f) (d^{(i')})^{\halb} g_{0, \gamma}(f) \| & \stackrel{\eqref{wC9},\eqref{wC10}}{\le} & \| g_{0, \gamma}(f)^{2} (d^{(i)})^{\halb} (d^{i'})^{\halb} g_{0, \gamma}(f)^{2} \| \\
&& + 2 \frac{\delta}{36} \\
& = & 2  \frac{\delta}{36} 
\end{eqnarray*}
(using that $d^{i} \perp d^{i'}$), whence by our choice of $\delta$ there are pairwise orthogonal c.p.c.\ order zero maps
\[
\bar{\Phi}^{(i)}:M_{k} \to \her(f)
\]
such that
\begin{equation}
\label{wC12}
\| \bar{\Phi}^{(i)}(x) - g_{0,\gamma}(f) (d^{(i)})^{\halb} g_{0,\gamma}(f) \hat{\Phi}^{(i)}(x) g_{0,\gamma}(f) (d^{(i)})^{\halb} g_{0,\gamma}(f) \| \le \eta \|x\|,
\end{equation}
cf.\ \eqref{wC11}. We may define a c.p.c.\ order zero map
\[
\bar{\Phi}:M_{k} \to \her(f)
\]
by
\[
\bar{\Phi}(x):= \bigoplus_{i=0}^{m} \bar{\Phi}^{(i)}(x).
\]
We have
\begin{eqnarray*}
\lefteqn{\tau(\bar{\Phi}^{(i)}(\be_{M_{k}})f)}\\
& \stackrel{\eqref{wC12}, \eqref{wC8},\eqref{wC9},\eqref{wC10}}{\ge} & \tau((d^{(i)})^{\halb} \hat{\Phi}^{(i)}(\be_{M_{k}}) (d^{(i)})^{\halb} f) - \frac{6}{36}\delta  \\
& \stackrel{\eqref{wC9},\eqref{wC10}}{\ge} & \tau(d^{(i)} \hat{\Phi}^{(i)}(\be_{M_{k}}) f) - \frac{7}{36}\delta  \\
& \stackrel{\eqref{wC13}}{\ge} & \tau(d^{(i)} \hat{\Phi}^{(i)}(\be_{M_{k}})^{\halb} f  \hat{\Phi}^{(i)}(\be_{M_{k}})^{\halb}) - \eta - \frac{7}{36}\delta  \\
& \stackrel{\eqref{wC14},\eqref{wC15}}{\ge} & \omega_{m+1} \cdot \tau( \hat{\Phi}^{(i)}(\be_{M_{k}}) f  ) - 2 \eta  \\
& \stackrel{\eqref{wC16}}{\ge} & \frac{1}{2} \omega_{m+1} \cdot \tau(\varphi^{(i)} \psi^{(i)}(f)) - 5 \eta,
\end{eqnarray*}
whence
\begin{eqnarray*}
\lefteqn{\tau(\bar{\Phi}(\be_{M_{k}})f)}\\
& \ge & \frac{1}{2} \omega_{m+1} \cdot \tau(\varphi \psi(f)) - (m+1) 5 \eta \\
& \ge & \frac{1}{2} \omega_{m+1} \cdot \tau(f) - (m+1) 5 \eta - \eta \\
& \stackrel{\eqref{wC17}}{\ge} & \frac{1}{4} \omega_{m+1} \cdot \tau(f)
\end{eqnarray*}
for any $\tau \in T(A)$.

Finally, for $x \in M_{k}$ and $i \in \{0, \ldots,m\}$ we have 
\begin{eqnarray*}
\lefteqn{\|[\bar{\Phi}^{(i)}(x),f]\|} \\
& \stackrel{\eqref{wC12}}{\le} & \| [ g_{0,\gamma}(f) (d^{(i)})^{\halb} g_{0,\gamma}(f) \hat{\Phi}^{(i)}(x) g_{0,\gamma}(f)(d^{(i)})^{\halb} g_{0,\gamma}(f), f] \| \\
&& + 2 \eta \|x\| \\
& \le & 2 \|[(d^{(i)})^{\halb},f]\| \cdot \|x\| + \|[\hat{\Phi}^{(i)}(x),f]\| + 2 \eta \|x\| \\
& \stackrel{\eqref{wC9},\eqref{wC10},\eqref{wC15},\eqref{wC18},\eqref{wC19}}{\le} & 6 \eta \|x\| ,
\end{eqnarray*}
whence
\[
\|[\bar{\Phi}(x),f]\| \le 6 (m+1) \eta \|x\| \stackrel{\eqref{wC17}}{\le} \bar{\varepsilon} \|x\|.
\]
\end{nproof}
\en

\bn
\label{wA}
\begin{nprop}
Given $m,k \in \N$, there is $\beta >0$ such that the following holds:

Let $A$ be separable, simple, unital and nonelementary with $\dr A \le m$. Let $F$ be a finite-dimensional $\mathrm{C}^{*}$-algebra, $\varphi:F \to A$ a c.p.c.\ order zero map and $\varepsilon>0$. 

Then, there is a c.p.c.\ order zero map
\[
\Phi:M_{k} \to \her(\varphi(\be_{F})) \subset A
\]
such that
\[
\|[\Phi(x), \varphi(y)]\| \le \varepsilon \|x\| \|y\|
\]
and
\[
\tau(\Phi(\be_{M_{k}})\varphi(a)) \ge \beta  \cdot \tau(\varphi(a))
\]
for all $x \in M_{k}$, $y \in F$, $a \in F_{+}$ and $\tau \in T(A)$.
\end{nprop}

\begin{nproof}
Choose 
\begin{equation}
\label{wA8}
\beta:= \frac{\bar{\beta}}{2},
\end{equation}
where $\bar{\beta}$ comes from Proposition~\ref{wC}. 

Let $A$,  $F$, $\varphi$ and $\varepsilon$ be given. A moment's thought shows that it suffices to verify the assertion of \ref{wA} when $F$ is of the form $M_{r}$ for some $r \in \N$. Let 
\[
\pi: M_{r} \to A''
\]
be the canonical supporting $*$-homomorphism for $\varphi$. Let 
\[
\{f_{ij} \mid i,j=1, \ldots, r\}
\]
denote the canonical matrix units for $M_{r}$. Choose 
\begin{equation}
\label{wA3}
0< \eta < \frac{1}{4} \min \{ \varepsilon, \bar{\beta} \cdot \tau(\varphi(a)) \mid a \in (M_{r})_{+} \mbox{ with } \|a\| = 1, \, \tau \in T(A) \} 
\end{equation}
(the minimum exists and is nonzero by Proposition~\ref{trace-min}). Set
\begin{equation}
\label{wA2}
\bar{\varepsilon}:= \frac{\eta}{r^{2}}.
\end{equation}
By Proposition~\ref{wC}, there is a c.p.c.\ order zero map
\[
\bar{\Phi}: M_{k} \to \her(\varphi(f_{11}))
\]
such that
\begin{equation}
\label{wA1}
\|[\bar{\Phi}(x), \varphi(f_{11})]\| \le \bar{\varepsilon} \cdot \|x\|
\end{equation}
for $x \in M_{k}$ and
\begin{equation}
\label{wA7}
\tau(\bar{\Phi}(\be_{M_{k}})\varphi(f_{11})) \ge \bar{\beta} \cdot \tau(\varphi(f_{11})).
\end{equation}
Define a c.p.c.\ map
\[
\Phi: M_{k} \to A''
\]
by
\begin{equation}
\label{wA6}
\Phi(x):= \sum_{i=1}^{r} \pi(f_{i1}) \bar{\Phi}(x) \pi(f_{1i}).
\end{equation}
It follows from Proposition~\ref{order-zero-facts} that in fact
\[
\Phi(M_{k}) \subset \her(\be_{M_{r}}) \subset A;
\]
it is straightforward to verify (using that $\pi$ is a $*$-homomorphism and $\bar{\Phi}$ has order zero) that $\Phi$ is a c.p.c.\ order zero map. Moreover, for $j,k \in \{1, \ldots, r\}$ and $0 \neq x \in M_{k}$ we have
\begin{eqnarray*}
\lefteqn{\|[\Phi(x),\varphi(f_{jk})]\|}\\
& = & \|\pi(f_{j1}) \bar{\Phi}(x) \varphi(f_{11}) \pi(f_{1k}) - \pi(f_{j1})  \varphi(f_{11}) \bar{\Phi}(x) \pi(f_{1k}) \\
& \stackrel{\eqref{wA1}}{\le} & \bar{\varepsilon} \cdot \|x\|, 
\end{eqnarray*}
whence
\begin{equation}
\label{wA4}
\|[\Phi(x),\varphi(y)]\| \le r^{2} \bar{\varepsilon} \|x\| \|y\| \stackrel{\eqref{wA2}}{\le} \eta \|x\| \|y\| \stackrel{\eqref{wA3}}{\le} \varepsilon \|x\| \|y\|
\end{equation}
for $x \in M_{k}$ and $y \in M_{r}$. 

Finally, for $a \in (M_{r})_{+}$ with norm 1 choose a unitary $u \in M_{r}$ such that 
\begin{equation}
\label{wA5}
u^{*} a u = \sum_{j=1}^{r} \alpha_{j} \cdot f_{jj},
\end{equation}
with $0 \le \alpha_{j}\le 1$. Now, let $ \tau \in T(A)$; denoting the canonical extension of $\tau$ to $A''$ again by $\tau$, we then obtain 
\begin{eqnarray*}
\lefteqn{\tau(\Phi(\be_{M_{k}}) \varphi(a))}\\
& = & \tau(\pi(u^{*}) \Phi(\be_{M_{k}}) \varphi(a)\pi(u)) \\
& = &  \tau(\pi(u^{*}) \Phi(\be_{M_{k}}) \pi(au) \varphi(\be_{M_{r}})) \\
& = &  \tau(\varphi(u^{*}) \Phi(\be_{M_{k}}) \pi(au)) \\
& \stackrel{\eqref{wA4}}{\ge} & \tau(\Phi(\be_{M_{k}}) \varphi(u^{*}au)) - \eta\\
& \stackrel{\eqref{wA5},\eqref{wA6}}{=} & \sum_{j=1}^{r} \alpha_{j} \cdot \tau(\bar{\Phi}(\be_{M_{k}}) \varphi(f_{11})) - \eta \\
& \stackrel{\eqref{wA7}}{\ge} & \sum_{j=1}^{r} \alpha_{j} \cdot \bar{\beta} \cdot \tau(\varphi(f_{11})) - \eta \\
& \stackrel{\eqref{wA5}}{\ge} & \bar{\beta} \cdot \tau(\varphi(a)) - \eta \\
& \stackrel{\eqref{wA3}}{\ge} & \frac{\bar{\beta}}{2} \cdot \tau(\varphi(a)) \\
& \stackrel{\eqref{wA8}}{=} & \beta \cdot \tau(\varphi(a)). 
\end{eqnarray*}
Since $\tau$ and $\varphi$ are both linear, the estimate holds for arbitrary $a \in (M_{r})_{+}$.
\end{nproof}
\en

\bn
\label{3.1aux}
\begin{nprop}
Given $m,k \in \N$, there is $\alpha_{k}>0$ such that the following holds:

Let $A$ be separable, simple, unital and nonelementary with $\dr A \le m$; let $\Eh \subset A$ be a compact subset and let $\gamma >0$ be given.

Then, there is a c.p.c.\ order zero map
\[
\Phi:M_{k} \to A
\]
such that
\[
\|[\Phi(x),a]\| \le \gamma \|x\|
\]
and
\[
\tau(\Phi(\be_{M_{k}})a) \ge \alpha_{k} \cdot \tau(a)
\]
for all $a \in \Eh \cap A_{+}$ and $\tau \in T(A)$.
\end{nprop}

\begin{nproof}
From Propositions~\ref{wA} and \ref{wB}  we obtain $\beta >0$ and $\omega_{k}>0$; set 
\begin{equation}
\label{3.1auxw6}
\alpha_{k}:= \frac{\beta \cdot \omega_{k}}{2}.
\end{equation}
Now, let $A$, $\Eh$ and $\gamma$ as in \ref{3.1aux} be given. We may assume the elements of $\Eh$ to be positive and normalized. Set
\begin{equation}
\label{3.1auxw7}
\bar{\Eh}:= \Eh \cup \{a^{2} \mid a \in \Eh\}
\end{equation}
and
\[
\zeta:= \min \{ \tau(a) \mid a \in \Eh, \, \tau \in T(A) \};
\]
by Proposition~\ref{trace-min}, $\zeta$ exists and is strictly positive. 

Choose $D \in \N$ such that 
\[
D \ge 17 (m+1) 
\]
and
\begin{equation}
\label{3.1auxw5}
\alpha_{k} \cdot \zeta > \frac{((m+1)10 +1) \gamma}{D}.
\end{equation}
Let $(F=F^{(0)}\oplus \ldots \oplus F^{(m)}, \psi, \varphi)$ be an $m$-decomposable c.p.c.\ approximation of $\bar{\Eh}$ to within $(\frac{\gamma}{3D})^2$. Note that, for $a \in \Eh$ and $i \in \{0, \ldots,m\}$, we have 
\begin{equation}
\label{3.1auxw3}
\| \varphi^{(i)}(\be_{F^{(i)}}) \varphi \psi(a) - \varphi^{(i)} \psi^{(i)} (a) \| < \frac{\gamma}{D}.
\end{equation}
Obtain
\[
0< \delta
\]
from Proposition~\ref{almost-order-zero}, applied to $M_{k}$ in place of $F$ and $\frac{\gamma}{D}$ in place of $\gamma$. Choose 
\[
0< \theta <\frac{\gamma}{D}
\]
such that, if $a,b$ are positive elements of norm at most 1 in some $\mathrm{C}^{*}$-algebra which satisfy $\|[a,b]\| \le \varepsilon$, then $\|[a,b^{\halb}]\|, \|[a^{\halb},b^{\halb}]\| < \frac{\delta}{2}$; cf.\ Proposition~\ref{I}.

From Proposition~\ref{wA}, for $i=0, \ldots,m$ we obtain c.p.c.\ order zero maps
\[
\widetilde{\Phi}^{(i)}: M_{k} \to \her(\varphi^{(i)}(\be_{F^{(i)}}))
\]
such that 
\begin{equation}
\label{3.1auxw4}
\|[\widetilde{\Phi}^{(i)}(x),\varphi^{(i)}(y)]\| \le \theta \|x\| \|y\|
\end{equation}
for all $x \in M_{k}$, $y \in F^{(i)}$, and such that
\[
\tau(\widetilde{\Phi}^{i}(\be_{F^{(i)}}) \varphi^{(i)}(b)) \ge \beta \cdot \tau(\varphi^{(i)}(b))
\]
for all $\tau \in T(A)$, $b \in (M_{k})_{+}$. Note that, by our choice of $\theta$, 
\[
\|[\varphi^{(i)}(\be_{F^{(i)}})^{\halb}, \widetilde{\Phi}^{(i)}(x)]\| \le \frac{\delta}{2} \|x\|
\] 
for $x \in M_{k}$.

Next, define a compact subset of $A$ by 
\begin{eqnarray*}
\Fh & := & \{ \varphi \psi(a), \, \widetilde{\Phi}^{(i)}(x), \, \varphi^{(i)}(\be_{F^{(i)}})^{\halb}, \,  \widetilde{\Phi}^{(i)}(\be_{M_{k}})^{\halb} \varphi^{(i)} \psi^{(i)}(a) \widetilde{\Phi}^{(i)}(\be_{M_{k}})^{\halb} \mid \\
&& a \in \bar{\Eh}, \, x \in M_{k} \mbox{ with } \|x \| \le 1, \, i =0, \ldots, m\}.
\end{eqnarray*}

Let $d^{(0)}, \ldots, d^{(m)} \in A_{+}$ be pairwise orthogonal positive elements of norm at most 1 as in Proposition~\ref{wB}. Note that 
\begin{equation}
\label{3.1auxw1}
\|[(d^{(i)})^{\halb},b]\| \le \frac{\delta}{2}
\end{equation}
for all $b \in \Fh$ and $i \in \{0, \ldots, m\}$, so that
\[
\|[(d^{(i)})^{\halb} \varphi^{(i)}(\be_{F^{(i)}})^{\halb}, \widetilde{\Phi}^{(i)}(x)]\| \le \delta \|x\|, \, x \in M_{k}.
\]
By our choice of $\delta$ and Proposition~\ref{almost-order-zero}, for each $i=0, \ldots,m$ there is a c.p.c.\ order zero map 
\[
\Phi^{(i)}:M_{k} \to \her(d^{(i)}) \subset A
\]
such that 
\begin{equation}
\label{3.1auxw2}
\|\Phi^{(i)}(x) - (d^{(i)})^{\halb} \varphi^{(i)}(\be_{F^{(i)}})^{\halb} \widetilde{\Phi}^{(i)}(x) \varphi^{(i)}(\be_{F^{(i)}})^{\halb} (d^{(i)})^{\halb} \| \le \frac{\gamma}{D} \|x\| , \, x \in M_{k}.
\end{equation}
Note that 
\[
\Phi(x):= \sum_{i=0}^{m} \Phi^{(i)}(x), \, x \in M_{k},
\]
defines a  c.p.c.\ order zero map 
\[
\Phi: M_{k} \to A,
\]
since the images of the $\Phi^{(i)}$ are pairwise orthogonal. We proceed to check that $\Phi$ has the right properties. 

For each $i \in \{0, \ldots, m\}$, $a \in \Eh$ and $x \in M_{k}$ with $\|x \| = 1$, we have
\begin{eqnarray*}
\lefteqn{\|[\Phi^{(i)}(x),a]\|}\\
& \le & \| [\Phi^{(i)}(x), \varphi\psi(a)]\| + \frac{2 \gamma}{D} \\
& \stackrel{\eqref{3.1auxw2}}{\le} & \| [(d^{(i)})^{\halb} \varphi^{(i)}(\be_{F^{(i)}})^{\halb} \widetilde{\Phi}^{(i)}(x) \varphi^{(i)}(\be_{F^{(i)}})^{\halb} (d^{(i)})^{\halb}, \varphi \psi(a)] \| +  \frac{4 \gamma}{D} \\
& \stackrel{\eqref{3.1auxw1}}{\le} & \| \widetilde{\Phi}^{(i)}(x) \varphi^{(i)}(\be_{F^{(i)}}) \varphi\psi(a) d^{(i)} - \varphi \psi(a) \varphi^{(i)}(\be_{F^{(i)}}) \widetilde{\Phi}^{(i)}(x) d^{(i)}\| +  \frac{14 \gamma}{D} \\
& \stackrel{\eqref{3.1auxw3}}{\le} &  \| \widetilde{\Phi}^{(i)}(x)  \varphi^{(i)}\psi^{(i)}(a)  - \varphi^{(i)} \psi^{(i)}(a)  \widetilde{\Phi}^{(i)}(x)\| +  \frac{16 \gamma}{D} \\
& \stackrel{\eqref{3.1auxw4}}{\le} & \frac{17 \gamma}{D},
\end{eqnarray*}
whence
\[
\|[\Phi(x),a]\| < \frac{(m+1) 17 \gamma}{D} \le \gamma.
\]
We next observe that, for $a \in \Eh$, $\tau \in T(A)$ and $i \in \{0, \ldots, m\}$, 
\begin{eqnarray*}
\lefteqn{\tau(\Phi^{(i)}(\be_{M_{k}})a)}\\
& \ge & \tau(\Phi^{(i)}(\be_{M_{k}}) \varphi \psi(a)) - \frac{\gamma}{D} \\
& \stackrel{\eqref{3.1auxw2}}{\ge} & \tau((d^{(i)})^{\halb} \varphi^{(i)}(\be_{F^{(i)}})^{\halb} \widetilde{\Phi}^{(i)}(\be_{M_{k}}) \varphi^{(i)}(\be_{F^{(i)}})^{\halb} (d^{(i)})^{\halb} \varphi \psi(a)) - \frac{2\gamma}{D} \\
& \stackrel{\eqref{3.1auxw1}}{\ge} & \tau(d^{(i)} \widetilde{\Phi}^{(i)}(\be_{M_{k}})^{\halb} \varphi^{(i)} \psi^{(i)}(a) \widetilde{\Phi}^{(i)}(\be_{M_{k}})^{\halb}) - \frac{10\gamma}{D} \\
& \stackrel{\ref{wB}}{\ge} &  \omega_{k} \cdot \tau( \widetilde{\Phi}^{(i)}(\be_{M_{k}})^{\halb} \varphi^{(i)} \psi^{(i)}(a) \widetilde{\Phi}^{(i)}(\be_{M_{k}})^{\halb}) - \frac{10\gamma}{D} \\
& \stackrel{\ref{wA}}{\ge} & \beta \omega_{k} \cdot \tau(\varphi^{(i)} \psi^{(i)}(a)) - \frac{10 \gamma}{D},
\end{eqnarray*}
whence
\begin{eqnarray*}
\tau(\Phi(\be_{M_{k}})a) & \ge &  \beta \omega_{k} \cdot \tau(\varphi \psi(a)) - \frac{(m+1)10 \gamma}{D} \\
& \ge & \beta \omega_{k} \cdot \tau(a ) - \frac{((m+1)10 + 1)\gamma}{D} \\
& \stackrel{\eqref{3.1auxw5},\eqref{3.1auxw6},\eqref{3.1auxw7}}{\ge} & \frac{\beta \omega_{k}}{2} \cdot \tau(a) \\
& = & \alpha_{k} \cdot \tau(a).
\end{eqnarray*}
\end{nproof}
\en

\bn
\label{F}
\begin{nlemma}
Let $A$ be a  separable, simple, nonelementary, unital $\mathrm{C}^{*}$-algebra with $\dr A = m < \infty$. Let $\Eh \subset A$ be a compact subset, and let $\gamma >0$ and $k \in \N$ be given. 

Then, there is a c.p.c.\ order zero map 
\[
\Phi : M_{k} \to A
\]
such that
\[
\|[\Phi(x),a]\| \le \gamma \cdot \|x\| \mbox{ for } x \in M_{k}, \, a \in \Eh
\]
and
\[
\tau(\Phi(\be_{M_{k}})) > 1 - \gamma \mbox{ for any } \tau \in T(A).
\]
\end{nlemma}

\begin{nproof}
Set 
\[
\alpha:= \min \{\alpha_{2}, \alpha_{k}\},
\]
where $\alpha_{2}$ and $\alpha_{k}$ come from Proposition~\ref{3.1aux}. We may clearly assume the elements of $\Eh$ to be positive and normalized. Choose $L \in \N$ such that 
\begin{equation}
\label{wF3}
\frac{\alpha^{2}}{2} \sum_{j=0}^{L} \left(1- \frac{\alpha^{2}}{2}\right)^{j} > 1 - \frac{\gamma}{4}.
\end{equation}
(We may clearly assume that $\alpha <1$, so that $\frac{\alpha^{2}}{2} \sum_{j=0}^{\infty} (1- \frac{\alpha^{2}}{2})^{j} =1$, and $L$ exists.)

For $i=L-1, \ldots, 0$ we inductively construct 
\begin{equation}
\label{wF5}
0<\eta_{i}<\eta_{i+1}, \,  0< \zeta_{i}< \eta_{i} \mbox{ and } 0<\gamma_{i}<\gamma_{i+1}
\end{equation} 
as follows.

First, we set 
\begin{equation}
\label{wF4}
\eta_{L}:= \frac{\gamma}{8L} \mbox{ and } \gamma_{L}:= \gamma.
\end{equation}

Now, suppose that for some $i \in \{L-1, \ldots,0\}$, the numbers $\eta_{i+1}$ and $\gamma_{i+1}$ have been constructed. Set
\begin{equation}
\label{wF37}
\eta_{i}:= \min \left\{ \eta_{i+1}, \frac{\gamma_{i+1}}{2(8+3k^{2})} \right\}.
\end{equation}
Choose 
\begin{equation}
\label{wF29}
0< \zeta_{i}< \frac{4\eta_{i}}{k}
\end{equation} 
such that, if 
\[
\Phi^{(1)},\Phi^{(2)}:M_{k}\to A
\]
are c.p.c.\ order zero maps with 
\begin{equation}
\label{wF22}
\|\Phi^{(1)}(\be_{M_{k}}) \Phi^{(2)}(\be_{M_{k}})\| \le \zeta_{i},
\end{equation} 
then there is 
\begin{equation}
\label{wF23}
\ddot{\Phi}:M_{k}\oplus M_{k} \to A
\end{equation} 
c.p.c.\ order zero such that 
\begin{equation}
\label{wF24}
\| \ddot{\Phi}(x,y) - (\Phi^{(1)}(x) + \Phi^{(2)}(y))\| \le \eta_{i} \cdot \max\{\|x\|,\|y\|\}
\end{equation} 
for $x,y \in M_{k}$. This is possible by \cite[Proposition~2.5]{KirWinter:dr}. 

Use Proposition~\ref{almost-order-zero} (in connection with Proposition~\ref{I}) to  choose 
\begin{equation}
\label{wF21}
0< \gamma_{i}< \gamma_{i+1}, \frac{\zeta_{i}}{2(k+1)}
\end{equation} 
such that the following holds: If $$\Phi':M_{k} \to A$$ is a c.p.c.\ order zero map and $h,g \in A$ are positive elements of norm at most 1 satisfying 
\begin{equation}
\label{wF8}
\|[\Phi'(x),h]\|, \|[\Phi'(x),g]\| \le \gamma_{i} \|x\|,
\end{equation} 
then there is a c.p.c.\ order zero map $$\Phi'': M_{k} \to A$$ such that
\begin{equation}
\label{wF6}
\| \Phi'(x) - h^{\halb} (\be-g)^{\halb} \Phi'(x) (\be - g)^{\halb} h^{\halb}\| \le \frac{\zeta_{i}}{4}\cdot \|x\|
\end{equation}
for all $x \in M_{k}$.

By making $\gamma_{i}$ smaller, if necessary, again by Proposition~\ref{I} we may assume that, if $a,b \in A$ are positive with norm  at most 1, and if 
\begin{equation}
\label{wF16}
\|[a,b]\| \le \gamma_{i},
\end{equation}
then 
\begin{equation}
\label{wF17}
\|[a^{\halb},b]\| \le \frac{\zeta_{i}}{2(k+1)} \mbox{ and } \|[g_{\eta_{i},2\eta_{i}}(a),b]\| < \frac{\zeta_{i}}{2(k+1)} \stackrel{\eqref{wF29}}{<} \eta_{i}.
\end{equation}
Moreover, by Proposition~\ref{order-0-almost-commuting} we may assume that, if
\[
\varphi: M_{k} \to A
\]
is a c.p.c.\ order zero map with 
\begin{equation}
\label{wF30}
\|[\varphi(x),b]\| \le \gamma_{i} \cdot \|x\|
\end{equation}
for $x \in M_{k}$, then 
\begin{equation}
\label{wF31}
\| [ ((g_{0,\eta_{i}} - g_{\eta_{i}, 2 \eta_{i}})^{\halb}(\varphi))(x),b]\| \le \eta_{i} \cdot \|x\|.
\end{equation}

Having constructed the numbers $\eta_{i}$, $\zeta_{i}$ and $\gamma_{i}$, we inductively construct c.p.c.\ order zero maps $$\Phi_{0}, \ldots, \Phi_{L}: M_{k} \to A$$ such that
\begin{equation}
\label{wF1}
\tau(\Phi_{i}(\be_{M_{k}})) \ge \frac{\alpha^{2}}{2} \sum_{j=0}^{i} \left(1-\frac{\alpha^{2}}{2}\right)^{j} - \sum_{j=0}^{i-1} 4 \eta_{j}
\end{equation}
and 
\begin{equation}
\label{wF2}
\|[\Phi_{i}(x),a]\| \le \gamma_{i} \cdot \|x\|
\end{equation}
for $i=0, \ldots, L$, $x \in M_{k}$, $a \in \Eh$ and $\tau \in T(A)$. The map $$\Phi:= \Phi_{L}$$ will then have the desired properties, since
\begin{eqnarray*}
\tau(\Phi_{L}(\be_{M_{k}})) &  \stackrel{\eqref{wF1}}{\ge} & \frac{\alpha^{2}}{2} \sum_{j=0}^{L} \left(1-\frac{\alpha^{2}}{2}\right)^{j} - \sum_{j=0}^{L-1} 4 \eta_{j} \\
& \stackrel{\eqref{wF3},\eqref{wF5}}{\ge} & 1 - \frac{\gamma}{4} - 4L \eta_{L} \\
& \stackrel{\eqref{wF4}}{>} & 1 - \gamma
\end{eqnarray*}
and 
\[
\|[\Phi_{L}(x),a]\| \stackrel{\eqref{wF2}}{\le} \gamma_{L} \cdot \|x\| \stackrel{\eqref{wF4}}{=} \gamma \cdot \|x\|
\]
for $x \in M_{k}$, $a \in \Eh$ and $\tau \in T(A)$.

So, let us construct $\Phi_{0}, \ldots, \Phi_{L}$.  We obtain $\Phi_{0}$ from Proposition~\ref{3.1aux} (with $\gamma_{0}$ in place of $\gamma$). Note that $\Phi_{0}$ will satisfy
\[
\tau(\Phi_{0}(\be_{M_{k}})a) \ge \alpha \cdot \tau(a)
\]
for $a \in \Eh$, $\tau \in T(A)$. Suppose $\Phi_{i}$ has been constructed for some $i \in \{0, \ldots, L-1\}$. Let
\[
\pi_{i}: M_{k} \to A''
\]
be the canonical supporting $*$-homomorphism for $\Phi_{i}$. Set
\begin{equation}
\label{wF19}
\bar{\Phi}_{i}:= g_{\eta_{i},2\eta_{i}}(\Phi_{i})
\end{equation}
and (cf.\ \ref{order-zero-notation})
\begin{equation}
\label{wF28}
\hat{\Phi}_{i}:= (g_{0,\eta_{i}}- g_{\eta_{i},2\eta_{i}})(\Phi_{i}).
\end{equation}
By Proposition~\ref{3.1aux}, there is 
\[
\Phi_{i}':M_{k}\to A
\]
c.p.c.\ of order zero such that 
\begin{equation}
\label{wF7}
\|[\Phi_{i}'(x),b]\| \le \gamma_{i} \cdot \|x\|
\end{equation}
and 
\begin{equation}
\label{wF13}
\tau(\Phi_{i}'(\be_{M_{k}})b)\ge \alpha \cdot \tau(b)
\end{equation}
for $x \in M_{k}$, $\tau \in T(A)$ and 
\begin{equation}
\label{wF14}
b \in \Eh \cup \{\be_{A} - g_{0,\eta_{i}}(\Phi_{i}(\be_{M_{k}})), \, g_{0,\eta_{i}}(\Phi_{i}(\be_{M_{k}}))\}.
\end{equation}
Again by Proposition~\ref{3.1aux} (this time in connection with Proposition~\ref{order-0-almost-commuting}), there is
\[
\Psi_{i}: M_{2} \to A
\]
c.p.c.\ order zero such that 
\begin{equation}
\label{wF8a}
\|[\Psi_{i}(x),b]\| \le \gamma_{i} \cdot \|x\|
\end{equation}
and
\begin{equation}
\label{wF12}
\tau(\Psi_{i}(e_{jj})c) \ge \frac{\alpha}{2} \cdot \tau(c) - \frac{\zeta_{i}}{4}
\end{equation}
for $x \in M_{2}$, $j=1,2$, $\tau \in T(A)$,  
\begin{eqnarray}
b & \in & \Eh \cup \{ \Phi_{i}'(x),\, \hat{\Phi}_{i}(e_{11}), \, \bar{\Phi}_{i}(\be_{M_{k}}), \nonumber \\
&& (\be_{A} - g_{0,\eta_{i}}(\Phi_{i}(\be_{M_{k}})))^{\halb}   \Phi'_{i}(y)  (\be_{A} - g_{0,\eta_{i}}(\Phi_{i}(\be_{M_{k}})))^{\halb}, \nonumber \\
&& \pi_{i}(e_{j1}) \hat{\Phi}_{i}(e_{11})^{\halb} \mid \nonumber \\
&& \mbox{for } j = 1, \ldots, k, \mbox{ and } y \in M_{k} \mbox{ with } \|y\| \le 1 \} \label{wF9}
\end{eqnarray}
and
\begin{eqnarray}
c & \in & \{ (\be_{A} - g_{0,\eta_{i}}(\Phi_{i}(\be_{M_{k}})))^{\halb}   \Phi'_{i}(y)  (\be_{A} - g_{0,\eta_{i}}(\Phi_{i}(\be_{M_{k}})))^{\halb} \mid  \nonumber \\
&& y \in (M_{k})_{+} \mbox{ with } \|y\| \le 1\}. \label{wF11}
\end{eqnarray}
Define a c.p.c.\ order zero  map 
\[
\check{\Phi}_{i}: M_{k} \to A
\]
by
\begin{equation}
\label{wF20}
\check{\Phi}_{i}(e_{jl}):= \pi_{i}(e_{j1}) \hat{\Phi}_{i}(e_{11})^{\halb}\Psi_{i}(e_{22}) \hat{\Phi}_{i}(e_{11})^{\halb} \pi_{i}(e_{1l})
\end{equation}
for $j,l \in \{1, \ldots,k\}$. That $\check{\Phi}_{i}$ indeed maps $M_{k}$ to $A$ (as opposed to $A''$), follows from the fact that $\pi_{i}(x) \hat{\Phi}_{i}^{\halb}(y) \in A $ for $x,y \in M_{k}$. Define 
\[
\tilde{\Phi}_{i}:M_{k} \to A
\]
by
\begin{equation}
\label{wF15}
\tilde{\Phi}_{i}(x) := \bar{\Phi}_{i}(x) + \check{\Phi}_{i}(x) \mbox{ for } x \in M_{k};
\end{equation}
one immediately checks that $\tilde{\Phi}_{i}$ is a c.p.c.\ order zero map. Note also that
\begin{eqnarray}
\label{wF27}
\tau(\check{\Phi}_{i}(\be_{M_{k}})) & = &  k \cdot \tau(\Psi_{i}(e_{22}) \hat{\Phi}_{i}(e_{11})) \nonumber \\
& \stackrel{\eqref{wF12}}{\ge} & k  \frac{\alpha}{2} \cdot \tau(\hat{\Phi}_{i}(e_{11})) - k \cdot \frac{\zeta_{i}}{4} \nonumber \\
& = & \frac{\alpha}{2} \cdot \tau(\hat{\Phi}_{i}(\be_{M_{k}})) - k \cdot \frac{\zeta_{i}}{4} .
\end{eqnarray}
By  our choice of $\gamma_{i}$ (cf.\ \eqref{wF8}, \eqref{wF6}), there is a c.p.c.\ order zero map
\[
\Phi_{i}'': M_{k} \to A
\]
such that
\begin{eqnarray}
 \|\Psi_{i}(e_{11})^{\halb} (\be_{A}- g_{0,\eta_{i}}(\Phi_{i}(\be_{M_{k}})))^\halb \Phi_{i}'(x)  && \nonumber \\ 
(\be_{A}- g_{0,\eta_{i}}(\Phi_{i}(\be_{M_{k}})))^\halb \Psi_{i}(e_{11})^{\halb}  -   \Phi''_{i}(x)  \| & \le &  \frac{\zeta_{i}}{4} \cdot \|x\|   \label{wF10}
\end{eqnarray}
for $x \in M_{k}$; cf.\ \eqref{wF6}, \eqref{wF8a}, \eqref{wF7} and  \eqref{wF9}. Note that 
\begin{eqnarray}
\label{wF26}
\lefteqn{\tau(\Phi''_{i}(\be_{M_{k}}))} \nonumber \\
 & \stackrel{\eqref{wF10}}{\ge} & \tau( \Psi_{i}(e_{11}) (\be_{A}- g_{0,\eta_{i}}(\Phi_{i}(\be_{M_{k}})))^\halb \Phi_{i}'(\be_{M_{k}})  (\be_{A}- g_{0,\eta_{i}}(\Phi_{i}(\be_{M_{k}})))^\halb) \nonumber \\
&&  - \frac{\zeta_{i}}{4} \nonumber \\
& \stackrel{\eqref{wF12},\eqref{wF11}}{>} & \frac{\alpha}{2} \cdot \tau(  (\be_{A}- g_{0,\eta_{i}}(\Phi_{i}(\be_{M_{k}})))^\halb \Phi_{i}'(\be_{M_{k}})  (\be_{A}- g_{0,\eta_{i}}(\Phi_{i}(\be_{M_{k}})))^\halb)   - \frac{\zeta_{i}}{2}\nonumber  \\
& = & \frac{\alpha}{2} \cdot \tau( \Phi_{i}'(\be_{M_{k}})  (\be_{A}- g_{0,\eta_{i}}(\Phi_{i}(\be_{M_{k}}))))  - \frac{\zeta_{i}}{2} \nonumber \\
& \stackrel{\eqref{wF13},\eqref{wF14}}{>} & \frac{\alpha^{2}}{2} \cdot \tau(\be_{A}- g_{0,\eta_{i}}(\Phi_{i}(\be_{M_{k}})))  - \frac{\zeta_{i}}{2}.
\end{eqnarray}
We have
\begin{eqnarray*}
\lefteqn{\| \Phi''_{i}(\be_{M_{k}}) \tilde{\Phi}_{i}(\be_{M_{k}})\|}\\
 & \stackrel{\eqref{wF10},\eqref{wF15}}{\le} & \| \Psi_{i}(e_{11})^{\halb} (\be_{A}- g_{0,\eta_{i}}(\Phi_{i}(\be_{M_{k}})))^\halb \Phi_{i}'(\be_{M_{k}})  (\be_{A}- g_{0,\eta_{i}}(\Phi_{i}(\be_{M_{k}})))^\halb \\
 &&  \Psi_{i}(e_{11})^{\halb} (\bar{\Phi}_{i} (\be_{M_{k}}) + \check{\Phi}_{i}(\be_{M_{k}})) \| + \frac{\zeta_{i}}{4} \\
& \le & \| (\be_{A}- g_{0,\eta_{i}}(\Phi_{i}(\be_{M_{k}})))^\halb \Psi_{i}(e_{11})^{\halb} \bar{\Phi}_{i} (\be_{M_{k}}) \| \\
& & +  \|  \Psi_{i}(e_{11})^{\halb} \check{\Phi}_{i}(\be_{M_{k}}))\| + \frac{\zeta_{i}}{4} \\
& \le & (k+1) \gamma_{i} + \frac{\zeta_{i}}{4} \\
& \stackrel{\eqref{wF21}}{<} & \zeta_{i}
\end{eqnarray*}
(where for the second last inequality we have used \eqref{wF20}, \eqref{wF8a}, \eqref{wF9}, \eqref{wF19}, \eqref{wF16} and \eqref{wF17}), whence, by our choice of $\zeta_{i}$, there is a c.p.c.\ order zero map 
\[
\ddot{\Phi}_{i}: M_{k} \oplus M_{k} \to A
\]
such that 
\begin{equation}
\label{wF25}
\| \ddot{\Phi}_{i}(x,y) - (\Phi_{i}''(x) + \tilde{\Phi}_{i}(y))\| \le \eta_{i} \cdot \max\{\|x\|, \|y\|\}
\end{equation}
for $x,y \in M_{k}$; cf.\ \eqref{wF22}, \eqref{wF23}, \eqref{wF24}. Now define
\[
\Phi_{i+1}: M_{k} \to A
\]
by 
\begin{equation}
\label{wF35}
\Phi_{i+1}(x) := \ddot{\Phi}_{i}(x,x). 
\end{equation}
We check that $\Phi_{i+1}$ satisfies \eqref{wF1} and \eqref{wF2}. First, we estimate
\begin{eqnarray*}
\lefteqn{\tau(\Phi_{i+1}(\be_{M_{k}}))} \\
& \stackrel{\eqref{wF25}}{\ge} & \tau(\Phi_{i}''(\be_{M_{k}})) + \tau(\tilde{\Phi}_{i}(\be_{M_{k}})) - \eta_{i} \\
& \stackrel{\eqref{wF26}}{\ge} & \frac{\alpha^{2}}{2} \cdot \tau(\be_{A} - g_{0, \eta_{i}}(\Phi_{i}(\be_{M_{k}}))) + \tau(\tilde{\Phi}_{i}(\be_{M_{k}})) - 2\eta_{i} \\
& \stackrel{\eqref{wF27},\eqref{wF15}}{\ge} &  \frac{\alpha^{2}}{2} \cdot \tau(\be_{A} - g_{0, \eta_{i}}(\Phi_{i}(\be_{M_{k}}))) + \tau(\bar{\Phi}_{i}(\be_{M_{k}})) + \frac{\alpha}{2} \cdot \tau(\hat{\Phi}_{i}(\be_{M_{k}})) \\
&& - 2\eta_{i} - k \cdot \frac{\zeta_{i}}{4}\\
& \stackrel{\eqref{wF19},\eqref{wF28},\eqref{wF29},\eqref{wF21}}{\ge} & \frac{\alpha^{2}}{2} \cdot \tau(\be_{A} - \bar{\Phi}_{i}(\be_{M_{k}})) + \tau(\bar{\Phi}_{i}(\be_{M_{k}})) - 3 \eta_{i} \\
& = & \frac{\alpha^{2}}{2} + \left(1 - \frac{\alpha^{2}}{2}\right) \cdot \tau(\bar{\Phi}_{i}(\be_{M_{k}})) - 3 \eta_{i} \\
& \stackrel{\eqref{wF19}}{\ge} &  \frac{\alpha^{2}}{2} + \left(1 - \frac{\alpha^{2}}{2}\right) \cdot \tau(\Phi_{i}(\be_{M_{k}})) - 4 \eta_{i} \\
& \stackrel{\eqref{wF1}}{\ge} &  \frac{\alpha^{2}}{2} + \left( 1 -  \frac{\alpha^{2}}{2}\right)  \frac{\alpha^{2}}{2}  \sum_{j=0}^{i} \left(1-  \frac{\alpha^{2}}{2}\right)^{j} - \sum_{j=0}^{i-1} 4 \eta_{j} - 4 \eta_{i} \\
& = &  \frac{\alpha^{2}}{2} \sum_{j=0}^{i+1} \left(1-  \frac{\alpha^{2}}{2}\right)^{j} - \sum_{j=0}^{i} 4 \eta_{j}.
\end{eqnarray*}
Next, we check that for $x \in M_{k}$ with $\|x\|=1$ and $a \in \Eh$, 
\begin{eqnarray}
\label{wF36}
\lefteqn{\| [\Phi''_{i}(x),a]\|} \nonumber \\
& \stackrel{\eqref{wF10},\eqref{wF29}}{\le} & \| [\Psi_{i}(e_{11})^{\halb} ( \be_{A} - g_{0,\eta_{i}}(\Phi_{i}(\be_{M_{k}})))^{\halb} \Phi'_{i}(x) \nonumber \\
&& ( \be_{A} - g_{0,\eta_{i}}(\Phi_{i}(\be_{M_{k}})))^{\halb} \Psi_{i}(e_{11})^{\halb},a]\|  + \eta_{i}  \nonumber \\
& < & 6 \eta_{i},
\end{eqnarray}
where for the last inequality we have used \eqref{wF8a}, \eqref{wF9}, \eqref{wF2}, \eqref{wF7}, \eqref{wF16} and \eqref{wF17}. Similarly, we obtain
\begin{equation}
\label{wF32}
\|[\bar{\Phi}_{i}(x),a]\| < \eta_{i}.
\end{equation}
Moreover,
\begin{eqnarray}
\label{wF33}
\lefteqn{ \|[\pi_{i}(e_{j,1}) \hat{\Phi}_{i}(e_{11})^{\halb} \Psi_{i}(e_{22}) \hat{\Phi}_{i}(e_{11})^{\halb} \pi_{i}(e_{1l}),a]\| } \nonumber \\
& \stackrel{\eqref{wF8},\eqref{wF9}}{\le} &   \| [\pi_{i}(e_{j1}) \hat{\Phi}_{i}(e_{11})\pi_{i}(e_{1l}) \Psi_{i}(e_{22}),a] \| + 2 \gamma_{i} \nonumber \\
& = &  \| [ \hat{\Phi}_{i}(e_{jl}) \Psi_{i}(e_{22}),a] \| + 2 \gamma_{i} \nonumber \\
& \le & 2 \gamma_{i} + \eta_{i} \nonumber \\
& \stackrel{\eqref{wF21},\eqref{wF5}}{<} & 3 \eta_{i}
\end{eqnarray}
for $j,l=1, \ldots,k$ (where for the second last inequality we have used \eqref{wF8}, \eqref{wF9}, \eqref{wF28}, \eqref{wF2}, \eqref{wF30} and \eqref{wF31}), whence
\begin{equation}
\label{wF34}
\|[\tilde{\Phi}_{i}(x),a]\| < \eta_{i} + 3 k^{2} \eta_{i}.
\end{equation}
We obtain
\begin{eqnarray*}
\|[\Phi_{i+1}(x),a]\| & \stackrel{\eqref{wF35},\eqref{wF25}}{\le} & \|[\Phi_{i}''(x),a] \| + \|[\tilde{\Phi}_{i}(x),a]\| + \eta_{i} \\
& \stackrel{\eqref{wF36},\eqref{wF34}}{<} & (8+3 k^{2}) \eta_{i} \\
& \stackrel{\eqref{wF37}}{<} & \gamma_{i+1}
\end{eqnarray*}
for $x \in M_{k}$ with $\|x\|=1$ and $a \in \Eh$. 

Induction yields maps $\Phi_{0}, \ldots, \Phi_{L}$, as desired. We are done.
\end{nproof}
\en

\bn
\label{n-comparison}
Let us recall Lemma~6.1 from \cite{TomsWinter:VI}. Here, $d_{\tau}$ denotes the dimension function associated to the trace $\tau$ and $\langle a \rangle$ is the class of $a \in A_{+}$ in the Cuntz semigroup $W(A)$.

\begin{nlemma}
Let $A$ be a  separable, simple, unital and nonelementary $\mathrm{C}^{*}$-algebra with $\dr A = m < \infty$. If $a,d^{(0)},\ldots,d^{(m)} \in A_{+}$ satisfy 
\[
d_{\tau}(a) < d_{\tau}(d^{(i)})
\]
for $i=0,\ldots,m$ and every $\tau \in T(A)$, then 
\[
\langle a \rangle \le \langle d^{(0)} \rangle + \ldots + \langle d^{(m)} \rangle.
\]
\end{nlemma}
\en

\bn
\label{cor-dr-comparison}
\begin{ncor}
Let $A$ be a simple, separable, unital $\mathrm{C}^{*}$-algebra with decomposition rank $m < \infty$. 

Let $a,b \in A$ be positive elements satisfying
\[
d_{\tau}(a) < \frac{1}{m+1} \cdot \tau(b) \mbox{ for any } \tau \in T(A).
\]
Then, $\langle a \rangle \le \langle b \rangle$ in $W(A)$.
\end{ncor}

\begin{nproof}
If $\langle a \rangle = \langle p \rangle$ for some projection $p \in A$, then $d_{\tau}(a) = \tau(p)$ for each $\tau \in T(A)$ and 
\[
\min \left\{ \frac{1}{m+1} \cdot  \tau(b) - \tau(p) \mid \tau \in T(A) \right\}
\]
exists and is strictly positive; cf.\  Proposition~\ref{trace-min}.

Furthermore, for any $\eta>0$ we have
\[
d_{\tau}((a-\eta)_{+}) + \tau((g_{0,\eta/2} - g_{\eta/2,\eta})(a)) \le d_{\tau}(a)
\]
for each $\tau \in T(A)$. Now if $a$ is not Cuntz equivalent to  a projection, then, as above, 
\[
\min\{ \tau((g_{0,\eta/2}- g_{\eta/2,\eta})(a)) \mid \tau \in T(A) \}
\]
exists and is strictly positive.

Thus, in both cases for any $\eta>0$ we may assume that there is $\zeta_{\eta}>0$ such that 
\[
d_{\tau}((a-\eta)_{+}) + \zeta_{\eta} < \frac{1}{m+1} \tau(b) 
\]
for each $\tau \in T(A)$. 

Since $\langle a \rangle \le \langle b \rangle$ iff $\langle (a-\eta)_{+} \rangle \le \langle b \rangle$ for any $\eta >0$, it will be enough to prove the corollary under the hypothesis that there is $\zeta>0$ such that 
\begin{equation}
\label{cdcw3}
d_{\tau}(a) + \zeta < \frac{1}{m+1} \cdot \tau(b)
\end{equation}
for any $\tau \in T(A)$. 

Now by Proposition~\ref{almost-order-zero} (in connection with Proposition~\ref{order-0-almost-commuting}), there is $\gamma>0$ such that, if $\Phi:M_{m+1} \to A$ is c.p.c.\ of order zero with $\|[\Phi(x),b]\|\le \gamma \|x\|$ for $x \in M_{k}$, then there is 
\[
\widetilde{\Phi}: M_{m+1} \to \her (b) \subset A
\] 
with 
\begin{equation}
\label{cdcw1}
\| \widetilde{\Phi}(x) - b^{\halb} \Phi(x) b^{\halb}\| \le \frac{\zeta}{2} \|x\|
\end{equation}
for $x \in M_{k}$. We may assume that $\gamma < \frac{\zeta}{2}$. 

By Lemma~\ref{F}, there is 
\[
\Phi:M_{m+1} \to A 
\]
c.p.c.\ of order zero such that
\[
\|[\Phi(x),b]\| < \gamma \|x\|
\]
for $x \in M_{k}$ and such that 
\begin{equation}
\label{cdcw2}
\tau(\Phi(\be_{m+1}))> 1-\gamma
\end{equation}
for any $\tau \in T(A)$. Thus, we may choose $\widetilde{\Phi}$ as above. For $i \in \{1, \ldots, m+1\}$ and $\tau \in T(A)$ we then have
\begin{eqnarray*}
\tau(\widetilde{\Phi}(e_{ii})) & = & \frac{1}{m+1}\cdot  \tau(\widetilde{\Phi}(\be_{m+1})) \\
& \stackrel{\eqref{cdcw1}}{\ge} & \frac{1}{m+1} \cdot  \left(\tau(b^{\halb}\Phi(\be_{m+1})b^{\halb}) - \frac{\zeta}{2}\right) \\
& = & \frac{1}{m+1}\cdot  \left(\tau(b) - \tau(b^{\halb}(\be_{A} - \Phi(\be_{m+1}))b^{\halb}) - \frac{\zeta}{2}\right) \\
& \stackrel{\eqref{cdcw2}}{>} & \frac{1}{m+1}\cdot  \tau(b) - \gamma - \frac{\zeta}{2}\\
& \ge & \frac{1}{m+1} \cdot \tau(b) - \zeta \\
& \stackrel{\eqref{cdcw3}}{>} & d_{\tau}(a).
\end{eqnarray*}
Now by Lemma~\ref{n-comparison}, we have $$a \precsim  \widetilde{\Phi}(e_{11}) + \ldots + \widetilde{\Phi}(e_{m+1,m+1}) = \widetilde{\Phi}(\be_{m+1}) \in \her (b).$$ Therefore, $a \precsim b$ and $\langle a \rangle \le \langle b \rangle$.
\end{nproof}
\en

\bn
\begin{nremark}
As mentioned above, with Lemma~\ref{F} and Corollary~\ref{cor-dr-comparison} at hand, at this point one could directly prove that,  for simple, unital $\mathrm{C}^{*}$-algebras, finite decomposition rank implies strict comparison of positive elements. However, this would require some small additional technicalities (as well as the reproduction of some of the arguments of \cite{Ror:Z-absorbing}), so we rather proceed in our proof of \ref{finite-con} (i) $\Longrightarrow$ (ii) and recall that \ref{finite-con} (ii) $\Longrightarrow$ (iii) was shown in \cite{Ror:Z-absorbing}.  
\end{nremark}
\en

\section{From tracial matrix cone absorption to almost central dimension drop embeddings}
\label{almost-central-dimension-drop-embeddings}

\noindent
The purpose of this section is to prove Proposition~\ref{D} (which in turn refines Proposition~\ref{E}); this will allow us to construct $\varphi$ and $v$ as in relations $\mathcal{R}(n,\mathcal{F},\eta)$ of \ref{R-relations} in a simple unital $\mathrm{C}^{*}$-algebra with finite decompostion rank. While existence of the map $\varphi$ from $\mathcal{R}(n,\mathcal{F},\eta)$ essentially follows from Lemma~\ref{F}, the element $v$ will be provided by Proposition~\ref{D}.

\bn
For the sequel, see \ref{d-tau-notation} for the definition of dimension functions associated to positive tracial functionals.

\label{G}
\begin{nprop}
Let $A$ be a unital $\mathrm{C}^{*}$-algebra, and let $a \in A_{+}$ be a positive element of norm at most 1 satisfying $\tau(a) > 0$ for any $\tau \in T(A)$. 

Then, for any $k \in \N$ and $0 < \beta < 1$ there is $\gamma >0$ such that the following holds: If  $h \in A_{+}$ is another positive element of norm at most 1 satisfying 
\[
\tau(h)> 1-\gamma \mbox{ for any } \tau \in T(A),
\]
then
\[
d_{\tau}((a^{\frac{1}{2}}(\be_{A}-h)a^{\frac{1}{2}})- \beta)_{+} < \frac{1}{k} \cdot \tau(a) \mbox{ for any } \tau \in T(A).
\]
\end{nprop}

\begin{nproof}
Set 
\begin{equation}
\label{G2}
\zeta:= \min \{\tau(a) \mid \tau \in T(A)\}.
\end{equation}
Since $A$ is unital, $T(A)$ is compact whence $\zeta$ exists and is strictly bigger than zero; cf.\ the proof of Proposition~\ref{trace-min}. Choose some 
\begin{equation}
\label{G1}
\gamma < \frac{\zeta  \beta}{2 k}.
\end{equation}
We then compute
\begin{eqnarray*}
d_{\tau}((a^{\frac{1}{2}}(\be_{A}-h)a^{\frac{1}{2}})- \beta)_{+} & \stackrel{\eqref{dtnw1}}{\le} & \tau(g_{\beta/2,\beta}((a^{\frac{1}{2}}(\be_{A}-h)a^{\frac{1}{2}}))\\
& \le & \frac{1}{\frac{\beta}{2}} \cdot  \tau(a^{\frac{1}{2}}(\be_{A}-h)a^{\frac{1}{2}})\\
& \le & \frac{1}{\frac{\beta}{2}} \cdot  \tau(\be_{A}-h)\\
& < & \frac{1}{\frac{\beta}{2}} \gamma\\
& \stackrel{\eqref{G1}}{<} & \frac{1}{k} \zeta \\
& \stackrel{\eqref{G2}}{\le} &  \frac{1}{k} \tau(a)
\end{eqnarray*}
for any $\tau \in T(A)$.
\end{nproof}
\en

\bn
\label{E-aux}
\begin{nprop}
Let $A$ be a simple, separable, unital $\mathrm{C}^{*}$-algebra with $\dr A = m < \infty$. Let 
\[
F  = M_{r_{1}} \oplus \ldots \oplus M_{r_{s}}
\] 
be a finite dimensional $\mathrm{C}^{*}$-algebra, $n \in \N$ and let 
\[
\psi: F \to A
\]
be a c.p.c.\ order zero map. For any $0<\bar{\gamma}_{0}$ and $0< \zeta < 1$ there are $0<\bar{\gamma}_{1}$ and $0<\bar{\beta}$ such that the following holds: 

If $ i_{0} \in \{0, \ldots,m\}$ and 
\[
\Phi: M_{m+1} \otimes M_{n} \otimes M_{2} \to A
\]
is a c.p.c.\ order zero map satisfying
\[
\| [\psi(x), \Phi(y)]\| \le \bar{\beta} \|x\| \|y\| \mbox{ for any }  x \in F \mbox{ and }  y \in M_{m+1} \otimes M_{n} \otimes M_{2}
\]
and
\begin{equation}
\label{E-auxw1}
\tau(\Phi(\be_{M_{m+1} \otimes M_{n} \otimes M_{2}})) > 1 - \bar{\beta},
\end{equation}
then there are $\bar{v}^{(j)} \in A$, $j=1, \ldots,s$, of norm at most one with the following properties:

\begin{enumerate}
\item $\|(\bar{v}^{(j)})^{*} \bar{v}^{(j)} - \psi^{(j)}(f_{11}^{(j)})^{\halb}(\be_{A}- \Phi(\be_{M_{m+1} \otimes M_{n} \otimes M_{2}}))  \psi^{(j)}(f_{11}^{(j)})^{\halb}\| < \bar{\gamma}_{0}$
\item $\|[\bar{v}^{(j)}, \psi^{(j)}(f_{11}^{(j)})]\| < \bar{\gamma}_{0}$
\item $\|\bar{v}^{(j)} - g_{0,\bar{\gamma}_{1}}(\psi^{(j)}(f_{11}^{(j)})) \bar{v}^{(j)}\| < \bar{\gamma}_{0}$
\item $\bar{v}^{(j)} = g_{\zeta,\zeta+\bar{\gamma}_{1}}(\Phi(d_{i_{0}i_{0}}\otimes e_{11} \otimes \be_{M_{2}})) \bar{v}^{(j)}$.
\end{enumerate}
Here, $\{f_{ii'}^{(j)} \mid j=1, \ldots,s; \, i,i' =1, \ldots,r^{j}\}$ and $\{d_{ii'} \mid i,i' = 1, \ldots,m+1\}$ for $F$ and $M_{m+1}$, respectively; the canonical matrix units for $M_{n}$ and $M_{2}$ are both denoted by $\{e_{ii'}\}$.  
\end{nprop}

\begin{nproof}
Choose $\bar{t} \in \N$ such that 
\begin{equation}
\label{ww22}
\bar{t} > \frac{8}{\bar{\gamma}_{0}^{2}}
\end{equation}
and define $h_{t,2\bar{t}} \in \Ch([0,1])$ for $t = 1, \ldots, 2 \bar{t}$ by
\[
h_{t,2\bar{t}}(x) := \left\{
\begin{array}{ll}
0, & x \le t - \frac{1}{2\bar{t}}\\
1, & x=t\\
0, & x \ge t + \frac{1}{2\bar{t}}\\
\mbox{linear } & \mbox{else.} 
\end{array}
\right.
\]
For $j=1, \ldots, s$ and $t=1, \ldots, 2\bar{t}$ set
\begin{equation}
\label{ww2}
a^{(j)}_{t} := h_{t,2\bar{t}}(\psi^{(j)}(f_{11}^{(j)}))
\end{equation}
and 
\[
\alpha_{t} := \frac{t}{2 \bar{t}};
\]
note that
\begin{equation}
\label{ww13}
\psi^{(j)}(f_{11}^{(j)}) = \sum_{t=1}^{2 \bar{t}} \alpha_{t} \cdot a^{(j)}_{t}.
\end{equation}
Set
\[
J:= \{(j,t) \in \{1, \ldots, s\} \times \{1, \ldots, 2\bar{t}\} \mid a^{(j)}_{t} \neq 0 \},
\]
\[
J_{1} := \{(j,t) \in J \mid t \mbox{ odd}\} \mbox{ and } J_{2} := \{(j,t) \in J \mid t \mbox{ even}\}.
\]
Choose  
\begin{equation}
\label{ww16}
0 < \delta_{1} < \frac{\bar{\gamma_{0}}^{4}}{(2 \bar{t})^{8}}.
\end{equation}
such that $f_{\delta_{1},2 \delta_{1}}(a^{(j)}_{t}) \neq 0$ if $(j,t) \in J$.  Set 
\begin{equation}
\label{ww31}
\bar{a}^{(j)}_{t}:= f_{\delta_{1}, 2 \delta_{1}}(a^{(j)}_{t});
\end{equation}
note that 
\begin{equation}
\label{www1}
\bar{a}^{(j)}_{t} \neq 0 \mbox{ if } (j,t) \in J.
\end{equation} 
By functional calculus, \eqref{ww2} and \eqref{ww31} we also have 
\begin{equation}
\label{ww6}
a_{t}^{(j)} \perp a_{t'}^{(j)} \mbox{ and } \bar{a}_{t}^{(j)} \perp \bar{a}_{t'}^{(j)}
\end{equation}
if $t \neq t'$ are both odd or both even; one checks that 
\begin{equation}
\label{ww14}
\|  (a_{t}^{(j)})^{\halb} - (\bar{a}_{t}^{(j)})^{\halb}\| < (2 \delta_{1})^{\halb}.
\end{equation}
Set
\begin{equation}
\label{E-auxw3}
k:= 24(m+1)^{2}n.
\end{equation}
Choose 
\begin{equation}
\label{ww15b}
0< \delta_{2} < \delta_{1}
\end{equation} 
such that 
\begin{equation}
\label{ww3a}
\tau(\bar{a}^{(j)}_{t})> 30 \delta_{2} \; \forall \, \tau \in T(A), \, (j,t) \in J;
\end{equation}
this is possible by \eqref{www1} and Proposition~\ref{trace-min}.

Choose   
\begin{equation}
\label{ww15a}
0 < \bar{\gamma}_{1} < \frac{1 - \zeta}{2}, \, \delta_{2}.
\end{equation}
Choose  
\begin{equation}
\label{ww10a}
0 < \bar{\beta}< \bar{\gamma}_{1}, \, \frac{\delta_{2}(1-\zeta)}{2}, \, \frac{\delta_{1}}{k} \cdot \min\{\tau(\bar{a}^{(j)}_{t}) \mid (j,t) \in J, \, \tau \in T(A)\},
\end{equation}
\[
\bar{\beta} < \bar{\gamma}_{1} \cdot (1-(\zeta +\bar{\gamma}_{1})) (\stackrel{\eqref{ww15a}}{>} 0)
\]
and such that Proposition~\ref{G} holds with $\delta_{1}$ in place of $\delta$ and for every $\bar{a}^{(j)}_{t}$ in place of $a$, for $(j,t) \in J$. By making $\bar{\beta}$ smaller, if necessary, we may also assume that the assertion of Proposition~\ref{order-0-almost-commuting} holds with $\bar{\beta}$ in place of $\beta$, $\delta_{1}$ in place of $\gamma$ and 
\begin{eqnarray*}
\lefteqn{\{f_{\zeta+\bar{\gamma}_{1},\zeta+2 \bar{\gamma}_{1}}((\, . \,)^{\halb}), \, (\, . \, -\delta_{2})_{+}, \, f_{\delta_{1},2\delta_{1}}(h_{t,2\bar{t}}(\, . \,)^{\halb}),  }\\
&& \, (\, . \, )^{\halb}, \, h_{t,2\bar{t}}(\, .\,)^{\halb}, \, g_{0,\delta_{1}}(h_{t,2\bar{t}}(\, .  \,)^{\halb}),  \, g_{0,\delta_{2}}(\, .\,)\}
\end{eqnarray*}
in place of $\Gh$.  

Now, let $\Phi$ be as in the Proposition; we need to construct $\bar{v}^{(j)}$, $j=1, \ldots, s$, satisfying properties (i) through (iv). That the $\bar{v}^{(j)}$ can be chosen to be normalized then follows from (i) (which entails that they will have norm not much greater than one) and a standard rescaling argument.

For $(j,t) \in J$, $i=1,2$ and $\tau \in T(A)$ we obtain
\begin{eqnarray}
\lefteqn{2 \tau(f_{\zeta+\bar{\gamma}_{1},\zeta+2\bar{\gamma}_{1}}(\Phi(\be_{M_{m+1}\otimes M_{n}} \otimes e_{ii}))^{\halb} \bar{a}_{t}^{(j)}   f_{\zeta+\bar{\gamma}_{1},\zeta+2\bar{\gamma}_{1}}(\Phi(\be_{M_{m+1}\otimes M_{n}} \otimes e_{ii}))^{\halb})} \nonumber \\
& \stackrel{\eqref{order-0-almost-commuting}\mathrm{(iii)}}{\ge} & \tau(f_{\zeta+\bar{\gamma}_{1},\zeta+2\bar{\gamma}_{1}}(\Phi(\be_{M_{m+1} \otimes M_{n} \otimes M_{2}}))^{\halb} \bar{a}_{t}^{(j)} \nonumber \\
&&  f_{\zeta+\bar{\gamma}_{1},\zeta+2\bar{\gamma}_{1}}(\Phi(\be_{M_{m+1} \otimes M_{n} \otimes M_{2}}))^{\halb})  - \bar{\gamma}_{1} \nonumber \\
& \ge & \tau(\bar{a}_{t}^{(j)}) - \tau(\be_{A} - f_{\zeta+\bar{\gamma}_{1},\zeta+2\bar{\gamma}_{1}}(\Phi(\be_{M_{m+1}\otimes M_{n} \otimes M_{2}}))) - \bar{\gamma}_{1} \nonumber \\
& \ge & \tau(\bar{a}_{t}^{(j)}) - \frac{1}{1- (\zeta + \bar{\gamma}_{1})} \tau(\be_{A} -\Phi(\be_{M_{m+1}\otimes M_{n} \otimes M_{2}})) - \bar{\gamma}_{1} \nonumber \\
& \stackrel{\eqref{ww10a}}{>} & \tau(\bar{a}_{t}^{(j)}) - \frac{\bar{\beta}}{1- (\zeta + \bar{\gamma}_{1})} - \bar{\gamma}_{1} \nonumber \\
& \stackrel{\eqref{ww15a},\eqref{ww15b}}{\ge} &  \tau(\bar{a}_{t}^{(j)}) - 2 \delta_{2}. \label{ww3}
\end{eqnarray}
By \eqref{ww3a}, this in particular implies 
\begin{eqnarray}
& &\tau(f_{\zeta+\bar{\gamma}_{1},\zeta+2\bar{\gamma}_{1}}(\Phi(\be_{M_{m+1}\otimes M_{n}} \otimes e_{ii}))^{\halb} \bar{a}_{t}^{(j)}   f_{\zeta+\bar{\gamma}_{1},\zeta+2\bar{\gamma}_{1}}(\Phi(\be_{M_{m+1}\otimes M_{n}} \otimes e_{ii}))^{\halb})  \nonumber \\
& &>  14 \delta_{2} \label{ww4} .
\end{eqnarray}
We obtain
\begin{eqnarray}
\tau(\bar{a}_{t}^{(j)}) & \stackrel{\eqref{ww3}}{\le} & 2 \tau(f_{\zeta+\bar{\gamma}_{1},\zeta+2\bar{\gamma}_{1}}(\Phi(\be_{M_{m+1}\otimes M_{n}} \otimes e_{ii}))^{\halb} \bar{a}_{t}^{(j)}   \nonumber\\
&& f_{\zeta+\bar{\gamma}_{1},\zeta+2\bar{\gamma}_{1}}(\Phi(\be_{M_{m+1}\otimes M_{n}} \otimes e_{ii}))^{\halb})  \nonumber\\
&& + 2 \bar{\gamma}_{1} \nonumber \\
& \stackrel{\eqref{ww4}}{<} & 4 \tau(f_{\zeta+\bar{\gamma}_{1},\zeta+2\bar{\gamma}_{1}}(\Phi(\be_{M_{m+1}\otimes M_{n}} \otimes e_{ii}))^{\halb} \bar{a}_{t}^{(j)}  \nonumber \\
&&  f_{\zeta+\bar{\gamma}_{1},\zeta+2\bar{\gamma}_{1}}(\Phi(\be_{M_{m+1}\otimes M_{n}} \otimes e_{ii}))^{\halb}) \nonumber \\
& \le & 4 \tau((f_{\zeta+\bar{\gamma}_{1},\zeta+2\bar{\gamma}_{1}}(\Phi(\be_{M_{m+1}\otimes M_{n}} \otimes e_{ii}))^{\halb} \bar{a}_{t}^{(j)} \nonumber \\
&&  f_{\zeta+\bar{\gamma}_{1},\zeta+2\bar{\gamma}_{1}}(\Phi(\be_{M_{m+1}\otimes M_{n}} \otimes e_{ii}))^{\halb} - \delta_{2})_{+})  \nonumber\\
& & + 4 \delta_{2} \nonumber\\
& \stackrel{\eqref{ww4}}{\le} & 8 \tau((f_{\zeta+\bar{\gamma}_{1},\zeta+2\bar{\gamma}_{1}}(\Phi(\be_{M_{m+1}\otimes M_{n}} \otimes e_{ii}))^{\halb} \bar{a}_{t}^{(j)} \nonumber \\
&&  f_{\zeta+\bar{\gamma}_{1},\zeta+2\bar{\gamma}_{1}}(\Phi(\be_{M_{m+1}\otimes M_{n}} \otimes e_{ii}))^{\halb} - \delta_{2})_{+})  \nonumber\\
& \stackrel{\ref{order-0-almost-commuting}\mathrm{(iii)}}{\le} & 8 (m+1) n \tau((f_{\zeta+\bar{\gamma}_{1},\zeta+2\bar{\gamma}_{1}}(\Phi(d_{i_{0}i_{0}} \otimes e_{11} \otimes e_{ii}))^{\halb} \bar{a}_{t}^{(j)}  \label{www30} \\
&&  f_{\zeta+\bar{\gamma}_{1},\zeta+2\bar{\gamma}_{1}}(\Phi(d_{i_{0}i_{0}}\otimes e_{11} \otimes e_{ii}))^{\halb} - \delta_{2})_{+}) \nonumber \\
&& + 8 \delta_{2} \nonumber ,
\end{eqnarray}
so
\begin{eqnarray}
\tau(\bar{a}^{(j)}_{t}) & \stackrel{\eqref{ww3a}}{\le} & 2 (\tau(\bar{a}^{(j)}_{t})- 8 \delta_{2}) \nonumber \\
& \stackrel{\eqref{www30}}{\le} & 16 (m+1) n  \cdot \tau((f_{\zeta+\bar{\gamma}_{1},\zeta+2\bar{\gamma}_{1}}(\Phi(d_{i_{0}i_{0}} \otimes e_{11} \otimes e_{ii}))^{\halb} \bar{a}_{t}^{(j)}    \label{ww3b} \\
&& f_{\zeta+\bar{\gamma}_{1},\zeta+2\bar{\gamma}_{1}}(\Phi(d_{i_{0}i_{0}}\otimes e_{11} \otimes e_{ii}))^{\halb} - \delta_{2})_{+}). \nonumber 
\end{eqnarray}
Consequently, for $(j,t) \in J$ and $\tau \in T(A)$ we have
\begin{eqnarray}
\lefteqn{d_{\tau}(((\bar{a}^{(j)}_{t})^{\halb} (\be_{A} - \Phi(\be_{M_{m+1} \otimes M_{n} \otimes M_{2}})) (\bar{a}^{(j)}_{t})^{\halb} -\delta_{1})_{+})} \nonumber \\
& \stackrel{\ref{d-tau-notation}}{\le} & \frac{1}{\delta_{1}} \cdot \tau((\bar{a}^{(j)}_{t})^{\halb} (\be_{A} - \Phi(\be_{M_{m+1} \otimes M_{n} \otimes M_{2}})) (\bar{a}^{(j)}_{t})^{\halb} ) \nonumber \\
& \stackrel{\eqref{E-auxw1}}{\le} & \frac{1}{\delta_{1}} \cdot \bar{\beta} \nonumber \\
& \stackrel{\eqref{ww10a}}{<} & \frac{1}{k} \cdot \tau(\bar{a}^{(j)}_{t}) \nonumber \\
& \stackrel{\eqref{ww3b}}{<} & \frac{16 (m+1) n}{k}  \cdot \tau((f_{\zeta+\bar{\gamma}_{1},\zeta+2\bar{\gamma}_{1}}(\Phi(d_{i_{0}i_{0}} \otimes e_{11} \otimes e_{ii}))^{\halb} \bar{a}_{t}^{(j)}    \label{ww5a} \\
&& f_{\zeta+\bar{\gamma}_{1},\zeta+2\bar{\gamma}_{1}}(\Phi(d_{i_{0}i_{0}}\otimes e_{11} \otimes e_{ii}))^{\halb} - \delta_{2})_{+}). \nonumber
\end{eqnarray}
For convenience, set
\begin{eqnarray}
c_{t,i}^{(j)} & := & f_{\zeta+\bar{\gamma}_{1},\zeta+2\bar{\gamma}_{1}}(\Phi(d_{i_{0}i_{0}} \otimes e_{11} \otimes e_{ii}))^{\halb} \bar{a}_{t}^{(j)} \nonumber \\
&& f_{\zeta+\bar{\gamma}_{1},\zeta+2\bar{\gamma}_{1}}(\Phi(d_{i_{0}i_{0}}\otimes e_{11} \otimes e_{ii}))^{\halb} \label{E-auxw5}
\end{eqnarray}
and
\begin{eqnarray}
\bar{c}_{t,i}^{(j)}& := & (f_{\zeta+\bar{\gamma}_{1},\zeta+2\bar{\gamma}_{1}}(\Phi(d_{i_{0}i_{0}} \otimes e_{11} \otimes e_{ii}))^{\halb} \bar{a}_{t}^{(j)} \nonumber \\
&& f_{\zeta+\bar{\gamma}_{1},\zeta+2\bar{\gamma}_{1}}(\Phi(d_{i_{0}i_{0}}\otimes e_{11} \otimes e_{ii}))^{\halb} - \delta_{2})_{+}.  \label{E-auxw2}
\end{eqnarray}
Note that we have 
\begin{equation}
\label{E-auxw4}
d_{\tau}(((\bar{a}_{t}^{(j)})^{\halb} (\be_{A} - \Phi(\be)) (\bar{a}_{t}^{(j)})^{\halb} - \delta_{1})_{+}) \stackrel{\eqref{ww5a}}{<} \frac{16(m+1)n}{k} \tau(\bar{c}^{(j)}_{t,i})  \stackrel{\eqref{E-auxw3}}{<} \frac{1}{m+1}\tau(\bar{c}^{(j)}_{t,i})
\end{equation}
for $(j,t) \in J$, $i=1,2$ and $\tau \in T(A)$. 

By Corollary~\ref{cor-dr-comparison} and \eqref{E-auxw4}, for $(j,t) \in J$ there are 
\[
\bar{v}_{t}^{(j)} \in A
\]
such that
\begin{equation}
\label{ww5}
\|(\bar{v}^{(j)}_{t})^{*} \bar{v}_{t}^{(j)} - ((\bar{a}^{(j)}_{t})^{\halb} (\be_{A} - \Phi(\be_{M_{m+1}\otimes M_{n} \otimes M_{2}})) (\bar{a}^{(j)}_{t})^{\halb} - \delta_{1})_{+} \| < \delta_{1}
\end{equation}
and such that
\begin{equation}
\label{ww7}
\bar{v}^{(j)}_{t}(\bar{v}_{t}^{(j)})^{*} \in \her(\bar{c}_{t,i}^{(j)})
\end{equation}
for $i=1$ if $t$ is odd and $i=2$ if $t$ is even. We define
\begin{equation}
\label{ww21}
\bar{v}_{t}^{(j)} := 0 \mbox{ if } (j,t) \notin J.
\end{equation}
Set
\begin{equation}
\label{ww30}
\bar{v}^{(j)}:= \sum_{t=1}^{2\bar{t}} \alpha^{\halb}_{t} \bar{v}^{(j)}_{t}
\end{equation}
for $j=1, \ldots, s$. These are the desired elements, but we need some more preparation to verify conditions (i) through (iv) of the proposition.

For $t \neq t'$ both odd or both even and $(j,t), \, (j,t') \in J$ we have
\begin{eqnarray}
\lefteqn{\| ((\bar{v}^{(j)}_{t})^{*} \bar{v}^{(j)}_{t})^{\halb}  ((\bar{v}^{(j)}_{t'})^{*} \bar{v}^{(j)}_{t'})^{\halb} \| } \nonumber  \\
& = & \|   ((\bar{v}^{(j)}_{t'})^{*} \bar{v}^{(j)}_{t'})^{\halb}       ((\bar{v}^{(j)}_{t})^{*} \bar{v}^{(j)}_{t})           ((\bar{v}^{(j)}_{t'})^{*} \bar{v}^{(j)}_{t'})    ((\bar{v}^{(j)}_{t})^{*} \bar{v}^{(j)}_{t})    ((\bar{v}^{(j)}_{t'})^{*} \bar{v}^{(j)}_{t'})^{\halb} \|^{\frac{1}{4}}  \nonumber \\ 
& \stackrel{\eqref{ww5},\eqref{ww6}}{<} & (2 \delta_{1})^{\frac{1}{4}}. \label{ww23} 
\end{eqnarray}
Note that, by functional calculus, \eqref{ww7}, \eqref{E-auxw5} and \eqref{E-auxw2} we have 
\begin{equation}
\label{ww8}
g_{0, \delta_{2}} (c_{t,i}^{(j)}) \bar{v}^{(j)}_{t} =  \bar{v}^{(j)}_{t}
\end{equation}
for $(j,t) \in J$ (and  $i=1$ if $t$ is odd and $i=2$ if $t$ is even); by Proposition~\ref{order-0-almost-commuting} and our choice of $\bar{\beta}$, we have 
\begin{eqnarray}
\lefteqn{\| g_{0,\delta_{1}}(a_{t}^{(j)}) g_{0,\delta_{2}}(c_{t,i}^{(j)}) - g_{0,\delta_{2}}(c_{t,i}^{(j)}) \|} \nonumber \\
& \stackrel{\eqref{E-auxw5}}{=} & \|g_{0,\delta_{1}}(a_{t}^{(j)}) g_{0,\delta_{2}}(f_{\zeta+\bar{\gamma}_{1},\zeta+2\bar{\gamma}_{1}}(\Phi(d_{i_{0}i_{0}} \otimes e_{11} \otimes e_{ii}))^{\halb} \bar{a}_{t}^{(j)}  \nonumber \\
&& f_{\zeta+\bar{\gamma}_{1},\zeta+2\bar{\gamma}_{1}}(\Phi(d_{i_{0}i_{0}}\otimes e_{11} \otimes e_{ii}))^{\halb}) \nonumber \\
&& - g_{0,\delta_{2}}(f_{\zeta+\bar{\gamma}_{1},\zeta+2\bar{\gamma}_{1}}(\Phi(d_{i_{0}i_{0}} \otimes e_{11} \otimes e_{ii}))^{\halb} \bar{a}_{t}^{(j)} \nonumber \\
&& f_{\zeta+\bar{\gamma}_{1},\zeta+2\bar{\gamma}_{1}}(\Phi(d_{i_{0}i_{0}}\otimes e_{11} \otimes e_{ii}))^{\halb})\|  \nonumber \\
& \stackrel{\ref{order-0-almost-commuting}\mathrm{(iv)},\eqref{ww31}}{\le} &  \bar{\gamma}_{1}. \label{ww9}
\end{eqnarray}
We thus obtain 
\begin{equation}
\label{ww20}
\| g_{0,\delta_{1}}(a_{t	}^{(j)}) \bar{v}_{t}^{(j)} - \bar{v}_{t}^{(j)}\| \stackrel{\eqref{ww8},\eqref{ww9}}{<} \bar{\gamma}_{1}
\end{equation}
for any $(j,t) \in J$ (and hence for any $j,t$ by \eqref{ww21}).

Moreover, by \eqref{ww7}, \eqref{E-auxw2} and \eqref{ww15a} for all $(j,t) \in J$ we have
\begin{equation}
\label{ww11}
\bar{v}_{t}^{(j)} = g_{\zeta,\zeta+\bar{\gamma}_{1}}(\Phi(d_{i_{0}i_{0}} \otimes e_{11} \otimes e_{ii})) \bar{v}_{t}^{(j)}
\end{equation}
if $t$ is odd and $i=1$, and if $t$ is even and $i=2$. From this and \eqref{ww30} it follows that, for all $j=1, \ldots,s$, 
\begin{equation}
\label{E-auxw6}
\bar{v}^{(j)} = g_{\zeta,\zeta+\bar{\gamma}_{1}}(\Phi(d_{i_{0}i_{0}} \otimes e_{11} \otimes \be_{M_{2}}))  \bar{v}^{(j)},
\end{equation} 
so we have verified property (iv) of the proposition.

Next, we estimate 
\begin{eqnarray}
\lefteqn{\left\| (\bar{v}^{(j)})^{*}\bar{v}^{(j)} -  \sum_{t=1}^{2 \bar{t}} \alpha_{t} (\bar{v}_{t}^{(j)})^{*}   \bar{v}_{t}^{(j)}\right\|} \nonumber \\
& \stackrel{\eqref{ww30}}{=} & \left\| \sum_{t=1}^{2 \bar{t}} \alpha_{t}^{\halb} (\bar{v}_{t}^{(j)})^{*} \sum_{t=1}^{2\bar{t}} \alpha_{t}^{\halb}  \bar{v}_{t}^{(j)} \right. \nonumber \\
& & \left.- \sum_{t=1}^{2 \bar{t}} \alpha_{t} (\bar{v}_{t}^{(j)})^{*}   \bar{v}_{t}^{(j)}\right\| \nonumber \\
& \stackrel{\eqref{ww11}}{=}&  \left\| \left(\sum_{t \, \mathrm{odd}} \alpha_{t}^{\halb} (\bar{v}_{t}^{(j)})^{*} g_{\zeta,\zeta+\bar{\gamma}_{1}}(\Phi(d_{i_{0}i_{0}} \otimes e_{11} \otimes e_{11}))\right. \right. \nonumber \\
&& \left.+ \sum_{t \, \mathrm{even}} \alpha_{t}^{\halb} (\bar{v}_{t}^{(j)})^{*}  g_{\zeta,\zeta+\bar{\gamma}_{1}}(\Phi(d_{i_{0}i_{0}} \otimes e_{11} \otimes e_{22})) \right)\nonumber \\
&&  \left( \sum_{t \, \mathrm{odd}} \alpha_{t}^{\halb} g_{\zeta,\zeta+\bar{\gamma}_{1}}(\Phi(d_{i_{0}i_{0}} \otimes e_{11} \otimes e_{11}))\bar{v}_{t}^{(j)} \right.\nonumber \\
&& \left.+ \sum_{t \, \mathrm{even}} \alpha_{t}^{\halb} g_{\zeta,\zeta+\bar{\gamma}_{1}}(\Phi(d_{i_{0}i_{0}} \otimes e_{11} \otimes e_{22})) \bar{v}_{t}^{(j)}\right) \nonumber \\
& & \left.- \sum_{t=1}^{2 \bar{t}} \alpha_{t} (\bar{v}_{t}^{(j)})^{*}   \bar{v}_{t}^{(j)}\right\| \nonumber \\
& \stackrel{\eqref{E-auxw6},\eqref{ww23}}{\le} & 2 \bar{t} (\bar{t}-1) (2 \delta_{1})^{\frac{1}{4}}, \label{ww12}
\end{eqnarray}
where for the last estimate we have also used that $\Phi$ has order zero, whence 
\[
\Phi(d_{i_{0}i_{0}} \otimes e_{11} \otimes e_{11}) \perp \Phi(d_{i_{0}i_{0}} \otimes e_{11} \otimes e_{22}).
\]
At this point, one sees how the additional copy of $M_{2}$ is used to handle the `odd' and the `even' parts of $\bar{v}^{(j)}$ separately.  

We are now prepared to check property (i) of the proposition:
\begin{eqnarray*}
\lefteqn{\| (\bar{v}^{(j)})^{*} \bar{v}^{(j)} - \psi^{(j)}(f_{11}^{(j)})^{\halb} (\be_{A} - \Phi(\be)) \psi^{(j)}(f_{11}^{(j)})^{\halb} \|} \\
& \stackrel{\eqref{ww12},\ref{order-0-almost-commuting}}{=} & \left\|\sum_{ t \, \mathrm{odd}} \alpha_{t}  (\bar{v}^{(j)}_{t})^{*} \bar{v}^{(j)}_{t} + \sum_{ t \, \mathrm{even}} \alpha_{t}  (\bar{v}^{(j)}_{t})^{*} \bar{v}^{(j)}_{t} \right. \\
&& \left.- \psi^{(j)}(f_{11}^{(j)}) (\be_{A} - \Phi(\be))  \right\| + \bar{\gamma}_{1} + 2 \bar{t}^{2} (2 \delta_{1})^{\frac{1}{4}}\\
& \stackrel{\eqref{ww13}}{=}&   \left\|\sum_{ t \, \mathrm{odd}} \alpha_{t}  (\bar{v}^{(j)}_{t})^{*} \bar{v}^{(j)}_{t} + \sum_{ t \, \mathrm{even}} \alpha_{t}  (\bar{v}^{(j)}_{t})^{*} \bar{v}^{(j)}_{t} \right. \\
&& \left.- \left(\sum_{t \, \mathrm{odd}} \alpha_{t} a_{t}^{(j)} + \sum_{t \, \mathrm{even}} \alpha_{t} a_{t}^{(j)} \right)(\be_{A} - \Phi(\be))  \right\| + \bar{\gamma}_{1} + 2 \bar{t}^{2} (2 \delta_{1})^{\frac{1}{4}}\\
& \le & \left\| \sum_{t \, \mathrm{odd}} \alpha_{t} ((\bar{v}^{(j)}_{t})^{*} \bar{v}^{(j)}_{t} - a^{(j)}_{t} (\be_{A} - \Phi(\be))) \right\| \\
&& + \left\| \sum_{t \, \mathrm{even}} \alpha_{t} ((\bar{v}^{(j)}_{t})^{*} \bar{v}^{(j)}_{t} - a^{(j)}_{t} (\be_{A} - \Phi(\be))) \right\|  +   \bar{\gamma}_{1} + 2 \bar{t}^{2} (2 \delta_{1})^{\frac{1}{4}}\\
& \stackrel{\ref{order-0-almost-commuting}}{\le} & \left\| \sum_{t \, \mathrm{odd}} \alpha_{t} ((\bar{v}^{(j)}_{t})^{*} \bar{v}^{(j)}_{t} - (a^{(j)}_{t})^{\halb} (\be_{A} - \Phi(\be))  (a^{(j)}_{t})^{\halb}) \right\| \\
&& + \left\| \sum_{t \, \mathrm{even}} \alpha_{t} ((\bar{v}^{(j)}_{t})^{*} \bar{v}^{(j)}_{t} - (a^{(j)}_{t})^{\halb} (\be_{A} - \Phi(\be))  (a^{(j)}_{t})^{\halb}) \right\| \\
&&  +   \bar{\gamma}_{1} + 2 \bar{t}^{2} (2 \delta_{1})^{\frac{1}{4}} + 2 \bar{t} \bar{\gamma}_{1}\\
& \stackrel{\eqref{ww14}}{\le} &  \left\| \sum_{t \, \mathrm{odd}} \alpha_{t} ((\bar{v}^{(j)}_{t})^{*} \bar{v}^{(j)}_{t} - (\bar{a}^{(j)}_{t})^{\halb} (\be_{A} - \Phi(\be))  (\bar{a}^{(j)}_{t})^{\halb}) \right\| \\
&& + \left\| \sum_{t \, \mathrm{even}} \alpha_{t} ((\bar{v}^{(j)}_{t})^{*} \bar{v}^{(j)}_{t} - (\bar{a}^{(j)}_{t})^{\halb} (\be_{A} - \Phi(\be))  (\bar{a}^{(j)}_{t})^{\halb}) \right\| \\
&&  +   \bar{\gamma}_{1} + 2 \bar{t}^{2} (2 \delta_{1})^{\frac{1}{4}} + 2 \bar{t} \bar{\gamma}_{1} + 2 \bar{t} \cdot 2(2 \delta_{1})^{\halb}\\
& \le &  \left\| \sum_{t \, \mathrm{odd}} \alpha_{t} ((\bar{v}^{(j)}_{t})^{*} \bar{v}^{(j)}_{t} - ((\bar{a}^{(j)}_{t})^{\halb} (\be_{A} - \Phi(\be))  (\bar{a}^{(j)}_{t})^{\halb}- \delta_{1})_{+}) \right\| \\
&& + \left\| \sum_{t \, \mathrm{even}} \alpha_{t} ((\bar{v}^{(j)}_{t})^{*} \bar{v}^{(j)}_{t} - ((\bar{a}^{(j)}_{t})^{\halb} (\be_{A} - \Phi(\be))  (\bar{a}^{(j)}_{t})^{\halb}- \delta_{1})_{+}) \right\| \\
&&  +   \bar{\gamma}_{1} + 2 \bar{t}^{2} (2 \delta_{1})^{\frac{1}{4}} + 2 \bar{t} \bar{\gamma}_{1} + 2 \bar{t} \cdot 2(2 \delta_{1})^{\halb} + 2 \bar{t} \delta_{1}\\
& \stackrel{\eqref{ww5}}{\le} &  \bar{\gamma}_{1} + 2 \bar{t}^{2} (2 \delta_{1})^{\frac{1}{4}} + 2 \bar{t} \bar{\gamma}_{1} + 2 \bar{t} \cdot 2(2 \delta_{1})^{\halb} + 2 \bar{t} \delta_{1} + 2 \bar{t} \delta_{1}\\
& \stackrel{\eqref{ww15a},\eqref{ww16}}{<} & \bar{\gamma}_{0}.
\end{eqnarray*} 

To confirm (iii), note that by \eqref{ww15a} and \eqref{ww2} for any $j$ we have
\begin{equation}
\label{ww19}
g_{0,\bar{\gamma}_{1}}(\psi^{(j)}(f_{11}^{(j)})) a^{(j)}_{t}= a_{t}^{(j)}
\end{equation} 
if $t \ge 2$, so 
\begin{eqnarray*}
\lefteqn{\| \bar{v}^{(j)} - g_{0,\bar{\gamma}_{1}}(\psi^{(j)}(f_{11}^{(j)})) \bar{v}^{(j)} \|   }\\
& \stackrel{\eqref{ww30}}{=} & \left\| \sum_{t=1}^{2\bar{t}} \alpha_{t}^{\halb} (\bar{v}_{t}^{(j)} - g_{0,\bar{\gamma}_{1}}(\psi^{(j)}(f_{11}^{(j)})) \bar{v}_{t}^{(j)}) \right\| \\
& \le & \alpha_{1}^{\halb} + (2 \bar{t} - 1) \cdot \max_{t \ge 2} \| \bar{v}_{t}^{(j)} - g_{0,\bar{\gamma}_{1}}(\psi^{(j)}(f_{11}^{(j)})) \bar{v}_{t}^{(j)} \| \\
& \stackrel{\eqref{ww20}}{ \le} &  \alpha_{1}^{\halb} + (2 \bar{t} - 1) \cdot (\max_{t \ge 2} \| \bar{v}_{t}^{(j)} - g_{0,\bar{\gamma}_{1}}(\psi^{(j)}(f_{11}^{(j)}))  g_{0,\delta_{1}}(a_{t}^{(j)}) \bar{v}_{t}^{(j)} \|  + \bar{\gamma}_{1} )\\
& \stackrel{\eqref{ww19}}{=} &  \alpha_{1}^{\halb} + (2 \bar{t} - 1) \cdot (\max_{t \ge 2} \| \bar{v}_{t}^{(j)} -   g_{0,\delta_{1}}(a_{t}^{(j)}) \bar{v}_{t}^{(j)} \|  + \bar{\gamma}_{1} )\\
& \stackrel{\eqref{ww20}}{\le} &  \alpha_{1}^{\halb} + (2 \bar{t} - 1) \cdot (\bar{\gamma}_{1}  + \bar{\gamma}_{1} )\\
& \stackrel{\eqref{ww15a},\eqref{ww16},\eqref{ww22}}{<} & \bar{\gamma}_{0}.
\end{eqnarray*}
 
 It remains to verify (ii). To this end, we define step functions $h_{\mathrm{odd}}$ and $h_{\mathrm{even}}$ on $[0,1]$ by 
\[
h_{\mathrm{odd}}(x):=
\frac{2l-1}{2\bar{t}}  \mbox{ if } x \in \left(\frac{2l-2}{2 \bar{t}}, \frac{2l}{2\bar{t}}\right] \cap [0,1] \mbox{ for } l=0,\ldots,\bar{t} \] 
and 
\[
h_{\mathrm{even}}(x):= \frac{2l}{2\bar{t}}  \mbox{ if } x \in \left(\frac{2l-1}{2 \bar{t}}, \frac{2l+1}{2\bar{t}}\right] \cap [0,1] \mbox{ for } l=0,\ldots,\bar{t},
\]
and, for $j=1, \ldots, s$, elements
\[
a^{(j)}_{\mathrm{odd}} := h_{\mathrm{odd}}(\psi^{(j)}(f_{11}^{(j)})), \, a^{(j)}_{\mathrm{even}} := h_{\mathrm{even}}(\psi^{(j)}(f_{11}^{(j)})) \in A''
\] 
of the envelopping von Neumann algebra of $A$. Note that for each $j$ we have
\begin{equation}
\label{ww24}
\| a^{(j)}_{\mathrm{odd}} - \psi^{(j)}(f_{11}^{(j)})\| \le \frac{1}{2\bar{t}} \mbox{ and } \| a^{(j)}_{\mathrm{even}} - \psi^{(j)}(f_{11}^{(j)})\| \le \frac{1}{2\bar{t}}
\end{equation}
in $A''$, and that 
\begin{equation}
\label{ww25}
a^{(j)}_{t} a^{(j)}_{\mathrm{odd}} = \alpha_{t} a_{t}^{(j)} \mbox{for } t \mbox{odd and } a^{(j)}_{t} a^{(j)}_{\mathrm{even}} = \alpha_{t} a_{t}^{(j)} \mbox{for } t \mbox{ even};
\end{equation} 
the same equalities hold with $\bar{a}^{(j)}_{t}$ in place of $a^{(j)}_{t}$. 

For $t$ odd we then have
\begin{eqnarray}
\lefteqn{\|((\bar{v}_{t}^{(j)})^{*}\bar{v}_{t}^{(j)})^{\halb}(a_{\mathrm{odd}}^{(j)} - \alpha_{t} \cdot \be_{A}) \|  } \nonumber \\
& \le & \| (a_{\mathrm{odd}}^{(j)} - \alpha_{t} \cdot \be_{A}) (\bar{v}_{t}^{(j)})^{*} \bar{v}_{t}^{(j)} (a_{\mathrm{odd}}^{(j)} - \alpha_{t} \cdot \be_{A})\|^{\halb} \nonumber \\
& \stackrel{\eqref{ww5}}{<} & ( \| (a^{(j)}_{\mathrm{odd}} -\alpha_{t} \cdot \be_{A}) ((\bar{a}_{t}^{(j)})^{\halb} (\be_{A} - \Phi(\be))  (\bar{a}_{t}^{(j)})^{\halb} - \delta_{1})  (a^{(j)}_{\mathrm{odd}} -\alpha_{t} \cdot \be_{A})  \| \nonumber \\
&& + \delta_{1}  )^{\halb} \nonumber \\
& \stackrel{\eqref{ww25}}{=} & \delta_{1}^{\halb},  \label{ww26}
\end{eqnarray} 
and the respective estimate holds for $t$ even: 
\begin{equation}
\label{ww27}
\|((\bar{v}_{t}^{(j)})^{*}\bar{v}_{t}^{(j)})^{\halb}(a_{\mathrm{even}}^{(j)} - \alpha_{t} \cdot \be_{A}) \| < \delta_{1}^{\halb}.
\end{equation}
Note that we can use \eqref{ww5} only if $(j,t) \in J$, but otherwise $\bar{v}_{t}^{(j)} = 0$ and the estimates are trivial.

Now, let $\hat{v}_{t}^{(j)} \in A''$ denote the partial isometry from the polar decomposition of $\bar{v}^{(j)}_{t}$, so that 
\begin{equation}
\label{ww28a}
\bar{v}_{t}^{(j)} = \hat{v}^{(j)}_{t}((\bar{v}^{(j)}_{t})^{*} \bar{v}^{(j)}_{t})^{\halb} 
\end{equation}
for all $j,t$. For each $j=1, \ldots,s$ we have 
\begin{eqnarray}
\lefteqn{\left\| \sum_{t \, \mathrm{odd}} \alpha_{t}^{\halb} \bar{v}_{t}^{(j)} \psi^{(j)}(f_{11}^{(j)}) - \alpha_{t}^{\frac{3}{2}} \bar{v}_{t}^{(j)} \right\|  } \nonumber \\
& \stackrel{\eqref{ww24},\eqref{ww29}}{\le} & \left\| \sum_{t \, \mathrm{odd}} \alpha_{t}^{\halb} \bar{v}_{t}^{(j)} a^{(j)}_{\mathrm{odd}} - \alpha_{t}^{\frac{3}{2}} \bar{v}_{t}^{(j)} \right\| + \frac{2}{\bar{t}} \nonumber \\
& \stackrel{\eqref{ww28a}}{=} & \left\|  \sum_{t \, \mathrm{odd}} \alpha_{t}^{\halb} \hat{v}^{(j)}_{t} ((\bar{v}^{(j)}_{t})^{*} \bar{v}^{(j)}_{t})^{\halb} (a^{(j)}_{\mathrm{odd}} - \alpha_{t} \cdot \be_{A}) \right\| + \frac{2}{\bar{t}} \nonumber \\
& \stackrel{\eqref{ww26}}{\le} & \bar{t} \delta_{1}^{\halb} + \frac{2}{\bar{t}}. \label{ww28}
\end{eqnarray} 
Here, we have used the very rough estimates
\begin{equation}
\label{ww29}
\left\| \sum_{t \, \mathrm{odd}} \alpha_{t}^{\halb} \bar{v}^{(j)}_{t}\right\| , \;  \left\| \sum_{t \, \mathrm{even}} \alpha_{t}^{\halb} \bar{v}^{(j)}_{t}\right\|  < 4,
\end{equation}
which follow directly from (i) and \eqref{ww30}. (In fact, the norms in \eqref{ww29} are at most little more than 1 by our construction, but the rough estimates above will do for our purposes.) 

Using \eqref{ww27} in place of \eqref{ww26}, we see that \eqref{ww28} also holds for the sum over even values of $t$, so that we have
\begin{eqnarray}
\lefteqn{ \left\| \bar{v}^{(j)} \psi^{(j)}(f_{11}^{(j)}) - \sum_{t=1}^{2\bar{t}} \alpha^{\frac{3}{2}}_{t} \bar{v}_{t}^{(j)} \right\| } \nonumber \\
& \stackrel{\eqref{ww30}}{=} & \left\| \sum_{t=1}^{2\bar{t}}\alpha_{t}^{\halb} \bar{v}_{t}^{(j)} \psi^{(j)}(f_{11}^{(j)}) -  \alpha^{\frac{3}{2}}_{t} \bar{v}_{t}^{(j)} \right\| \nonumber  \\
& \stackrel{\eqref{ww28}}{<} & 2(\bar{t} \delta_{1}^{\halb} + \frac{2}{\bar{t}}) \nonumber \\
& \stackrel{\eqref{ww22},\eqref{ww16}}{<} & \frac{\bar{\gamma}_{0}}{2}. \label{ww31a}
\end{eqnarray}
On the other hand, we compute
\begin{eqnarray}
\lefteqn{\left\| \psi^{(j)}(f_{11}^{(j)}) \bar{v}^{(j)} - \sum_{t=1}^{2\bar{t}} \alpha_{t}^{\frac{3}{2}} \bar{v}_{t}^{(j)} \right\| } \nonumber \\
& \stackrel{\eqref{ww30},\eqref{ww24},\eqref{ww29}}{\le} & \left\| a^{(j)}_{\mathrm{odd}} \sum_{t \, \mathrm{odd}} \alpha_{t}^{\halb} \bar{v}_{t}^{(j)} - \sum_{t \, \mathrm{odd}} \alpha_{t}^{\frac{3}{2}} \bar{v}_{t}^{(j)} \right\| \nonumber \\
&& +  \left\| a^{(j)}_{\mathrm{even}} \sum_{t \, \mathrm{odd}} \alpha_{t}^{\halb} \bar{v}_{t}^{(j)} - \sum_{t \, \mathrm{even}} \alpha_{t}^{\frac{3}{2}} \bar{v}_{t}^{(j)} \right\| + \frac{4}{\bar{t}} \nonumber \\
& \stackrel{\eqref{ww20}}{\le} &  \left\| a^{(j)}_{\mathrm{odd}} \sum_{t \, \mathrm{odd}} \alpha_{t}^{\halb} g_{0,\delta_{1}}(a^{(j)}_{t}) \bar{v}_{t}^{(j)} - \sum_{t \, \mathrm{odd}} \alpha_{t}^{\frac{3}{2}} \bar{v}_{t}^{(j)} \right\| \nonumber \\
&&  + \left\| a^{(j)}_{\mathrm{even}} \sum_{t \, \mathrm{odd}} \alpha_{t}^{\halb} g_{0,\delta_{1}}(a^{(j)}_{t}) \bar{v}_{t}^{(j)} - \sum_{t \, \mathrm{even}} \alpha_{t}^{\frac{3}{2}} \bar{v}_{t}^{(j)} \right\| + \frac{4}{\bar{t}} + 2 \bar{t} \bar{\gamma}_{1}  \nonumber \\
& \stackrel{\eqref{ww25}}{=} & \left\|  \sum_{t \, \mathrm{odd}} \alpha_{t}^{\frac{3}{2}} g_{0,\delta_{1}}(a^{(j)}_{t}) \bar{v}_{t}^{(j)} - \sum_{t \, \mathrm{odd}} \alpha_{t}^{\frac{3}{2}} \bar{v}_{t}^{(j)} \right\| \nonumber \\
&&  + \left\|  \sum_{t \, \mathrm{odd}} \alpha_{t}^{\frac{3}{2}} g_{0,\delta_{1}}(a^{(j)}_{t}) \bar{v}_{t}^{(j)} - \sum_{t \, \mathrm{even}} \alpha_{t}^{\frac{3}{2}} \bar{v}_{t}^{(j)} \right\| + \frac{4}{\bar{t}} + 2 \bar{t} \bar{\gamma}_{1}  \nonumber \\
& \stackrel{\eqref{ww20}}{\le} & \frac{4}{\bar{t}} + 2 \bar{t} \bar{\gamma}_{1} + 2 \bar{t} \bar{\gamma}_{1} \nonumber \\
& \stackrel{\eqref{ww15a},\eqref{ww15b},\eqref{ww16},\eqref{ww22}}{<} & \frac{\bar{\gamma}_{0}}{2}. \label{ww32} 
\end{eqnarray}
Combining \eqref{ww31a} and \eqref{ww32}, we obtain property (i) of the proposition. We are done.
\end{nproof}
\en

\bn
\label{E}
\begin{nprop}
Let $A$ be a simple, separable, unital $\mathrm{C}^{*}$-algebra with $\dr A = m < \infty$. Let $F$ be a finite dimensional $\mathrm{C}^{*}$-algebra, $n \in \N$ and
\[
\psi: F \to A
\]
a c.p.c.\ order zero map. For any $0<\theta$ and $0< \zeta < 1$ there is $\beta >0$ such that the following holds: 

If $ i_{0} \in \{0, \ldots,m\}$ and 
\[
\Phi: M_{m+1} \otimes M_{n} \otimes M_{2} \to A
\]
is a c.p.c.\ order zero map satisfying
\begin{equation}
\label{c}
\| [\psi(x), \Phi(y)]\| \le \beta \|x\| \|y\| \mbox{ for any }  x \in F \mbox{ and } y  \in M_{m+1} \otimes M_{n} \otimes M_{2}
\end{equation}
and
\[
\tau(\Phi(\be_{M_{m+1} \otimes M_{n} \otimes M_{2}})) > 1 - \beta \; \forall \, \tau \in T(A),
\]
then there is $v \in A$ of norm at most one such that
\[
\|v^{*}v - (\be_{A}- \Phi(\be_{M_{m+1} \otimes M_{n} \otimes M_{2}}))^{\frac{1}{2}} \psi(\be_{F})  (\be_{A}- \Phi(\be_{M_{m+1} \otimes M_{n} \otimes M_{2}}))^{\frac{1}{2}}\| < \theta,
\]
\[
v v^{*} \in \overline{(\Phi(d_{i_{0}i_{0}} \otimes e_{11} \otimes \be_{M_{2}})- \zeta)_{+}A(\Phi(d_{i_{0}i_{0}} \otimes e_{11} \otimes \be_{M_{2}})- \zeta)_{+}}
\]
(where $\{d_{kl} \mid k,l=0, \ldots,m\}$  and $\{e_{kl} \mid k,l=1, \ldots, n\}$ denote  sets of matrix units for $M_{m}$ and $M_{n}$, respectively) and
\[
\|[\psi(x),v]\| \le \theta \|x\| \mbox{ for all }  x \in F.
\]
\end{nprop}

\begin{nproof}
Let $F  = M_{r^{(1)}} \oplus \ldots \oplus M_{r^{(s)}}$ and let $\psi^{(j)}: M_{r^{(j)}} \to A$ denote the $j$-th component of $\psi$. Let 
\begin{equation}
\label{wwE5}
\pi^{(j)}:M_{r^{(j)}} \to A'' 
\end{equation}
be the canonical supporting $*$-homomorphisms. 

Choose $\bar{\gamma}_{0}>0$ such that 
\begin{equation}
\label{wwE7}
10 s \max_{j} \{r^{(j)}\} \bar{\gamma}_{0} < \theta.
\end{equation}
Choose $0<\bar{\gamma}_{1}$ and $0<\bar{\beta}$ as in Proposition~\ref{E-aux}. 

Using Proposition~\ref{order-0-almost-commuting}, we may choose 
\begin{equation}
\label{wwE1}
0<\beta<\bar{\beta}
\end{equation}  
such that, if \eqref{c}  holds for some $\Phi$ as in \ref{E}, then 
\begin{equation}
\label{wwE4}
\| [\pi^{(j)}(f_{1k}^{(j)}) g_{0,\bar{\gamma}_{1}}(\psi^{(j)}(\be_{M_{r^{(j)}}})) ,  g_{\zeta,\zeta +\bar{\gamma}_{1}}(\Phi(d_{ii} \otimes e_{11} \otimes \be_{M_{2}}))^{2}] \| < \bar{\gamma}_{0},
\end{equation}
\begin{equation}
\label{wwE6}
\| [\pi^{(j)}(f_{1k}^{(j)}) \psi^{(j)}(\be_{M_{r^{(j)}}})^{\halb} ,  (\be_{A} - \Phi(\be))^{\halb}] \| < \bar{\gamma}_{0},
\end{equation}
and
\begin{equation}
\label{wwE8}
\| [\pi^{(j)}(x^{(j)}) \psi^{(j)}(\be_{M_{r^{(j)}}}) ,  g_{\zeta,\zeta +\bar{\gamma}_{1}}(\Phi(d_{ii} \otimes e_{11} \otimes \be_{M_{2}}))] \| < \bar{\gamma}_{0} \|x\|
\end{equation}
for $i \in \{0, \ldots, m\}$, $j\in \{1,\ldots, s\}$, $k \in \{1, \ldots, r^{(j)}\}$ and $x = x^{(1)} \oplus \ldots \oplus x^{(s)} \in F$.

Now, suppose there is a map
\[
\Phi: M_{m+1} \otimes M_{n} \otimes M_{2} \to A
\]
as in Proposition~\ref{E}. By \eqref{wwE1}, this map will also satisfy the hypotheses of Proposition~\ref{E-aux}, so there are 
\[
\bar{v}^{(j)} \in A, \, j=1, \ldots, s,
\]
as in \ref{E-aux}. For $j=1, \ldots, s$ set
\begin{equation}
\label{wwE2}
v^{(j)} := \sum_{k=1}^{r^{(j)}} \pi^{(j)}(f_{k1}^{(j)}) \bar{v}^{(j)} \pi^{(j)}(f_{1k}^{(j)})
\end{equation}
and
\begin{equation}
\label{wwE3}
v:= g_{\zeta,\zeta+\bar{\gamma}_{1}}(\Phi(d_{i_{0}i_{0}} \otimes e_{11} \otimes \be_{M_{2}})) \sum_{j=1}^{s} v^{(j)}.
\end{equation}
Note that, for $j=1, \ldots, s$, 
\begin{eqnarray}
\lefteqn{\| v^{(j)} - g_{0,\bar{\gamma}_{1}}(\psi^{(j)}(\be_{M_{r^{(j)}}})) v^{(j)} \|} \nonumber \\
& \stackrel{\eqref{order-zero-facts-w1},\eqref{wwE2}}{=} & \left\| \sum_{k=1}^{r^{(j)}} \pi^{(j)} (f_{k1}^{(j)})(\bar{v}^{(j)} - g_{0,\bar{\gamma}_{1}}(\psi^{(j)}(f_{11}^{(j)})) \bar{v}^{(j)})  \pi^{(j)} (f_{1k}^{(j)}) \right\| \nonumber \\
& = & \| \bar{v}^{(j)} -  g_{0,\bar{\gamma}_{1}}(\psi^{(j)}(f_{11}^{(j)})) \bar{v}^{(j)} \| \nonumber \\
& \stackrel{\ref{E-aux}\mathrm{(iii)}}{<} & \bar{\gamma}_{0}. \label{a} 
\end{eqnarray}
We have 
\[
vv^{*} \in \her((\Phi(d_{i_{0}i_{0}} \otimes e_{11} \otimes \be_{M_{2}}) - \zeta)_{+})
\]
by \eqref{wwE3}.

We now compute
{
\allowdisplaybreaks
\begin{eqnarray*}
\lefteqn{\| v^{*}v - (\be_{A}- \Phi(\be))^{\halb} \psi(\be_{F}) (\be_{A}- \Phi(\be))^{\halb} \|  }  \\
& \stackrel{\eqref{wwE3}}{=} & \left\| \sum_{j=1}^{s} (v^{(j)})^{*} g_{\zeta,\zeta+\bar{\gamma}_{1}}(\Phi(d_{i_{0}i_{0}} \otimes e_{11} \otimes \be_{M_{2}}))^{2} \sum_{j=1}^{s} v^{(j)} \right. \\
&&\left. - (\be_{A}- \Phi(\be))^{\halb} \psi(\be_{F}) (\be_{A}- \Phi(\be))^{\halb} \right\|  \\
& \stackrel{\eqref{a}}{\le} &  \left\| \sum_{j=1}^{s} (v^{(j)})^{*} g_{0,\bar{\gamma}_{1}}(\psi^{(j)}(\be_{M_{r^{(j)}}}))   g_{\zeta,\zeta+\bar{\gamma}_{1}}(\Phi(d_{i_{0}i_{0}} \otimes e_{11} \otimes \be_{M_{2}}))^{2} \right.\\
&& \sum_{j=1}^{s}   g_{0,\bar{\gamma}_{1}}(\psi^{(j)}(\be_{M_{r^{(j)}}}))  v^{(j)}  \\
&& \left.- (\be_{A}- \Phi(\be))^{\halb} \psi(\be_{F}) (\be_{A}- \Phi(\be))^{\halb} \right\|   + 2s \bar{\gamma}_{0} \\
& \stackrel{\eqref{wwE2}}{=} &  \left\| \sum_{j=1}^{s} (\sum_{k=1}^{r^{(j)}} \pi^{(j)}(f_{k1}^{(j)}) (\bar{v}^{(j)})^{*} \pi^{(j)}(f_{1k}^{(j)})) g_{0,\bar{\gamma}_{1}}(\psi^{(j)}(\be_{M_{r^{(j)}}})) \right. \\
&& g_{\zeta,\zeta+\bar{\gamma}_{1}}(\Phi(d_{i_{0}i_{0}} \otimes e_{11} \otimes \be_{M_{2}}))^{2} \\
&& \sum_{j=1}^{s}     \sum_{k=1}^{r^{(j)}} g_{0,\bar{\gamma}_{1}}(\psi^{(j)}(\be_{M_{r^{(j)}}})) \pi^{(j)}(f_{k1}^{(j)}) \bar{v}^{(j)} \pi^{(j)}(f_{1k}^{(j)})  \\
&& \left. - (\be_{A}- \Phi(\be))^{\halb} \psi(\be_{F}) (\be_{A}- \Phi(\be))^{\halb} \right\|  \\
&& + 2s \bar{\gamma}_{0} \\
& \stackrel{\eqref{wwE4}}{\le} &  \left\| \sum_{j=1}^{s} (\sum_{k=1}^{r^{(j)}} \pi^{(j)}(f_{k1}^{(j)}) (\bar{v}^{(j)})^{*}  g_{\zeta,\zeta+\bar{\gamma}_{1}}(\Phi(d_{i_{0}i_{0}} \otimes e_{11} \otimes \be_{M_{2}}))^{2}\right. \\
&& \pi^{(j)}(f_{1k}^{(j)})) g_{0,\bar{\gamma}_{1}}(\psi^{(j)}(\be_{M_{r^{(j)}}}))  \\
&& \sum_{j=1}^{s}     \sum_{k=1}^{r^{(j)}} g_{0,\bar{\gamma}_{1}}(\psi^{(j)}(\be_{M_{r^{(j)}}})) \pi^{(j)}(f_{k1}^{(j)}) \bar{v}^{(j)} \pi^{(j)}(f_{1k}^{(j)})  \\
&& \left.- (\be_{A}- \Phi(\be))^{\halb} \psi(\be_{F}) (\be_{A}- \Phi(\be))^{\halb} \right\|  \\
&& + 2s \bar{\gamma}_{0} + s \bar{\gamma}_{0} \cdot \max_{j} \{r^{(j)}\} \\
& \stackrel{\ref{E-aux}\mathrm{(iv)}}{=} &  \left\| \sum_{j=1}^{s} (\sum_{k=1}^{r^{(j)}} \pi^{(j)}(f_{k1}^{(j)}) (\bar{v}^{(j)})^{*}   \pi^{(j)}(f_{1k}^{(j)})) g_{0,\bar{\gamma}_{1}}(\psi^{(j)}(\be_{M_{r^{(j)}}})) \right. \\
&& \sum_{j=1}^{s}     \sum_{k=1}^{r^{(j)}} g_{0,\bar{\gamma}_{1}}(\psi^{(j)}(\be_{M_{r^{(j)}}})) \pi^{(j)}(f_{k1}^{(j)}) \bar{v}^{(j)} \pi^{(j)}(f_{1k}^{(j)})  \\
&& \left.- (\be_{A}- \Phi(\be))^{\halb} \psi(\be_{F}) (\be_{A}- \Phi(\be))^{\halb} \right\|  \\
&& + 2s \bar{\gamma}_{0} + s \bar{\gamma}_{0} \cdot \max_{j} \{r^{(j)}\} \\
& \stackrel{\eqref{order-zero-facts-w1},\eqref{wwE5}}{=} & \left\| \sum_{j=1}^{s} (\sum_{k=1}^{r^{(j)}} \pi^{(j)}(f_{k1}^{(j)}) (\bar{v}^{(j)})^{*} g_{0,\bar{\gamma}_{1}}(\psi^{(j)}(f_{11}^{(j)}))^{2}   \bar{v}^{(j)} \pi^{(j)}(f_{1k}^{(j)}) \right. \\
&& \left.- (\be_{A}- \Phi(\be))^{\halb} \psi(\be_{F}) (\be_{A}- \Phi(\be))^{\halb} \right\|  \\
&& + 2s \bar{\gamma}_{0} + s \bar{\gamma}_{0} \cdot \max_{j} \{r^{(j)}\}  \\
& \stackrel{\ref{E-aux}\mathrm{(iii)}}{\le} & \left\| \sum_{j=1}^{s} (\sum_{k=1}^{r^{(j)}} \pi^{(j)}(f_{k1}^{(j)}) (\bar{v}^{(j)})^{*}    \bar{v}^{(j)} \pi^{(j)}(f_{1k}^{(j)}) \right.  \\
&& \left.- (\be_{A}- \Phi(\be))^{\halb} \psi(\be_{F}) (\be_{A}- \Phi(\be))^{\halb} \right\|  \\
&& + 2s \bar{\gamma}_{0} + s \bar{\gamma}_{0} \cdot \max_{j} \{r^{(j)}\} +  2s \bar{\gamma}_{0} \cdot \max_{j} \{r^{(j)}\}\\
& \stackrel{\eqref{wwE5}}{=} & \left\| \sum_{j=1}^{s} (\sum_{k=1}^{r^{(j)}} \pi^{(j)}(f_{k1}^{(j)}) (\bar{v}^{(j)})^{*}    \bar{v}^{(j)} \pi^{(j)}(f_{1k}^{(j)}) \right.  \\
&& \left.- (\be_{A}- \Phi(\be))^{\halb} \pi^{(j)}(f_{k1}^{(j)})    \psi^{(j)}(f_{11}^{(j)}) \pi^{(j)}(f_{1k}^{(j)})  (\be_{A}- \Phi(\be))^{\halb}) \right\|  \\
&& + 2s \bar{\gamma}_{0} + s \bar{\gamma}_{0} \cdot \max_{j} \{r^{(j)}\}  + 2s \bar{\gamma}_{0} \cdot \max_{j} \{r^{(j)}\}\\
& \stackrel{\eqref{wwE6}}{\le} & \left\| \sum_{j=1}^{s} \left(\sum_{k=1}^{r^{(j)}} \pi^{(j)}(f_{k1}^{(j)}) ((\bar{v}^{(j)})^{*}    \bar{v}^{(j)}\right. \right. \\
&& \left. \left. - \psi^{(j)}(f_{11}^{(j)})^{\halb} (\be_{A}- \Phi(\be))      \psi^{(j)}(f_{11}^{(j)})^{\halb} ) \pi^{(j)}(f_{1k}^{(j)}) \right) \right\|  \\
&& + \bar{\gamma}_{0}(2s  + 5s  \cdot \max_{j} \{r^{(j)}\}) \\
& \stackrel{\ref{E-aux}\mathrm{(i)}}{\le} &   \bar{\gamma}_{0}(2s  + 6s  \cdot \max_{j} \{r^{(j)}\}) \\
& \stackrel{\eqref{wwE7}}{<} & \theta.
\end{eqnarray*}
}

Finally, using that $\psi$ has order zero, for $x \in F= M_{r^{(1)}} \oplus \ldots \oplus M_{r^{(s)}}$ we check
{
\allowdisplaybreaks
\begin{eqnarray*}
\lefteqn{\| [\psi(x),v]\|} \\
& \stackrel{\eqref{order-zero-facts-w1},\eqref{wwE5}}{=} & \left\| \sum_{j=1}^{s} \pi^{(j)}(x^{(j)}) \psi^{(j)}(\be_{M_{r^{(j)}}}) v - v \psi^{(j)}(\be_{M_{r^{(j)}}})  \pi^{(j)}(x^{(j)}) \right\| \\
& \stackrel{\eqref{wwE3},\eqref{wwE8}}{\le} & \left\|  g_{\zeta,\zeta+\bar{\gamma}_{1}}(\Phi(d_{i_{0}i_{0}} \otimes e_{11} \otimes \be_{M_{2}})) \right.\\
&& \left. \left(\sum_{j=1}^{s} \pi^{(j)}(x^{(j)}) \psi^{(j)}(\be_{M_{r^{(j)}}}) v^{(j)} - v^{(j)}  \psi^{(j)}(\be_{M_{r^{(j)}}}) \pi^{(j)}(x^{(j)}) \right) \right\| \\
& & + s \bar{\gamma}_{0} \|x\| \\
& \stackrel{\ref{order-zero-facts}}{\le} &  \max_{j} \left\| \pi^{(j)}(x^{(j)}) \sum_{k=1}^{r^{(j)}} \pi^{(j)}(f^{(j)}_{k1}) \psi^{(j)}(f_{11}^{(j)}) \pi^{(j)}(f_{11}^{(j)}) \bar{v}^{(j)} \pi^{(j)}(f_{1k}^{(j)})  \right. \\
&& \left.- \pi^{(j)}(f^{(j)}_{k1})  \bar{v}^{(j)}  \pi^{(j)}(f_{11}^{(j)})  \psi^{(j)}(f_{11}^{(j)}) \pi^{(j)}(f_{1k}^{(j)}) \right\| \\
& & + s \bar{\gamma}_{0} \|x\| \\
& \stackrel{\eqref{wwE5}}{\le} & \|x\| \cdot \max_{j} \|\psi^{(j)}(f_{11}^{(j)}) \bar{v}^{(j)} -  \bar{v}^{(j)} \psi^{(j)}(f_{11}^{(j)})\| \\
& & + s \bar{\gamma}_{0} \|x\| \\
& \stackrel{\ref{E-aux}\mathrm{(ii)}}{\le} &  (\max_{j}\{r^{(j)}\} + s) \bar{\gamma}_{0}  \|x\| \\
& \stackrel{\eqref{wwE7}}{<} & \theta \|x\|.
\end{eqnarray*}
}
This completes the proof.
\end{nproof}
\en

\bn
\label{H}
\begin{nprop}
Let $A$ be a simple, separable, unital $\mathrm{C}^{*}$-algebra with $\dr A = m < \infty$. Given a finite subset $\Fh \subset A$, a positive normalized function $\bar{h} \in \mathcal{C}_{0}((0,1])$ and $\delta>0$, there are a finite subset $\Gh \subset A$ and $\alpha > 0$ such that the following holds:

Suppose $(F=F^{(0)} \oplus \ldots \oplus F^{(m)},\sigma,\varrho)$ is an $m$-decomposable c.p.c.\ approximation for $\Gh$ to within  $\alpha$. For $i=0, \ldots, m$, let $v_{i} \in A$ be normalized elements satisfying
\begin{equation}
\label{w5}
\|[\varrho^{(i)}(x),v_{i}]\| \le \alpha \|x\| \mbox{ for all } x \in F^{(i)}.
\end{equation}
Then, 
\begin{equation}
\label{wH5}
v:= \sum_{i=0}^{m} v_{i} \bar{h}(\varrho^{(i)}(\be_{F^{(i)}}))
\end{equation}
satisfies
\[
\|[v,a]\| < \delta \mbox{ for all } a \in \Fh.
\]
\end{nprop}

\begin{nproof}
We may assume the elements of $\mathcal{F}$ to be positive and normalized. For convenience, we set
\begin{equation}
\label{wH4}
\bar{\delta}:= \frac{\delta}{9(m+1)}.
\end{equation}
Choose $\tilde{h} \in \Ch_{0}((0,1])$ such that
\begin{equation}
\label{w2}
\| \id_{(0,1]} \cdot \tilde{h} - \bar{h}\| < \bar{\delta}
\end{equation}
Using Lemma~\ref{multiplicative-domain}, we find $\alpha>0$ and a finite subset $\mathcal{G} \subset A$ such that, whenever $(F,\sigma,\varrho)$ is a c.p.c.\ approximation of $\mathcal{G}$ to within  $\alpha$, we have 
\begin{equation}
\label{w4}
\|\varrho(x) \varrho \sigma(a) - \varrho(x \sigma(a))\| \le \frac{\bar{\delta}}{\|\tilde{h}\|}\|x\| 
\end{equation}
for all $x \in F$. Using Proposition~\ref{I} and making $\alpha$ smaller if necessary, we may even assume that
\begin{equation}
\label{wH2}
\|[b,\bar{h}(c)]\|, \, \|[b,\tilde{h}(c)]\| \le \bar{\delta}
\end{equation}
whenever $b,c \in A$ are elements of norm at most 1 with $c$ positive, and satisfying
\begin{equation}
\label{wH1}
\|[b,c]\| \le \alpha.
\end{equation}
We may further assume that 
\begin{equation}
\label{wH3}
\alpha < \bar{\delta},  \frac{\bar{\delta}}{\|\tilde{h}\|}. 
\end{equation}
Now let a c.p.c.\ approximation $(F,\sigma,\varrho)$ and $v_{i} \in A$ as in the proposition be given. Note that by \eqref{w5}, \eqref{wH1} and \eqref{wH2} we have
\begin{equation}
\label{w1}
\|[v_{i}, \bar{h}(\varrho^{(i)}(\be_{F^{(i)}}))]\| , \, \|[v_{i}, \tilde{h}(\varrho^{(i)}(\be_{F^{(i)}}))] \| \le \bar{\delta} 
\end{equation}
for $i = 0, \ldots, m$ and that
\begin{equation}
\label{w3}
\| \varrho^{(i)} (\be_{F^{(i)}}) \varrho \sigma (a) - \varrho^{(i)} \sigma^{(i)}(a) \| \stackrel{\eqref{w4}}{\le} \frac{\bar{\delta}}{\|\tilde{h}\|}
\end{equation}
for $i = 0, \ldots, m$ and $a \in \mathcal{F}$ (the elements of $\mathcal{F}$ are normalized). We obtain
\begin{eqnarray*}
\lefteqn{\|[v_{i} \bar{h}(\varrho^{(i)}(\be_{F^{(i)}})), a]\|} \\
& \le &  \|[v_{i} \bar{h}(\varrho^{(i)}(\be_{F^{(i)}})), \varrho \sigma(a)]\| +2 \alpha \\
& \stackrel{\eqref{w1}}{\le} & \| v_{i} \bar{h} (\varrho^{(i)}(\be_{F^{(i)}})) \varrho \sigma(a) - \varrho \sigma(a) \bar{h} (\varrho^{(i)}(\be_{F^{(i)}})) v_{i} \| + 2 \alpha +\bar{\delta} \\
& \stackrel{\eqref{w2}}{\le} & \| v_{i} \tilde{h} (\varrho^{(i)}(\be_{F^{(i)}})) \varrho^{(i)}(\be_{F^{(i)}})  \varrho \sigma(a) - \varrho \sigma(a)  \varrho^{(i)}(\be_{F^{(i)}}) \tilde{h} (\varrho^{(i)}(\be_{F^{(i)}})) v_{i} \|\\
&&  + 2 \alpha +\bar{\delta} + 2 \bar{\delta}\\
& \stackrel{\eqref{w3}}{\le} & \| v_{i} \tilde{h} (\varrho^{(i)}(\be_{F^{(i)}}))   \varrho^{(i)} \sigma^{(i)}(a) - \varrho^{(i)} \sigma^{(i)}(a)   \tilde{h} (\varrho^{(i)}(\be_{F^{(i)}})) v_{i} \| \\
&& + 2 \alpha +\bar{\delta} + 2 \bar{\delta} + 2 \|\tilde{h}\| \cdot \frac{\bar{\delta}}{\|\tilde{h}\|}\\
& \stackrel{\eqref{w1},\eqref{w5}}{\le} & 2 \alpha + \bar{\delta} + 2 \bar{\delta} + 2 \bar{\delta} + \bar{\delta} + \|\tilde{h}\| \cdot \alpha \\
& \stackrel{\eqref{wH3}}{<} & 9 \bar{\delta} \\
& \stackrel{\eqref{wH4}}{<} & \frac{\delta}{m+1}
\end{eqnarray*}
for $i = 0, \ldots, m$ and $ a \in \mathcal{F}$. Here, we have tacitly used that the $\varrho^{(i)}$ have order zero, whence $\varrho^{(i)}(\be_{F^{(i)}})$ commutes with $\varrho^{(i)}(F^{(i)})$. It follows that
\[
\|[v,a]\| \stackrel{\eqref{wH5}}{=} \left\|\left[\sum_{i=0}^{m} v_{i} \bar{h}(\varrho^{(i)}(\be_{F^{(i)}})),a\right]\right\| < (m+1) \frac{\delta}{m+1} = \delta
\]
for $a \in \mathcal{F}$.
\end{nproof}
\en

\bn
\label{D}
\begin{nprop}
Let $A$ be a separable, simple, unital $\mathrm{C}^{*}$-algebra with $\dr A = m < \infty$. Given a finite subset $\Fh \subset A$, $0<\delta$, $0<\zeta<1$ and $n \in \N$, there are $\Gh \subset A$ finite and $\alpha>0$ such that, whenever $(F,\sigma,\varrho)$ is an $m$-decomposable c.p.c.\ approximation of $\Gh$ to within  $\alpha$, there is $\gamma >0$ such that the following holds:

If 
\[
\Phi: M_{m+1} \otimes M_{n} \otimes M_{2} \to A
\]
is a c.p.c.\ order zero map satisfying 
\begin{equation}
\label{wwD7}
\|[\varrho(x), \Phi(y)]\| \le \gamma \|x\| \|y\| \mbox{ for all }  x \in F, \,  y \in M_{m+1} \otimes M_{n} \otimes M_{2}
\end{equation}
and 
\begin{equation}
\label{wwD12}
\tau( \Phi(\be_{M_{m+1} \otimes M_{n} \otimes M_{2}})) > 1 - \gamma \mbox{ for all } \tau \in T(A), 
\end{equation}
then  
\[
\|[\Phi(y),a]\| \le \delta \|y\| \mbox{ for all }  y \in M_{m+1} \otimes M_{n} \otimes M_{2}, \, a \in \Fh,
\]
and there is $v \in A$ such that
\[
\|v^{*} v - (\be_{A} - \Phi(\be_{M_{m+1} \otimes M_{n} \otimes M_{2}}))\| < \delta,
\]
\[
vv^{*} \in \overline{(\Phi(\be_{M_{m+1}} \otimes e_{11} \otimes \be_{M_{2}})-\zeta)_{+} A (\Phi(\be_{M_{m+1}} \otimes e_{11} \otimes \be_{M_{2}})-\zeta)_{+}}
\]
and
\[
\|[a,v]\| < \delta  \mbox{ for all } a \in \Fh.
\]
\end{nprop}

\begin{nproof}
We may clearly assume that the elements of $\Fh$ are normalized. Set 
\begin{equation}
\label{wwD9}
\bar{h}:= g_{0, \delta/(6(m+1))}
\end{equation}
From Proposition~\ref{H}, obtain a finite subset $\Gh \subset A$ and $\alpha>0$ such that the assertion of \ref{H} holds. We may assume that 
\begin{equation}
\label{wwD10}
\alpha < \frac{\delta}{6}
\end{equation} 
and that 
\begin{equation}
\label{wwD6}
\be_{A} \in \mathcal{G}.
\end{equation}  

Now suppose 
\begin{equation}
\label{wwD11}
(F,\sigma,\varrho)
\end{equation} 
is a c.p.c.\ approximation of $\Gh$ to within  $\alpha$ such that $\varrho$ is $m$-decomposable with respect to $F=F^{(0)} \oplus \ldots \oplus F^{(m)}$.

Choose some 
\begin{equation}
\label{2}
0< \theta<\frac{\alpha}{m+1}. 
\end{equation}
Obtain $\beta >0$ from Proposition~\ref{E} so that the assertion of \ref{E} holds for each 
\[
\varrho^{(i)}:= \varrho|_{F^{(i)}},
\]
$i=0,\ldots,m$, in place of $\psi$. Using Proposition~\ref{I}, we may choose some 
\begin{equation}
\label{wwD1}
0< \gamma < \min\{\theta,\beta\}
\end{equation}
such that, if $b$ and $c$ are elements of norm at most one in some $\mathrm{C}^{*}$-algebra, with $c \ge 0$ and 
\[
\|[b,c]\|< \gamma,
\]
then 
\begin{equation}
\label{wwD8}
\|[b,\bar{h}(c)]\| < \theta.
\end{equation} 

Now if $\Phi$ is as in the assertion of Proposition~\ref{D}, then (note \eqref{wwD1}) Proposition~\ref{E} yields normalized elements $v_{i} \in A$, $i=0,\ldots,m$, such that
\begin{eqnarray}
\|v_{i}^{*} v_{i} - (\be_{A} - \Phi(\be_{M_{m+1} \otimes M_{n} \otimes M_{2}}))^{\frac{1}{2}} \varrho^{(i)}(\be_{F^{(i)}}) (\be_{A} - \Phi(\be_{M_{m+1} \otimes M_{n} \otimes M_{2}}))^{\frac{1}{2}} \| && \nonumber \\
<  \theta, && \label{wwD5}
\end{eqnarray}
\begin{equation}
\label{wwD2}
v_{i}v_{i}^{*} \in \overline{(\Phi(d_{ii} \otimes e_{11} \otimes \be_{M_{2}}) - \zeta)_{+} A (\Phi(d_{ii} \otimes e_{11} \otimes \be_{M_{2}}) - \zeta)_{+}}
\end{equation}
and
\begin{equation}
\label{1}
\|[\varrho^{(i)}(x),v_{i}]\| \le \theta \|x\| \mbox{ for all }  x \in F^{(i)}
\end{equation}
for each $i$.  Note that 
\begin{equation}
\label{wwD3}
v_{i} v_{i}^{*} \perp v_{i'} v_{i'}^{*}
\end{equation}
by \eqref{wwD2} if $i \neq i'$, since $\Phi$ has order zero.   Set
\begin{equation}
\label{wwD4}
v:= \sum_{i=0}^{m} v_{i} \bar{h}(g^{(i)}(\be_{F})),
\end{equation}
then
\begin{eqnarray*}
\lefteqn{\|v^{*}v - (\be_{A}- \Phi(\be_{M_{m+1} \otimes M_{n} \otimes M_{2}})) \|}\\
& \stackrel{\eqref{wwD4},\eqref{wwD3}}{=} & \left\| \sum_{i=0}^{m} \bar{h}(\varrho^{(i)}(\be_{F^{(i)}})) v_{i}^{*} v_{i}  \bar{h}(\varrho^{(i)}(\be_{F^{(i)}})) -  (\be_{A}- \Phi(\be_{M_{m+1} \otimes M_{n} \otimes M_{2}}))  \right\| \\
& \stackrel{\eqref{wwD5},\eqref{wwD6}}{\le} & \left\| \sum_{i=0}^{m} \bar{h}(\varrho^{(i)}(\be_{F^{(i)}}))  (\be_{A}- \Phi(\be_{M_{m+1} \otimes M_{n} \otimes M_{2}}))^{\frac{1}{2}} \varrho^{(i)}(\be_{F^{(i)}}) \right. \\
&& \cdot (\be_{A}- \Phi(\be_{M_{m+1} \otimes M_{n} \otimes M_{2}}))^{\frac{1}{2}} \bar{h}(\varrho^{(i)}(\be_{F^{(i)}})) \\
&& \left. -  (\be_{A}- \Phi(\be_{M_{m+1} \otimes M_{n} \otimes M_{2}}))^{\frac{1}{2}} \varrho^{(i)}(\be_{F^{(i)}})   (\be_{A}- \Phi(\be_{M_{m+1} \otimes M_{n} \otimes M_{2}}))^{\frac{1}{2}}  \right\| \\
&&  + (m+1) \theta + \alpha \\
& \stackrel{\eqref{wwD7},\eqref{wwD8}}{\le} & \left\| \sum_{i=0}^{m} (\be_{A}- \Phi(\be_{M_{m+1} \otimes M_{n} \otimes M_{2}}))^{\frac{1}{2}}   (\bar{h}(\varrho^{(i)}(\be_{F^{(i)}}))   \varrho^{(i)}(\be_{F^{(i)}})    \bar{h}(\varrho^{(i)}(\be_{F^{(i)}})) \right. \\
&& \left. -    \varrho^{(i)}(\be_{F^{(i)}})   )      (\be_{A}- \Phi(\be_{M_{m+1} \otimes M_{n} \otimes M_{2}}))^{\frac{1}{2}} \right\| \\
&& + (m+1) \theta + \alpha + 2 (m+1) \theta \\
&\stackrel{\eqref{wwD9}}{\le} & + (m+1) \theta + \alpha + 2 (m+1) \theta + 2(m+1) \frac{\delta}{6(m+1)} \\
& \stackrel{\eqref{2},\eqref{wwD10}}{<} & \delta.
\end{eqnarray*}
Moreover,
\begin{eqnarray*}
vv^{*} & \stackrel{\eqref{wwD4},\eqref{wwD2}}{\in} & \her \left(\sum_{i=0}^{m} (\Phi(d_{ii} \otimes e_{11} \otimes \be_{M_{2}})-\zeta)_{+}\right)  \\
&\stackrel{\ref{order-zero-facts}}{=}& \her ((\Phi(\be_{M_{m+1}} \otimes e_{11} \otimes \be_{M_{2}})- \zeta)_{+}).
\end{eqnarray*}
We have 
\[
\|[v,a]\|  < \delta \mbox{ for all } a \in \Fh
\]
by Proposition~\ref{H} and our choice of $\alpha$ and $g$, using \eqref{1} and \eqref{2}. Finally, for $y \in M_{m+1} \otimes M_{n} \otimes M_{2}$ and $a \in \Fh$ we have
\begin{eqnarray*}
\|[\Phi(y),a]\| & \stackrel{\eqref{wwD11}}{\le}  & \|[\Phi(y), \varrho \sigma (a)] \| + 2 \alpha \|y\| \\
& \stackrel{\eqref{wwD7}}{<} & \gamma \|y\| + 2 \alpha \|y\| \\
& \stackrel{\eqref{wwD1},\eqref{2},\eqref{wwD10}}{<} & \delta \|y\|.
\end{eqnarray*}
\end{nproof}
\en

\section{The main result and its consequences}
\label{main-result}

\noindent
We are finally prepared to assemble the technical results of the preceding sections to prove the main result; we also derive a number of corollaries and explain some applications.

\bn
\label{A}
\begin{ntheorem}
Let $A$ be a separable, simple, nonelementary, unital $\mathrm{C}^{*}$-algebra with finite decomposition rank. 

Then, $A$ is $\Zh$-stable.
\end{ntheorem}

\begin{nproof}
Let $m:= \dr A$. We check that $A$ satisfies the hypotheses of Proposition~\ref{B}. So, let $n \in \N$, $\Fh \subset A$ finite and $\eta>0$ be given.

We have to find $\varphi:M_{n} \to A$ and $v \in A$ as in the hypotheses of \ref{B}. But by Proposition~\ref{C}, there are $0<\delta<1$ and $0<\zeta<1$ such that, if there are a c.p.c.\ order zero map
\[
\varphi': M_{n} \to A
\]
and $v' \in A$ satisfying the hypotheses in \ref{C}, then there are $\varphi$ and $v$ as desired. 

By Proposition~\ref{D}, there are a finite subset $\Gh \subset A$ and $\alpha>0$, such that, if  $(F, \sigma, \varrho)$ is an $m$-decomposable c.p.c.\ approximation of $\Gh$ to within  $\alpha$, there is $\gamma >0$ such that the following holds: If
\[
\Phi: M_{m+1} \otimes M_{n} \otimes M_{2} \to A
\]
is a c.p.c.\ order zero map satisfying \eqref{wwD7} and \eqref{wwD12} of \ref{D}, then there is $v' \in A$ which, together with 
\[
\varphi':= \Phi|_{\be_{M_{m+1}}\otimes M_{n} \otimes \be_{M_{2}}},
\]
satisfies the hypotheses of \ref{C}. 

So, choose an $m$-decomposable c.p.c.\  approximation $(F, \sigma, \varrho)$ for $\Gh$ to within  $\alpha$.  
The existence of $\Phi$ now follows from Lemma~\ref{F} with $(m+1)\cdot n \cdot 2$ in place of $k$ and $\Eh:= \varrho(\Bh_{1}(F))$ (where $\Bh_{1}(F)$ denotes the unit ball of $F$). 
\end{nproof}
\en

\bn
Combining our result with the classification theorem of \cite{Win:localizingEC}, we now see that finite decomposition rank entails classification, at least in the presence of the UCT and if there are enough projections to distinguish traces.  

\label{few-traces-classification}
\begin{ncor}
The class of simple, separable, nonelementary, unital $\mathrm{C}^{*}$-algebras with finite decomposition rank, which satisfy the UCT and for which projections separate tracial states satisfy the Elliott conjecture. 
\end{ncor}

\begin{nproof}
By \cite[Corollary~8.1]{Win:localizingEC} in connection with \cite{Lin:localizingECappendix} and \cite{LinNiu:KKlifting} (to remove the remaining $\mathrm{K}$-theory condition of \cite{Win:localizingEC}), the class of the corollary satisfies the Elliott conjecture up to $\Zh$-stability; but the latter is automatic by Theorem~\ref{A}.
\end{nproof}

The result considerably generalizes \cite[Corollary~6.5]{Winter:fintopdim} and \cite[Corollary~5.1]{Win:Z-class}. Note that it covers the real rank zero case as well as the monotracial projectionless case.
\en

\bn
As an application, we can now complete the classification of $\mathrm{C}^{*}$-algebras coming from uniquely ergodic, smooth, minimal dynamical systems.

\begin{ncor}
The $\mathrm{C}^{*}$-algebras associated to uniquely ergodic, smooth, minimal dynamical systems (with compact, finite dimensional, smooth manifolds as underlying spaces) are classified by their ordered $\mathrm{K}$-groups.
\end{ncor}

\begin{nproof}
By \cite[Corollary~8.4]{Win:localizingEC}, we have classification up to $\Zh$-stability. By \cite[Corollary~1.7]{Winter:subhomdr} (using \cite{LinPhi:mindifflimits}), the $\mathrm{C}^{*}$-algebras of the corollary have finite decomposition rank, hence are $\Zh$-stable by Theorem~\ref{A}. The invariant reduces to just ordered $\mathrm{K}$-groups in this case, since we have only one tracial state by unique ergodicity.
\end{nproof}
\en

\bn
\begin{nexamples}
(i) Theorem~\ref{A} shows that the examples of \cite{Vil:perforation}, \cite{Toms:classproblem} and \cite{Toms:example} all have infinite decomposition rank, since they are not $\mathcal{Z}$-stable. (Before, this was only known for those examples with small tracial state space, like the ones in \cite{Vil:sr=n}.)

(ii) Together with \cite{LinPhi:mindifflimits}, Theorem~\ref{A} shows that crossed products of the form $\Ch(S^{3}) \rtimes_{\alpha} \Z$ considered in \cite[Section~5, Example~4]{Con:Thom} are $\mathcal{Z}$-stable.

(iii)  Among many other examples, Corollary~\ref{few-traces-classification} in particular covers  UHF algebras, Bunce--Deddens algebras, irrational rotation algebras, the above-mentioned crossed products of odd spheres by minimal diffeomorphisms, and the Jiang--Su algebra itself.
\end{nexamples}
\en

\bn
We also obtain the following generalization of \cite[Theorem~7.6]{RorWin:Z-revisited}, which may be regarded as a finite version of Kirchberg's characterization of $\Oh_{\infty}$ as the uniquely determined purely infinite $\mathrm{C}^{*}$-algebra which is $\mathrm{KK}$-equivalent to $\C$; cf.\ \cite{Kir:fields}.
  
\begin{ncor}
The Jiang--Su algebra $\Zh$ is the uniquely determined separable, simple, nonelementary, unital $\mathrm{C}^{*}$-algebra with finite decomposition rank and unique tracial state, which satisfies the UCT and which is $\mathrm{KK}$-equivalent to the complex numbers.
\end{ncor}

\begin{nproof}
The Jiang--Su algebra is well known to satisfy the characterizing properties, so the statement follows from Corollary~\ref{few-traces-classification}.
\end{nproof}
\en

\bibliographystyle{amsplain}

\providecommand{\bysame}{\leavevmode\hbox to3em{\hrulefill}\thinspace}
\providecommand{\MR}{\relax\ifhmode\unskip\space\fi MR }
\providecommand{\MRhref}[2]{%
  \href{http://www.ams.org/mathscinet-getitem?mr=#1}{#2}
}
\providecommand{\href}[2]{#2}

\end{document}